\documentclass[10pt]{article}
    \usepackage[utf8]{inputenc}
    \usepackage{xurl}
    \usepackage{microtype}
    \usepackage{amsfonts}
    \usepackage{amsmath}
    \usepackage{amsthm}
    \usepackage{thmtools}
    \usepackage{amssymb}
    \usepackage{mathtools}
    \usepackage{subcaption}
    \usepackage[title]{appendix}
        \AtBeginEnvironment{appendices}{\crefalias{section}{appendix}}
    \PassOptionsToPackage{hyphens}{url}\usepackage{hyperref}
    \usepackage[dvipsnames]{xcolor}
    \usepackage[noabbrev,nameinlink,capitalise]{cleveref}

    \usepackage[backend=biber,style=alphabetic]{biblatex}
    \usepackage{tikz}
        \usetikzlibrary{cd}
        \usetikzlibrary{
            calc,
            arrows,
            }
        \tikzset{every picture/.style=thick}
        \tikzset{vertex/.style = {shape=circle,draw,inner sep=1pt, outer sep=1pt}}
        \tikzset{edge/.style = {->, >={latex[scale=1.5]}}}
        
\def\forkindep{\mathrel{\raise0.2ex\hbox{\ooalign{\hidewidth$\vert$\hidewidth\cr\raise-0.9ex\hbox{$\smile$}}}}}

    \theoremstyle{plain}
        \newtheorem{theorem}{Theorem}[section]
        \newtheorem{corollary}[theorem]{Corollary}
        \newtheorem{lemma}[theorem]{Lemma}
        \newtheorem{proposition}[theorem]{Proposition}
        
        \newtheorem{introtheorem}{Theorem}
        
        \newtheorem{introcorollary}{Corollary}
        
    \theoremstyle{definition}
        \newtheorem{definition}[theorem]{Definition}
        \newtheorem{example}[theorem]{Example}
        \newtheorem{introdefinition}{Definition}
        
    \theoremstyle{remark}
        \newtheorem{remark}[theorem]{Remark}
        \newtheorem{question}[theorem]{Question}
        \newtheorem{introremark}{Remark}

    \newcommand{\Aut}{\mathrm{Aut}}
    
    \newcommand{\Sym}{\mathrm{Sym}}
    \newcommand{\Cyc}{\mathrm{Cyc}}
    \newcommand{\Alt}{\mathrm{Alt}}
    
    \newcommand{\Homeo}{\mathrm{Homeo}}
    
    \newcommand{\Split}{\mathrm{Split}}
    \newcommand{\Prob}{\mathrm{Prob}}
    \newcommand{\A}{\mathcal{A}}
    \newcommand{\B}{\mathcal{B}}
    \newcommand{\R}{\mathcal{R}}
    \newcommand{\U}{\mathcal{U}}
    \newcommand{\K}{\mathcal{K}}
    \newcommand{\Chi}{\mathcal{X}}
    \newcommand{\Br}{\mathrm{Br}}
    \newcommand{\interior}[1]{\overset{\circ}{#1}}
    \newcommand{\C}{\mathbf{C}}

\title{Homeomorphism groups of basilica, rabbit and airplane Julia sets}
\author{Bruno Duchesne \and Matteo Tarocchi}
\date{}

\hypersetup{hypertexnames=false}

\begin{document}

\maketitle

\begin{abstract}
    The airplane, the basilica and the Douady rabbit (and, more generally, rabbits with more than two ears) are well-known Julia sets of complex quadratic polynomials.
    In this paper we study the groups of all homeomorphisms of such fractals and of all automorphisms of their laminations.
    In particular, we identify them with some kaleidoscopic group or universal groups and thus realize them as Polish permutation groups.
    From these identifications, we deduce algebraic, topological and geometric properties of these groups.
\end{abstract}


\section*{Introduction}

\subsection*{Abstract definitions of classical fractals and identification of the homeomorphism groups}

Since its introduction on the stage of mathematics, fractal geometry has had plenty of fruitful interactions with group theory.
Geometric group theory often fabricates fractal spaces associated to groups, such as boundaries of hyperbolic spaces (introduced in \cite[Section 1.8]{Gromov}, see also the survey \cite{HyperbolicBoundaries}) and limit spaces of contracting self-similar groups \cite[Chapter 3]{Nekrashevych}.
In the converse direction, Julia sets (fractals arising from complex dynamics) have been used as a main ingredient for defining certain groups known as iterated monodromy groups \cite[Chapter 5]{Nekrashevych}.
Recent interest in groups acting on fractals, such as rearrangement groups \cite{Belk-Forrest} (which are discrete countable examples) or dendrite and kaleidoscopic groups \cite{Kaleidoscopic} (uncountable Polish groups), is also noteworthy.

\begin{figure}
\centering
\includegraphics[width=\textwidth]{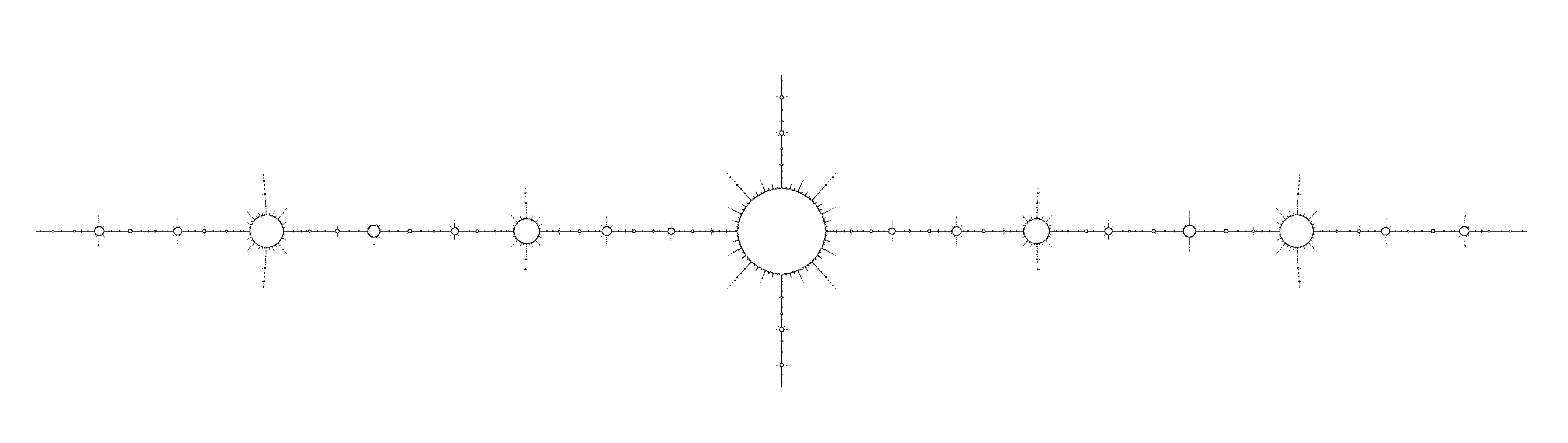}
\caption{The airplane Julia set.}
\label{fig.airplane}
\end{figure}

\begin{figure}
\centering
\includegraphics[width=.8\textwidth]{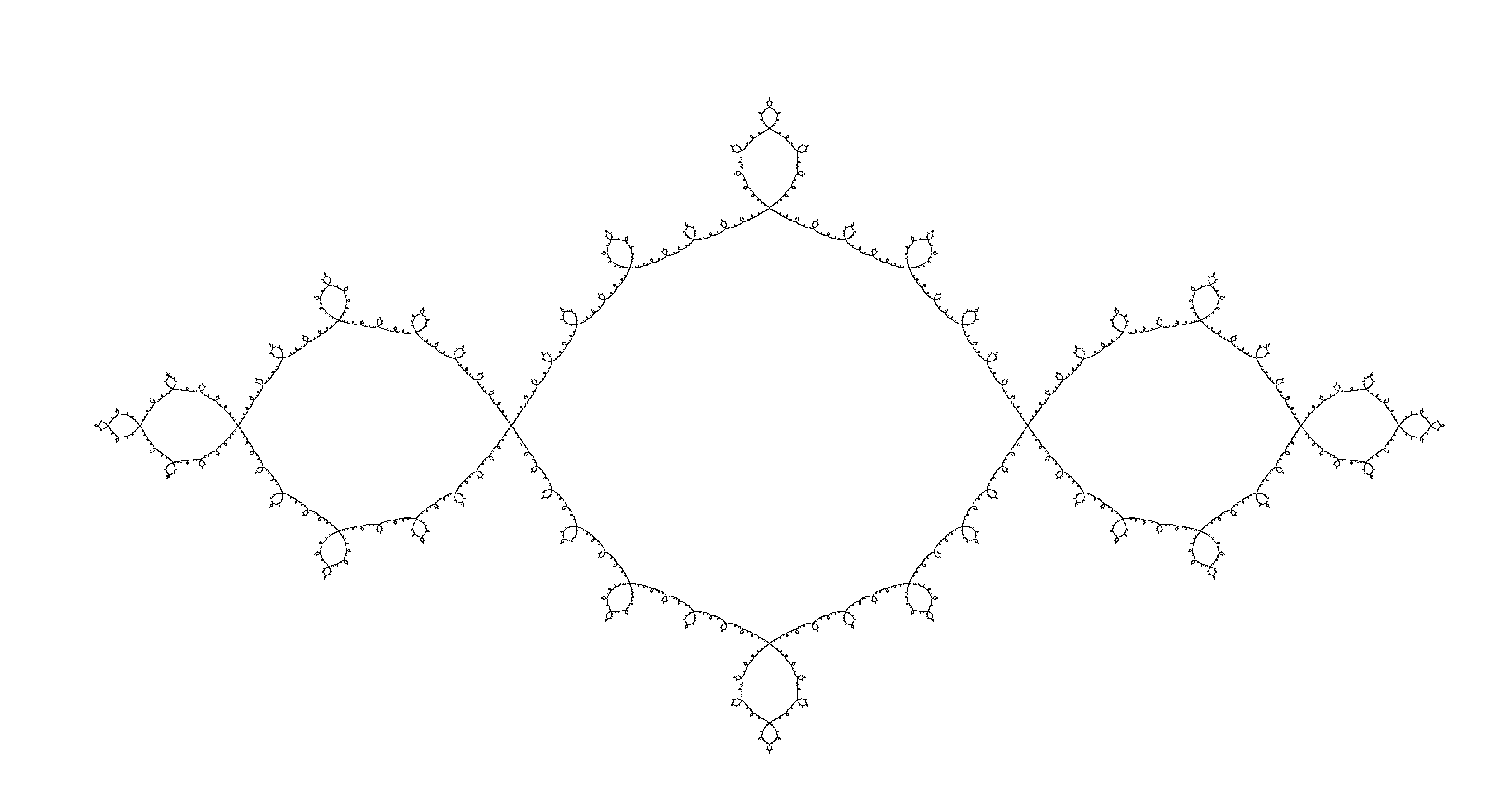}
\caption{The basilica fractal.}
\label{fig.basilica}
\end{figure}

\begin{figure}
\centering
\includegraphics[width=.7\textwidth]{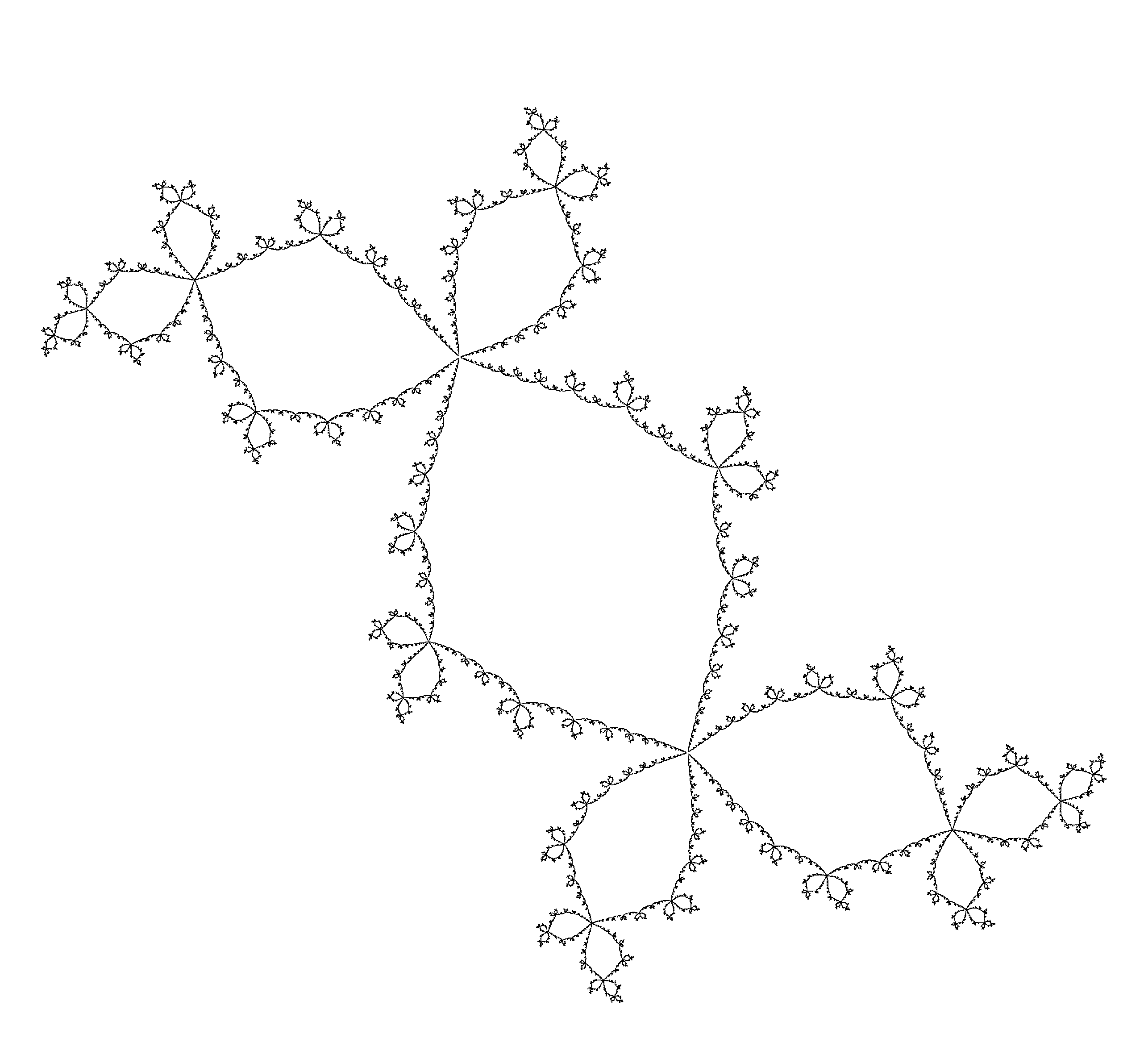}
\caption{The Douady rabbit fractal.}
\label{fig.douadyrabbit}
\end{figure}

\begin{figure}
\centering
\includegraphics[width=.6\textwidth]{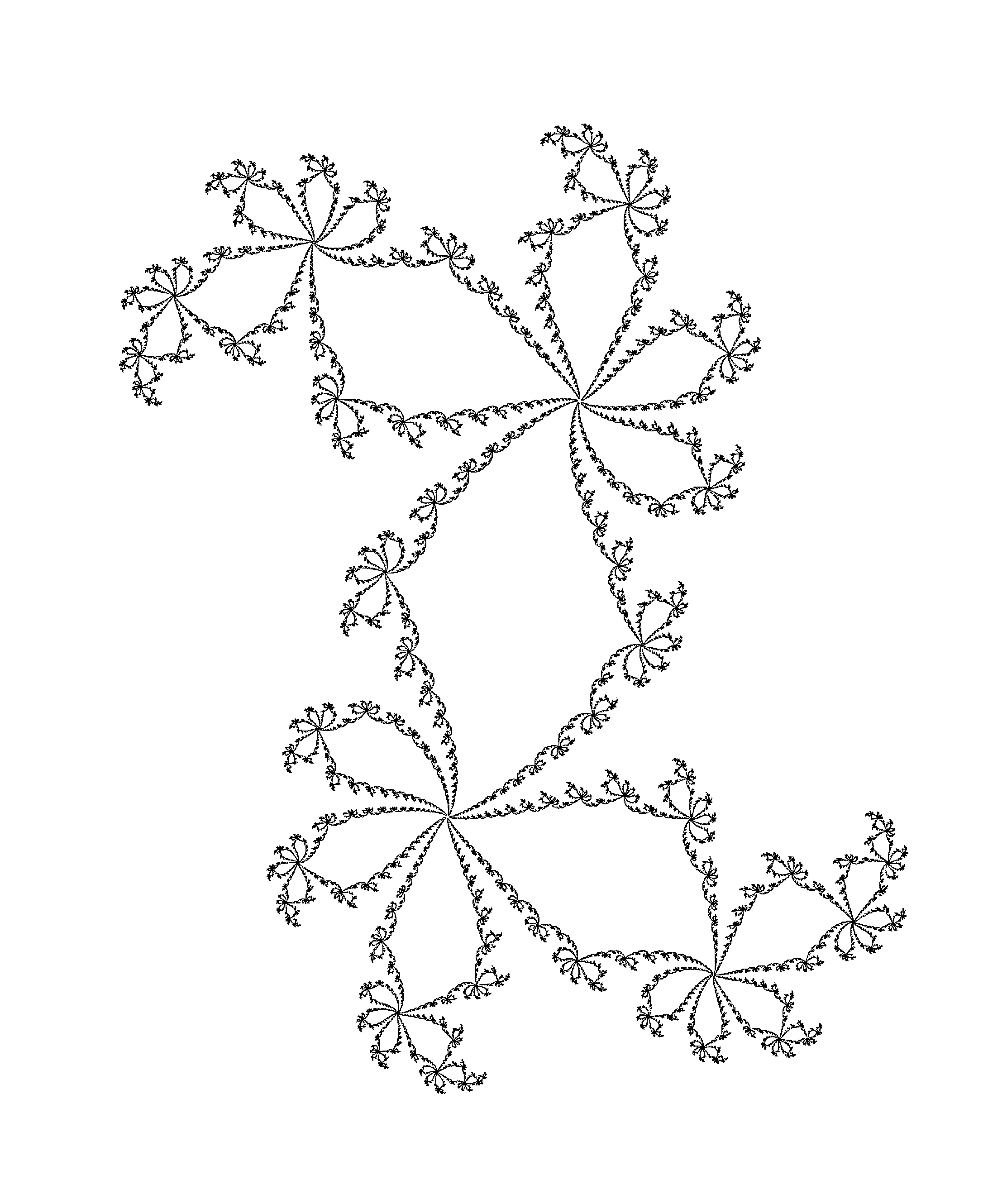}
\caption{A 5-regular rabbit.}
\label{fig.5-rabbit}
\end{figure}

Julia sets of quadratic polynomials are compact sets that exhibit numerous abstract homeomorphisms, independent of their embedding in the complex plane. However, the study of their homeomorphism groups appears to have been largely overlooked, with the notable exception of \cite{Neretin}.
In this paper, we study certain well-known fractal spaces: the airplane $\A$ (\cref{fig.airplane}) and the rabbits $\R_n$, which include the basilica $\B$ (\cref{fig.basilica}), the Douady rabbit $\R_3$ (\cref{fig.douadyrabbit}) and others (\cref{fig.5-rabbit}).
All of these are Julia sets of complex quadratic polynomials.

Our focus is on analyzing their full homeomorphism groups and certain natural positive subgroups that preserve the associated laminations.
We begin by providing topological characterizations of these fractals.
These characterizations apply to the basilica, rabbit, and airplane Julia sets, as well as to the limit spaces of specific edge replacement systems, as detailed in \cref{sec.julia.sets,sec.limit.spaces}, respectively.

Recall that a \textbf{Peano continuum} is a compact, metrizable, connected, and locally connected space.
Peano continua are path-connected, for example because of \cite[Theorem 8.23]{continua}.
By an \textbf{arc} of such a continuum, we mean a subset that is homeomorphic to the topological interval $[0,1]$ and by a \textbf{circle} we mean a subset that is homeomorphic to the topological circle $\mathbf{S}^1$.
Finally, recall that the \textbf{order} of a point is the number of connected components of its complement, and a point is a \textbf{cut point} if its order is at least $2$.
With this in mind, we present our abstract definition of rabbits.

\begin{introdefinition}
\label{def.rabbit}
A \textbf{rabbit} is a Peano continuum that enjoys the following properties.
\begin{enumerate}
    \item The set of points that belong to multiple circles is dense.
    \item If two points belong to distinct circles, they are separated by a cut point.
    \item For any two points each of which belongs to some circle, every arc between them is included in a union of finitely many circles.
\end{enumerate}
A rabbit is $n$-\textbf{regular} ($n\geq2$) if its cut points all have the same order $n \in \mathbf{N}_{\geq 2} \cup \{\infty\}$.
A \textbf{basilica} is a $2$-regular rabbit and a \textbf{Douady rabbit} is a $3$-regular rabbit.
\end{introdefinition}

We also define abstract airplanes with a very similar definition.
The main difference with rabbits is the relative position of circles:
no point belongs to more than one circle and there is a third circle between every two of them.

\begin{introdefinition}
\label{def.airplane}
An \textbf{airplane} is a Peano continuum that enjoys the following properties.
\begin{enumerate}
    \item The set of all cut points that belong to circles is dense.
    \item If two points belong to distinct circles, they are separated by a third circle.
    \item Any two distinct circles are disjoint.
    \item Every cut point has order $2$.
\end{enumerate}
\end{introdefinition}

We prove that rabbits and the airplane are unique in the following sense.

\begin{introtheorem}[\cref{thm.uniqueness.rabbit,,thm.uniqueness.airplane}]
\label{thm.uniqueness.rabbits.and.airplane}
Up to homeomorphism, there is a unique airplane and, for each $n\geq2$, there is a unique regular $n$-rabbit.
\end{introtheorem}

Thanks to this uniqueness, we denote by $\A$ the airplane and $\R_n$ the $n$-regular rabbit.
The key tool to prove \cref{thm.uniqueness.rabbits.and.airplane} is to associate another topological space that encodes the arrangements of circles: the dendrite or tree of circles. 
For the airplane, this dendrite is the universal Wa\.zewski dendrite $D_\infty$.
Let us recall that a dendrite is a Peano continuum where any two points are joined by a unique arc and can be thought as a compact metric tree.
The universal Wa\.zewski dendrite is the unique dendrite such that branch points have infinite order and are dense (it is universal in the sense that it includes homeomorphically any other dendrite, see \cite{MR3532906}).
For the regular $n$-rabbit, the tree of circles is the $(n,\infty)$-biregular tree $T_{n,\infty}$, i.e., any edge joins a vertex of degree $n$ and one of infinite degree.

Every homeomorphism of the fractal spaces that we consider induces a homeomorphism of the associated dendrite of circles in the case of the airplane or an automorphism of the tree of circles in the case of rabbits.
Moreover, the homeomorphism group of the fractal space acts faithfully on this dendrite or tree of circles.
However, since the topology of each circle in the fractal space is preserved under homeomorphisms, not every homeomorphism of the dendrite or tree of circles corresponds to a homeomorphism of the fractal space. 

The preservation of each circle's structure is captured by the separation relation:
two pairs of distinct points, \((x, y)\) and \((a, b)\), in a circle \(C\) are said to be \textbf{separated} if \(a\) and \(b\) lie in different connected components of \(C \setminus \{x, y\}\).
For a dense countable subset of a circle, we define \(\Aut(S)\) as the group of all bijections of this set that preserve the separation relation.
In fact, all such bijections induce homeomorphisms of the circle and either preserve or reverse its cyclic order.
For more details on separation relations, cyclic orders, and their automorphism groups, see \cref{sec.order.and.separation}.

In \cite{BM00}, Burger and Mozes introduced the concept of a universal group \(\U(\Gamma)\) for automorphisms of the \(n\)-regular tree \(T_n\), given a permutation group \(\Gamma \leq \Sym([n])\) (where \([n] = \{1, \dots, n\}\) if $n$ is finite and $[\infty]=\mathbf{N}$).
The universal group \(\U(\Gamma)\) is characterized by the property that the local action, i.e., the action induced on the neighborhood of each vertex, lies in \(\Gamma\).
This construction was later generalized by Smith in \cite{biregular} for biregular trees \(T_{m, n}\), where \(m, n\) are countable and possibly infinite.
In this case, two permutation groups \(\Gamma_m \leq \Sym([m])\) and \(\Gamma_n \leq \Sym([n])\) are used to construct the universal group \(\U(\Gamma_m, \Gamma_n)\).

For Ważewski dendrites, these universal groups inspired the construction of kaleidoscopic groups, introduced by the first author, Monod, and Wesolek in \cite{Kaleidoscopic}.

We prove the following identifications of the groups of homeomorphisms of rabbits and the airplane.

\begin{introtheorem}[\cref{thm.homeo.rabbits.are.auto.trees,,thm.homeo.airplane.is.kaleidoscopic}]
\label{thm.identification}
The homeomorphism group of the Airplane is the kaleidoscopic group $\K(\Aut(S))$ and the homeomorphism group of the regular $n$-rabbit is the universal group $\U(\Sym([n]),\Aut(S))$.
\end{introtheorem}

\begin{introremark}[\cref{rmk.BM.and.biregular}]
In the case of the regular $2$-rabbit, which is the basilica, since $\U(\Sym([2]),\Aut(S))$ and $\U(\Aut(S))$ are isomorphic groups, this result was already obtained by Neretin in \cite{Neretin}.
Moreover, he gave some hint about how his result should be modified for the airplane.
\end{introremark}

We emphasize that in \cref{thm.identification}, the identification is not only as abstract groups but also as Polish groups.
Recall that a \textbf{Polish group} is a separable, completely metrizable topological group.
For example, the compact-open topology is a Polish topology on the homeomorphism group of a compact metrizable space, while the permutation topology (also known as the topology of pointwise convergence) equips the group of all bijections of a countable set (such as the branch points of the Ważewski dendrite or the vertices of a biregular tree) with a Polish topology.
In particular, the homeomorphism groups of the airplane or an \(n\)-regular rabbit are non-Archimedean groups.

Although, in general, we do not consider fractal spaces as embedded in the plane, it is natural — when realizing them as Julia sets of quadratic polynomials — to endow them with an orientation arising from this embedding.
More precisely, we can equip each circle and each set of components in the complement of a branch point with a cyclic order and study the subgroups that preserve this orientation.
For a rigorous definition, see \cref{sec.lamination.automorphisms}, where we define the kaleidoscopic group \(\K(\Aut(O))\) and universal groups \(\U(\Cyc(n), \Aut(O))\) for finite $n$.
Here, \(\Aut(O)\) denotes the automorphism group of the dense countable cyclic order, and \(\Cyc(n)\) is the cyclic group of finite size \(n\) within \(\Sym([n])\).
We denote these groups by \(\Homeo^+(\A)\) in the case of the airplane and \(\Homeo^+(\R_n)\) in the case of an \(n\)-regular rabbit.
Throughout this paper, this last notation will be used only in the case where $n$ is finite, contrarily to \(\Homeo(\R_n)\) that implicitly includes the case $n=\infty$.

\begin{figure}
    \centering
    \includegraphics[width=.45\textwidth]{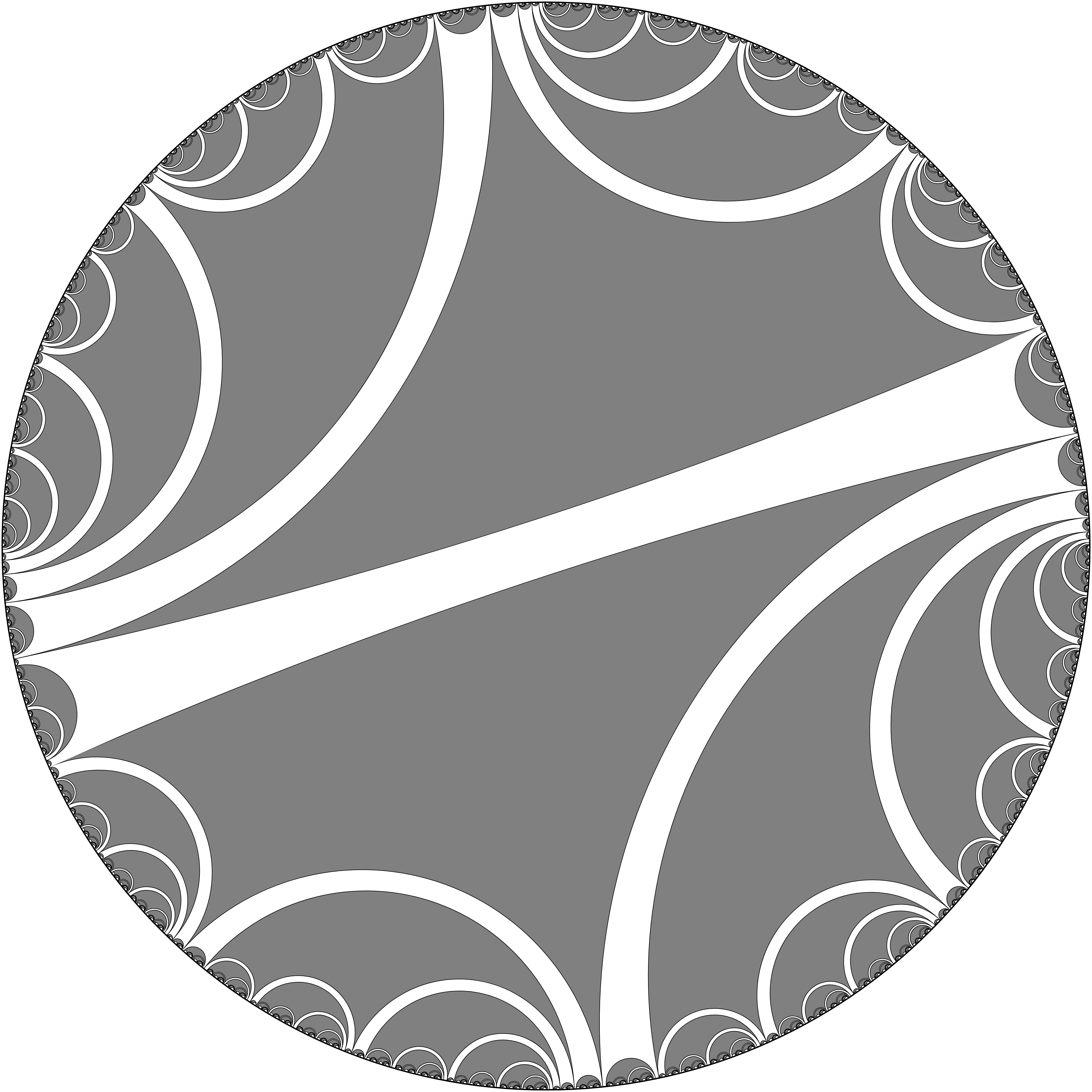}
    \hfill
    \includegraphics[width=.45\textwidth]{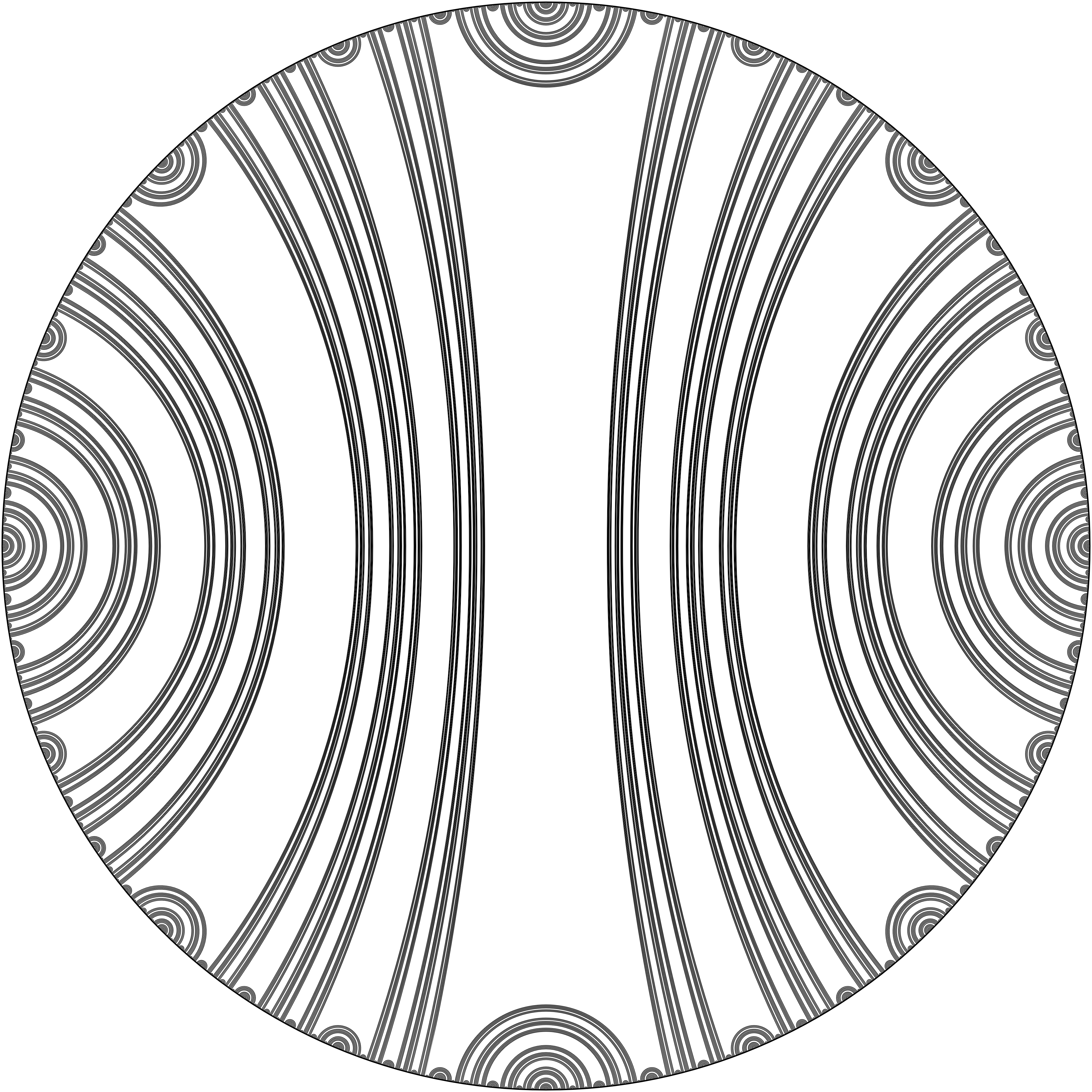}
    \caption{The $4$-rabbit lamination and the airplane lamination.}
    \label{fig.laminations}
\end{figure}

Following Thurston, it has become standard to associate a lamination of the unit disk with any Julia set (see \cref{fig.laminations} for examples).
When realizing \(\A\) or \(\R_n\) (for $n$ finite) as the Julia set of a quadratic polynomial, we naturally obtain an associated lamination.

\begin{introtheorem}[\cref{thm.auto.laminations,,cor.lamination.groups.closed.in.circle.group}]
\label{thm.positive.groups.lamination}
The group $\Homeo^+(\A)$ and the groups $\Homeo^+(\R_n)$ are the subgroups of $\Homeo^+(\mathbf{S}^1)$ that preserve the associated laminations, which are closed in $\Homeo^+(\mathbf{S}^1)$.
\end{introtheorem}

Moreover, we prove in \cref{sub.extending.homeomorphisms} that the groups $\Homeo^+(\A)$ and $\Homeo^+(\R_n)$ correspond to the groups of homeomorphisms of the Julia sets that can be extended to homeomorphisms of the entire complex plane.

\subsection*{Consequences of the identification of the groups}

Thanks to our identification, we can deduce algebraic, topological, dynamical and geometric properties of the groups. Let us start with simplicity results.

\begin{introtheorem}[\cref{cor.rabbits.simple,,cor.airplane.simple,,thm.simple.commutators}]
\label{thm.simplicity.results} The groups $\Homeo(\A)$, $\Homeo^+(\A)$ and $\Homeo(\R_n)$ for $n\geq3$ are simple.
The groups $\Homeo(\R_n)$ for $n=2$ and $\Homeo^+(\R_n)$ for $n\geq2$ have index-$n$ simple commutator subgroup.
\end{introtheorem}

Given a family of groups, it is natural to investigate whether one group can be embedded into another and whether they are isomorphic.
The following theorem gathers results addressing these questions, beginning with a general statement about universal groups and kaleidoscopic groups associated with the same permutation group.
For instance, an immediate consequence is that the homeomorphism group of the basilica can be embedded into the homeomorphism group of the airplane.

\begin{introtheorem}[\cref{sec.embedding}]
\label{thm.embeddings.results}
\begin{enumerate}
    \item[]
    \item For any permutation group $\Gamma$, the universal group $\mathcal{U}(\Gamma)$ embeds topologically in the kaleidoscopic group $\K(\Gamma)$.
    \item For every $2 \leq n < m$, there are topological group embeddings
    \[ \Homeo(\R_n) \hookrightarrow \Homeo(\R_m). \]
    \item The group $\Homeo(\R_m)$ does not embed into $\Homeo(\R_n)$ for $4 \leq n < m$.
    \item The abstract groups $\Homeo^+(\R_n)$ and $\Homeo^+(\R_m)$ are isomorphic if and only if $n=m$.
    \item As soon as $n | m$, there are topological group embeddings
    \[ \Homeo^+(\R_n) \hookrightarrow \Homeo^+(\R_m) \hookrightarrow \Homeo^+(\B). \]
\end{enumerate}
\end{introtheorem}

In particular, the first point implies that
$$\Homeo(\B) \hookrightarrow \Homeo(\A) \text{ and } \Homeo^+(\B) \hookleftarrow \Homeo^+(\A).$$

Up to this point, the results have been largely similar for the airplane and the regular rabbits.
The following results highlight the differences arising from the distinct nature of Ważewski dendrites and biregular trees.

Recently, Rosendal extended geometric group theory beyond discrete (and more generally, locally compact) groups to more general contexts, specifically Polish groups \cite{MR4327092}.
By introducing the left coarse structure for topological groups, Rosendal replaces relatively compact subsets with coarsely bounded subsets — those with bounded orbits under any continuous isometric action.
For a Polish group generated by a coarsely bounded subset, the group acquires a canonical quasi-metric structure, allowing meaningful discussion of the group’s quasi-isometry type (see \cite[\S 2.8 and 2.9]{MR4327092}).

Every topological group $G$ is equipped with two uniform structures: the left and the right uniformities.
The \textbf{Roelcke uniformity} is the meet of these two structures.
A subset of $G$ is termed \textbf{Roelcke precompact} if its completion under the Roelcke uniform structure is compact.
By definition, any Roelcke precompact subset is also coarsely bounded \cite[\S 3.1]{MR4327092}.  

\begin{introtheorem}[\cref{prop.rabbit.QItype,,cor.airplane.QItype}]
The topological groups $\Homeo(\R_n)$ and $\Homeo^+(\R_n)$ are locally Roelcke precompact, in particular locally bounded, and quasi-isometric to the infinite-order tree $T_\infty$.
On the other hand, the topological groups $\Homeo(\A)$ and $\Homeo^+(\A)$ are Roelcke precompact and thus coarsely bounded.
\end{introtheorem}

Another difference between the airplane and the regular rabbits comes from their unitary representation theory.
Recall that the Haagerup property and property (T) are opposite properties (see \cref{sec.haagerup.property.(T)} for definitions).

\begin{introcorollary}[\cref{thm.property.T}]
The Polish groups $\Homeo(\R_n)$ and $\Homeo^+(\R_n)$ have the Haagerup property.
On the other hand, $\Homeo(\A)$ and $\Homeo^+(\A)$ have property (T).
\end{introcorollary}

For a topological group $G$, a $G$-\textbf{flow} is a compact Hausdorff space equipped with a continuous action of $G$.
Among minimal $G$-flows (those without proper non-empty invariant subspaces), there exists a universal one: any other minimal $G$-flow is an equivariant factor of this \textbf{universal minimal flow}.
Understanding this flow provides significant insight into the topological dynamics of the group.
However, this space is often very large; for instance, for locally compact non-compact groups, the universal minimal flow is always non-metrizable.
Actually, one can find a copy of the Stone-Cech compactification $\beta X$ of an infinite discrete $X$ in such minimal flow by the discussion following \cite[Theorem 3.1.1]{Pestov_99}.

In the case of \(\Homeo(\A)\), this group can be identified with a kaleidoscopic group where the local action is oligomorphic.
As a result, the identification of the universal minimal flow of \(\Homeo(\A)\) follows directly from \cite[Theorem 1.1]{Kaleidoscopic2}.
For \(\Homeo^+(\A)\), however, the universal minimal flow has a particularly explicit description.
Following \cite{Glasner-Megrelishvili-2018}, for a subset \(A \subseteq \mathbf{S}^1\), let \(\Split(\mathbf{S}^1, A)\) denote the compact space obtained from \(\mathbf{S}^1\) (with its continuous cyclic order) by replacing each point \(a \in A\) with a pair of points \(a^-, a^+\), where \(a^+\) is the cyclic successor of \(a^-\).
Let \(P\) be the countable set that corresponds to the preimage of the set of cut points of \(\A\) under the Carathéodory loop.

\begin{introtheorem}[\cref{prop.UMF.homeo.airplane.is.metrizable,,thm.universal.minimal.flow.airplane}]
\label{thm.universal.minimal.flows}
The universal minimal flow of $\Homeo(\A)$ is metrizable and the universal minimal flow of $\Homeo^+(\A)$ is $\Split(\mathbf{S}^1,P)$. 
\end{introtheorem}

We observe in particular that $\Split(\mathbf{S}^1,P)$ has a $\Homeo^+(\A)$-invariant cyclic order, but the following sequence of equivariant factors
\[ \Split(\mathbf{S}^1,P)\to\mathbf{S}^1\to\A\to D_\infty \]
gives an example where the universal minimal flow is cyclically ordered while some other minimal flows are not, since $\A$ and $D_\infty$ have no invariant cyclic order under $\Homeo^+(\mathbf{S}^1)$ (\cref{prop.no.invariant.cyclic.order}).
This answers a question of Glasner and Megrelishvili, \cite[Question 4.9]{Glasner-Megrelishvili-2018}.

\subsection*{Figure attributions}

The pictures of Julia sets displayed in this article were produced using Python code developed by the first author, which is available at the following link:
\url{https://gitlab.imo.universite-paris-saclay.fr/bruno.duchesne/julia-sets-in-python}

The pictures of laminations (\cref{fig.laminations}) were instead made with software developed by Caleb Falcione, which is available at the following link:
\url{https://csfalcione.github.io/lamination-builder/}.

\subsection*{Acknowledgements}

We thank Marco Barbieri for discussions about actions of symmetric groups on trees, Gianluca Basso for asking if there is an abstract characterization of the airplane, basilica and rabbits and Francesco Matucci for conversations about homeomorphisms of fractals and universal groups.

We are grateful to Matthieu Joseph, Eli Glasner and Michael Megrelishvili and Yuri Neretin for their comments on a preliminary version of this text.

We are thankful to the anonymous referee for their comments, which led to corrections, improvements and a significant reduction of typos.

The first author is supported by the "Beach" project ANR-24-CE40-3137.
The second author is supported by the "GoFR" project ANR-22-CE40-0004 and is a member of the Gruppo Nazionale per le Strutture Algebriche, Geometriche e le loro Applicazioni (GNSAGA) of the Istituto Nazionale di Alta Matematica (INdAM).


\section{The rabbits and the basilica}
\label{sec.rabbits}

In this section we consider the family of fractals that we call \textit{rabbits}, which generalize the notorious basilica and the Douady rabbit fractals (\cref{fig.basilica,fig.douadyrabbit}).
We discuss their topological characterization, which we gave in \cref{def.rabbit} (it is shown in \cref{sec.julia.sets} that certain Julia sets satisfy this topological characterization).
We expand on the ideas of \cite{Neretin} and we prove that the homeomorphism groups of rabbits are topologically isomorphic to certain locally determined universal groups of automorphisms of biregular trees, which were introduced and studied in \cite{biregular}.

\subsection{Topological characterization of rabbits}

Throughout this work, given points $a,b \in X$ (or subsets $A,B \subseteq X$), we say that a point $p$ or a subset $S \subset X$ \textbf{separates} $a$ and $b$ if they lie in distinct connected components of $X \setminus \{p\}$ or $X \setminus S$ (or separates $A$ and $B$ if they are included in distinct connected components of $X \setminus \{p\}$ or $X \setminus S$).

Let us now focus on the family of rabbits, which we defined in \cref{def.rabbit}.

\begin{figure}
\centering
\includegraphics[width=.8\textwidth]{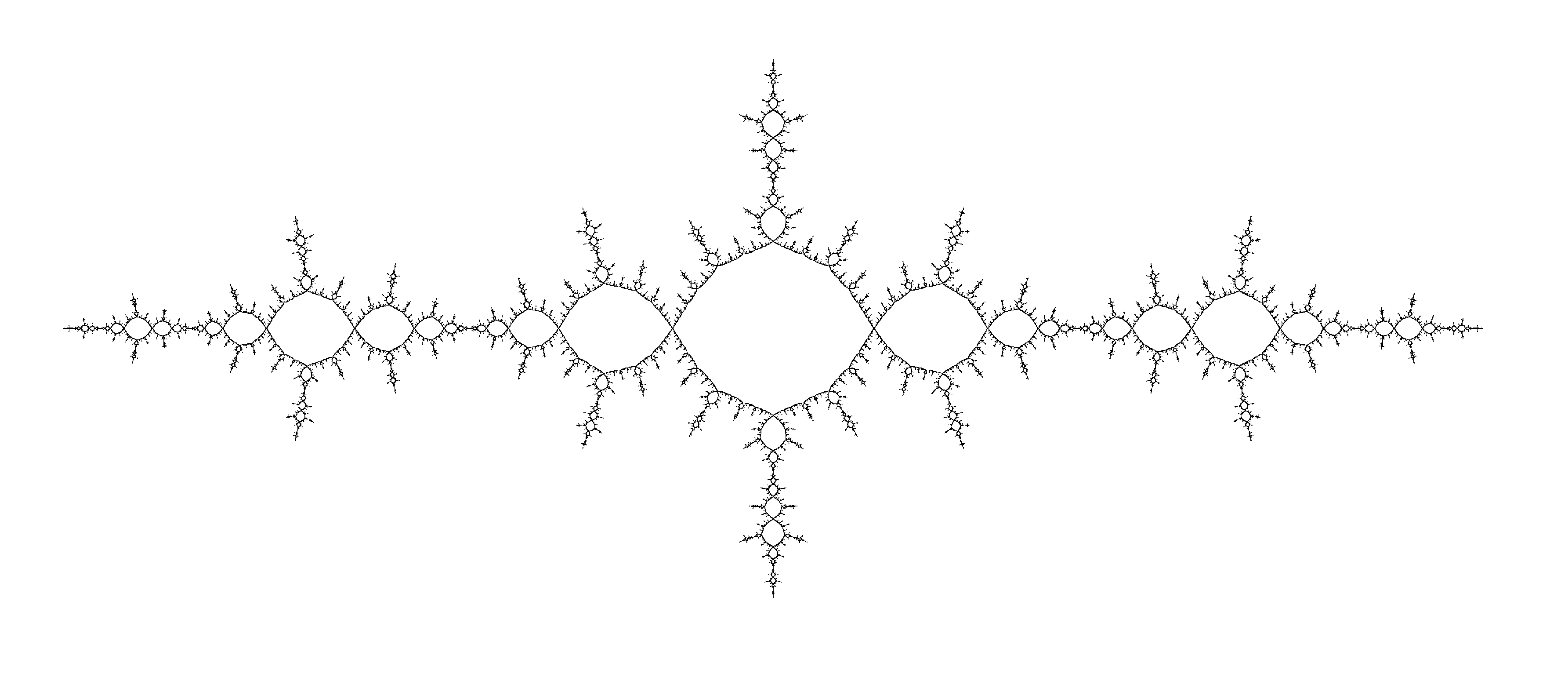}
\caption{A space that satisfies conditions 1 and 2 but not 3 of \cref{def.rabbit}.}
\label{fig.quasi.rabbit}
\end{figure}

\begin{remark}
Without condition 3, \cref{def.rabbit} would include Julia sets such as that for the map $z \mapsto z^2 + c$ where $c \approx -1.311$, which is depicted in \cref{fig.quasi.rabbit}.
\end{remark}

\begin{definition}
\label{def.rabbit.points}
A point of a rabbit is said to be:
\begin{itemize}
    \item an \textbf{end point} if it does not belong to any circle,
    \item a \textbf{cut point} if its complement is not connected,
    \item a \textbf{regular point} if it belongs to a circle and is not a cut point (equivalently thanks to \cref{lem.circles.intersection} right below, if it belongs to a unique circle).
\end{itemize}
\end{definition}

The reason behind the name of end points will be clarified with \cref{rmk.rabbit.end.points}, after the introduction of the tree of circles.

The next two lemmas show why all conditions of \cref{def.rabbit} are essential.

\begin{lemma}
\label{lem.circles.intersection}
The intersection between two distinct circles of a rabbit is either empty or consists of a sole cut point.
\end{lemma}

\begin{proof}
Consider two distinct circles $C_1$ and $C_2$ of a rabbit.
It is impossible that one includes the other, so there exist points $p \in C_1 \setminus C_2$ and $q \in C_2 \setminus C_1$.
By condition 2 of \cref{def.rabbit}, $p$ and $q$ are separated by a cut point.
Because of this, $C_1$ and $C_2$ cannot meet at more than one point, as otherwise $p$ and $q$ could not be separated by a cut point.
Then $C_1 \cap C_2$ is either empty or consists of a sole cut point.
\end{proof}

\begin{lemma}
\label{lem.cut.points.order}
The order of a cut point is equal to the number of circles that include it.
\end{lemma}

\begin{proof}
We need to show that the removal of a cut point $p$ from a rabbit $\R$ produces a distinct connected component for each circle $C$ that contains $p$. 

By condition 2 of \cref{def.rabbit} there is a distinct connected component for each circle that meets at $p$.
Indeed, if two circles $C_1$ and $C_2$ at $p$ are such that $C_1 \setminus \{p\}$ and $C_2 \setminus \{p\}$ lie in the same connected component of $\R \setminus \{p\}$, then it would be impossible to separate $C_1$ and $C_2$ with a single cut point in $\R$, which contradicts condition 2.

For the converse, if $A$ is a connected component of $\R \setminus \{p\}$, we need to show that $A \cup \{p\}$ meets $p$ in a circle.
Fix another connected component $B$ at $p$, which is going to be useful to apply condition 3 of \cref{def.rabbit}.
Since $\R$ is locally connected (by the definition of Peano continua), $A$ and $B$ are open, so there exist points $a \in A$ and $b \in B$ that belong to multiple circles by condition 1 of \cref{def.rabbit}.
Fix any arc $[a,b]$ and note that it must contain $p$, since $p$ separates $A \ni a$ and $B \ni b$.
By condition 3 of \cref{def.rabbit}, since both $a$ and $b$ belong to some circle, the arc $[a,b]$ is included in a union of finitely many circles.
In particular, its subset $[a,p)$ is included in a union of finitely many circles.
Then there is some $q \in [a,p)$ such that $[q,p)$ is included in some circle $C$.
Since $C$ is closed, we have that $p \in C$.
Finally, observe that if a circle intersects $A$ non-trivially then it must be entirely included in $A \cup \{p\}$.
Thus, $C \subseteq A$ and we are done.
\end{proof}

\begin{corollary}
\label{cor.rabbit.cut.points}
A point of a rabbit is a cut point if and only if it belongs to multiple circles.
\end{corollary}

\begin{remark}
\label{rmk.rabbit.countable}
Note that the order of a cut point of a rabbit is at most countably infinite.
Indeed, since rabbits are Peano continua they are locally connected second countable spaces, so the complement of a point is a union of disjoint open connected subsets and, by second countability, there are at most countably many such subsets.

In a similar fashion, the removal of a circle must produce at most countably many connected components, so each circle contains at most countably many cut points.
Because of condition 1 of \cref{def.rabbit}, it must contain at least countably many, so each circle has precisely countably many cut points.
\end{remark}

\begin{lemma}
\label{lem.rabbit.arcs.and.circles}
Any point of an arc minus its end points is contained in some circle.
\end{lemma}

\begin{proof}
Let $I = [a,b]$ be an arc in a rabbit $\R$ and let $p \in (a,b) \coloneqq [a,b] \setminus \{a,b\}$.
If $p$ is a cut point, then it belongs to multiple circles by \cref{cor.rabbit.cut.points}, so there is nothing to prove.
If $p$ is not a cut point, then $\R \setminus \{p\}$ is connected and thus path-connected because $\R$ is a Peano continuum.
It follows that there exists an arc $I^*$ that joins $a$ and $b$ while avoiding $p$.
Then $I \cup I^*$ must include some circle that contains $p$.
\end{proof}

\begin{lemma}
The set of cut points of a rabbit is arcwise dense, i.e., it has non-trivial intersection with any arc.
\end{lemma}

\begin{proof}
Let us show that any arc $I=[a,b]$ contains cut points.
Fix a point $p \in (a,b)$ and assume that it is not a cut point.
Then, up to taking a subinterval, we can assume that $(a,b)$ is entirely included in some circle $C$.
Let $U$ be a path-connected neighborhood of $p$ that does not contain $a$ and $b$ (which can be found because $\R$ is locally path-connected).
By condition 1 of \cref{def.rabbit}, $U$ contains some point $p'$ that belongs to some circle.
Let $I'$ be an arc in $U$ that joins $p$ and $p'$.
The first point $q$ in which $I'$ meets $C$ (coming from $p'$) must also belong to $I$, because $U$ does not contain $a$ and $b$.
The point $q$ must be a cut point, so we are done.
\end{proof}

\begin{proposition}
\label{prop.subbasis}
Let $X$ be a Hausdorff compact space and consider a point $x \in X$.
Assume that $\mathcal{O}_x$ is a collection of open subsets of $X$ containing $x$ and with the property that, for every $y \in X$, there exists $U_y \in \mathcal{O}_x$ such that $y \notin \overline{U_y}$.
Then $\mathcal{O}_x$ is a subbasis of neighborhoods at $x$.

Moreover, given such collections $\mathcal{O}_x$ for every $x \in X$, their union $\mathcal{O} = \bigcup_{x \in X} \mathcal{O}_x$ is a subbasis for the topology of $X$.
\end{proposition}

\begin{proof}
Let $V$ be an open subset of $X$ and $x \in V$.
We want to show that there exist finitely many elements of the collection $\mathcal{O}_x$ whose intersection contains $x$ and is included in $V$.
Observe that
\[ \bigcap \left\{ \overline{U_y} \mid y \in X \setminus V \right\} \cap (X \setminus V) = \emptyset, \]
otherwise there would exist $y \in X \setminus V$ contained in $\overline{U_y}$, which is a contradiction.
Since $X$ is compact and the infinite intersection in the previous equation is empty, there exist $y_1, \dots, y_n \in X \setminus V$ such that
\[ \overline{U_{y_1}} \cap \dots \cap \overline{U_{y_n}} \cap (X \setminus V) = \emptyset, \]
meaning that $\overline{U_{y_1}} \cap \dots \cap \overline{U_{y_n}} \subseteq V$, as needed.
\end{proof}

Using this proposition, it is easy to show the following useful fact.

\begin{corollary}
\label{cor.rabbit.subbasis}
Given a rabbit $\R$, the collection of all connected components of $\R$ minus finitely many of its circles is a subbasis for the topology of $\R$.
\end{corollary}

\subsection{The tree of circles of a rabbit}\label{sec.tree.circles.rabbit}

In this subsection we will see how the circles of a rabbit are arranged in a tree-like manner, which is going to be essential when studying their homeomorphism groups.

\begin{definition}
\label{def.traverse}
Given an arc $I$ in a rabbit $\R$ and a subset $S \subseteq \R$, we say that $I$ \textbf{traverses} $S$ if the intersection $I \cap S$ includes a non-degenerate arc (i.e., an arc that is not a singleton).
\end{definition}

We will often be using the previous definition in the case in which $S$ is a circle.
In that case, because of condition 2 of \cref{def.rabbit}, we can equivalently say that an arc $I$ traverses a circle $C$ if $I \cap C$ contains at least two points.

\begin{definition}
\label{def.tree.of.circles}
Given a rabbit $\R$, its \textbf{tree of circles} $T_\R$ is the (undirected) graph defined as follows.
\begin{itemize}
    \item The set of vertices consists of the circles and the cut points of $\R$.
    \item There is an edge between a circle $C$ and a cut point $p$ when $p \in C$.
\end{itemize}
\end{definition}

The fact that the graph $T_\R$ is really a tree is proved in \cref{lem.tree.of.circles.is.tree}.

\begin{remark}
\label{rmk.arcs.and.paths}
Because of \cref{lem.circles.intersection} and by condition 3 of \cref{def.rabbit}, there is a natural surjective map from the set of arcs of a rabbit $\R$ joining cut points or regular points onto the set of paths of the graph $T_\R$.
An arc $I = [p,q]$ is mapped to the sequence of circles traversed by $I$ (\cref{def.traverse}), interposed by the cut points between them.
If $p$ is a cut point, then the first vertex of the path is $p$ itself;
if $p$ is a regular point that belongs to a unique circle $C$, then the first vertex is $C$.
The same goes for $q$.
The preimages of this map are only distinguished by the choice of one of the two orientations of each traversed circle and possibly by a choice of regular point of $C$ when the initial (terminal) vertex of the path is a circle $C$.
\end{remark}

\begin{lemma}
\label{lem.tree.of.circles.is.tree}
Given a rabbit $\R$, the graph $T_\R$ is a tree.
Moreover, if a vertex is a circle then its degree is countably infinite and if it is instead a cut point then its degree is the order of the cut point.
\end{lemma}

\begin{proof}
Since rabbits are path-connected (because they are Peano continua), by the previous \cref{rmk.arcs.and.paths} the graph $T_\R$ is connected.
Now, a cycle in $T_\R$ would need to contain at least two distinct cut points and two distinct circles as vertices.
This would correspond to two distinct circles being joined by two arcs that travel through distinct cut points.
Then the two circles would not be separated by a cut point, which would violate condition 2 of \cref{def.rabbit}.
Hence, $T_\R$ is a tree.

By \cref{lem.cut.points.order}, each cut point $p$ belongs to one circle for each of the connected components of $\R \setminus \{p\}$.
Thus, the degree of $p$ in $T_\R$ is the order of the cut point.
As noted in \cref{rmk.rabbit.countable}, each circle contains countably many cut points, so its degree in $T_\R$ is infinite.
\end{proof}

\begin{remark}
\label{rmk.rabbit.end.points}
The set of end points of rabbits (\cref{def.rabbit.points}) corresponds to the set $\partial T$ of ends of the tree of circles (i.e., the set of semi-infinite paths from an arbitrary fixed vertex).

Indeed, if $e$ is an end point and $C_0$ is an arbitrary circle, consider an arc $I$ from $C_0$ to $e$. By \cref{lem.rabbit.arcs.and.circles}, $I \setminus \{e\}$ is included in a union of circles. The arc $I$ thus traverses (\cref{def.traverse}) an infinite sequence of distinct circles $(C_1, C_2, \dots)$ which corresponds, by \cref{rmk.arcs.and.paths}, to a semi-infinite path in $T_\R$ from the vertex $C_0$.
Conversely, a semi-infinite path in $T_\R$ from $C_0$ corresponds to an infinite sequence of distinct circles $(C_1, C_2, \dots)$.
If $B_k$ is the unique connected component of $\R \setminus C_k$ that meets $C_{k+1}$, then the intersection of all the $B_k$'s is non-empty (because each finite intersection $B_0 \cap \dots \cap B_k$ is non-empty and $\R$ is compact).
Points in the intersection of all the $B_k$'s are end points, as they cannot belong to a circle by condition 3 of \cref{def.rabbit}, so the correspondence between end points and semi-infinite paths from $C_0$ is surjective.
Finally, given two distinct end points $e_1$ and $e_2$, since $\R$ is a metrizable space the arcs $I_1$ and $I_2$ from $C_0$ to $e_1$ and $e_2$, respectively, must eventually diverge.
Then the sequences of circles that they meet cannot be the same, so neither can the semi-infinite paths in $T_\R$.

This argument also shows that, given an end point $e$ that corresponds to a sequence $(C_0, C_1, \dots)$ of circles, the set $\{ B_k \}_{k \in \mathbf{N}}$ (each $B_k$ as described above) is a basis of neighborhoods at $e$.
\end{remark}

\subsection{Legal coloring of the tree of circles}

In this and the following sections, we focus our attention on the regular rabbits, which include the basilica and the Douady rabbit.

\begin{definition}
\label{def.biregular.tree}
Any tree $T$ has a unique natural bipartition $\{A,B\}$ of vertices in which any two vertices lie in the same part of the partition if and only if their distance is even.
A tree $T$ is $(n,m)$\textbf{-biregular} if every vertex of $A$ has degree $n$ and every vertex of $B$ has degree $m$, where $n, m \in \mathbf{N} \cup \{\infty\}$.
\end{definition}

\begin{remark}
For every rabbit $\R$, its tree of circles $T_\R$ is bipartitioned into the sets of circles and of cut points of $\R$.
By \cref{lem.tree.of.circles.is.tree}, $T_\R$ is an $(n,\infty)$-biregular tree if and only if $\R$ is an $n$-regular rabbit.
\end{remark}

For each circle $C$ of a rabbit, we want to encode the cyclic structure of the cut points belonging to $C$ in the tree $T_\R$.
We will do this by using certain colorings of $T_\R$.
We cannot employ the renowned construction of \cite{BM00} because the trees of circles are not regular trees of finite degree.
We instead use the technology developed in \cite{biregular}, which generalizes the construction of Burger and Mozes to the case of biregular trees and allows infinite degree.
Let us first describe how this works.
 
\begin{definition}
\label{def.half.edges}
Given a graph $T$ (undirected and without loops nor parallel edges), a \textbf{half-edge} of $T$ is an ordered pair of adjacent vertices.
If $T$ is a graph, we denote by $H(T)$ the set of half-edges of $T$.
Given a vertex $v$ of a graph, we denote by $\mathrm{out}(v)$ the set of half-edges that originate from $v$ and by $\mathrm{in}(v)$ the set of half-edges that terminate at $v$.
\end{definition}

In \cite{biregular} half-edges are instead called \textit{arcs}, but we avoid this term because we are already using it for the arcs of topological spaces.
The term half-edge is inspired by the fact that each edge between $v$ and $w$ corresponds precisely to two half-edges, namely $(v,w)$ and $(w,v)$.

\begin{remark}
\label{rmk.rabbit.tree.separation.relation}
Let $T_\R$ be the tree of circles of a rabbit.
If $v = C$ is a circle, then we can naturally identify $\mathrm{out}(C)$ with the set of cut points that belong to $C$.
If instead $v = p$ is a cut point, then we identify $\mathrm{out}(p)$ with the set of circles that $p$ belongs to, or equivalently with the set of connected components of $\R \setminus \{p\}$.

Each circle $C$ has two natural cyclic orders, each inducing the same separation relation (see \cref{sec.order.and.separation}).
Thanks to condition 1 of \cref{def.rabbit}, restricting these relations to the set of cut points belonging to $C$ yields a countable dense cyclic order.
With the identification of $\mathrm{out}(C)$ with the set of cut points that belong to $C$, we thus have two natural dense cyclic orders inducing a unique separation relation on $\mathrm{out}(C)$, for each circle $C$ of the rabbit $\R$
\end{remark}

In order to equip $T_\R$ with the dense cyclic order on $\mathrm{out}(C)$, we rely on the tool defined right below.

\begin{definition}[\cite{biregular}]
\label{def.legal.coloring}
Let $T$ be an $(n,m)$-biregular tree ($n, m \in \mathbf{N} \cup \{\infty\}$) with bipartition $\{A,B\}$.
A \textbf{legal coloring} of $T$ with colors $K_n$ and $K_m$ is a coloring $k$ of its half-edges, i.e., a map
\[ k \colon H(T) \to K_n \cup K_m, \]
that enjoys the following properties:
\begin{enumerate}
    \item for all $v \in A$, $k$ maps $\mathrm{out}(v)$ bijectively to $K_n$;
    \item for all $v \in B$, $k$ maps $\mathrm{out}(v)$ bijectively to $K_m$;
    \item for all $v \in V(T)$, $k$ is constant on $\mathrm{in}(v)$.
\end{enumerate}
\end{definition}

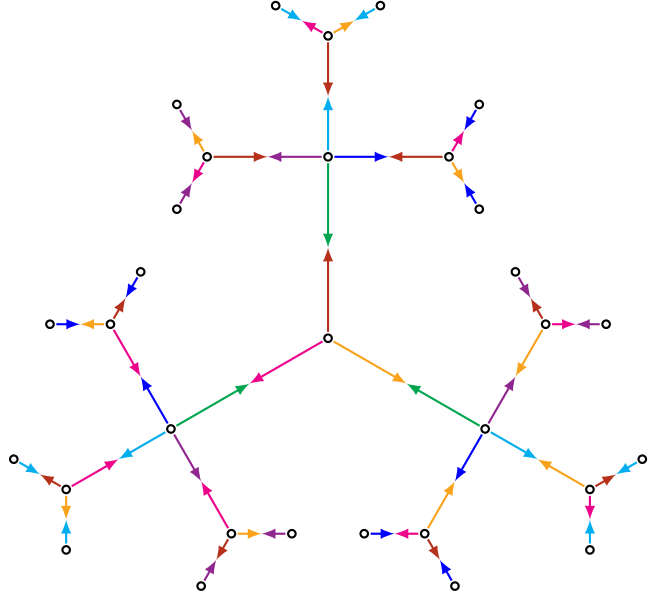
\begin{figure}
\centering
\begin{tikzpicture}[scale=.8]
    \node[vertex] (0) at (0:0) {};
    \node[vertex] (1) at (90:3) {};
        \node[vertex] (11) at ($(1)+(0:2)$) {};
            \node[vertex] (111) at ($(11)+(60:1)$) {};
            \node[vertex] (112) at ($(11)+(-60:1)$) {};
        \node[vertex] (12) at ($(1)+(90:2)$) {};
            \node[vertex] (121) at ($(12)+(150:1)$) {};
            \node[vertex] (122) at ($(12)+(30:1)$) {};
        \node[vertex] (13) at ($(1)+(180:2)$) {};
            \node[vertex] (131) at ($(13)+(240:1)$) {};
            \node[vertex] (132) at ($(13)+(120:1)$) {};
    \node[vertex] (2) at (210:3) {};
        \node[vertex] (21) at ($(2)+(120:2)$) {};
            \node[vertex] (211) at ($(21)+(180:1)$) {};
            \node[vertex] (212) at ($(21)+(60:1)$) {};
        \node[vertex] (22) at ($(2)+(210:2)$) {};
            \node[vertex] (221) at ($(22)+(270:1)$) {};
            \node[vertex] (222) at ($(22)+(150:1)$) {};
        \node[vertex] (23) at ($(2)+(300:2)$) {};
            \node[vertex] (231) at ($(23)+(0:1)$) {};
            \node[vertex] (232) at ($(23)+(240:1)$) {};
    \node[vertex] (3) at (330:3) {};
        \node[vertex] (31) at ($(3)+(240:2)$) {};
            \node[vertex] (311) at ($(31)+(300:1)$) {};
            \node[vertex] (312) at ($(31)+(180:1)$) {};
        \node[vertex] (32) at ($(3)+(330:2)$) {};
            \node[vertex] (321) at ($(32)+(30:1)$) {};
            \node[vertex] (322) at ($(32)+(270:1)$) {};
        \node[vertex] (33) at ($(3)+(60:2)$) {};
            \node[vertex] (331) at ($(33)+(120:1)$) {};
            \node[vertex] (332) at ($(33)+(0:1)$) {};
    \draw[edge,BrickRed] (0) to ($(0)!0.5!(1)$);
    \draw[edge,Magenta] (0) to ($(0)!0.5!(2)$);
    \draw[edge,YellowOrange] (0) to ($(0)!0.5!(3)$);
    \draw[edge,Green] (1) to ($(1)!0.5!(0)$);
    \draw[edge,blue] (1) to ($(1)!0.5!(11)$);
    \draw[edge,Cyan] (1) to ($(1)!0.5!(12)$);
    \draw[edge,Plum] (1) to ($(1)!0.5!(13)$);
        \draw[edge,BrickRed] (11) to ($(11)!0.5!(1)$);
        \draw[edge,Magenta] (11) to ($(11)!0.5!(111)$);
        \draw[edge,YellowOrange] (11) to ($(11)!0.5!(112)$);
        \draw[edge,blue] (111) to ($(11)!0.5!(111)$);
        \draw[edge,blue] (112) to ($(11)!0.5!(112)$);
        \draw[edge,BrickRed] (12) to ($(12)!0.5!(1)$);
        \draw[edge,Magenta] (12) to ($(12)!0.5!(121)$);
        \draw[edge,YellowOrange] (12) to ($(12)!0.5!(122)$);
        \draw[edge,Cyan] (121) to ($(12)!0.5!(121)$);
        \draw[edge,Cyan] (122) to ($(12)!0.5!(122)$);
        \draw[edge,BrickRed] (13) to ($(13)!0.5!(1)$);
        \draw[edge,Magenta] (13) to ($(13)!0.5!(131)$);
        \draw[edge,YellowOrange] (13) to ($(13)!0.5!(132)$);
        \draw[edge,Plum] (131) to ($(131)!0.5!(13)$);
        \draw[edge,Plum] (132) to ($(132)!0.5!(13)$);
    \draw[edge,Green] (2) to ($(2)!0.5!(0)$);
    \draw[edge,blue] (2) to ($(2)!0.5!(21)$);
    \draw[edge,Cyan] (2) to ($(2)!0.5!(22)$);
    \draw[edge,Plum] (2) to ($(2)!0.5!(23)$);
        \draw[edge,Magenta] (21) to ($(21)!0.5!(2)$);
        \draw[edge,YellowOrange] (21) to ($(21)!0.5!(211)$);
        \draw[edge,BrickRed] (21) to ($(21)!0.5!(212)$);
        \draw[edge,blue] (211) to ($(211)!0.5!(21)$);
        \draw[edge,blue] (212) to ($(212)!0.5!(21)$);
        \draw[edge,Magenta] (22) to ($(22)!0.5!(2)$);
        \draw[edge,YellowOrange] (22) to ($(22)!0.5!(221)$);
        \draw[edge,BrickRed] (22) to ($(22)!0.5!(222)$);
        \draw[edge,Cyan] (221) to ($(221)!0.5!(22)$);
        \draw[edge,Cyan] (222) to ($(222)!0.5!(22)$);
        \draw[edge,Magenta] (23) to ($(23)!0.5!(2)$);
        \draw[edge,YellowOrange] (23) to ($(23)!0.5!(231)$);
        \draw[edge,BrickRed] (23) to ($(23)!0.5!(232)$);
        \draw[edge,Plum] (231) to ($(231)!0.5!(23)$);
        \draw[edge,Plum] (232) to ($(232)!0.5!(23)$);
    \draw[edge,Green] (3) to ($(3)!0.5!(0)$);
    \draw[edge,blue] (3) to ($(3)!0.5!(31)$);
    \draw[edge,Cyan] (3) to ($(3)!0.5!(32)$);
    \draw[edge,Plum] (3) to ($(3)!0.5!(33)$);
        \draw[edge,YellowOrange] (31) to ($(31)!0.5!(3)$);
        \draw[edge,BrickRed] (31) to ($(31)!0.5!(311)$);
        \draw[edge,Magenta] (31) to ($(31)!0.5!(312)$);
        \draw[edge,blue] (311) to ($(311)!0.5!(31)$);
        \draw[edge,blue] (312) to ($(312)!0.5!(31)$);
        \draw[edge,YellowOrange] (32) to ($(32)!0.5!(3)$);
        \draw[edge,BrickRed] (32) to ($(32)!0.5!(321)$);
        \draw[edge,Magenta] (32) to ($(32)!0.5!(322)$);
        \draw[edge,Cyan] (321) to ($(321)!0.5!(32)$);
        \draw[edge,Cyan] (322) to ($(322)!0.5!(32)$);
        \draw[edge,YellowOrange] (33) to ($(33)!0.5!(3)$);
        \draw[edge,BrickRed] (33) to ($(33)!0.5!(331)$);
        \draw[edge,Magenta] (33) to ($(33)!0.5!(332)$);
        \draw[edge,Plum] (331) to ($(331)!0.5!(33)$);
        \draw[edge,Plum] (332) to ($(332)!0.5!(33)$);
\end{tikzpicture}
\caption{A legal coloring of a finite portion of a $(3,4)$-regular tree.}
\label{fig.colored.tree}
\end{figure}

For example, \cref{fig.colored.tree} portrays the legal coloring of a finite portion of a $(3,4)$-biregular tree, where each half-edge $(v,w)$ is represented by the half of the edge $\{v,w\}$ that originates from $v$ and whose arrow points towards $w$.

Note that, by \cite[Proposition 11]{biregular}, legal colorings form an orbit under the action of $\Aut(T)$, i.e., every two legal colorings with the same set of colors differ from an automorphism of $T$ and the image of a legal coloring is another legal coloring.
The third condition of \cref{def.legal.coloring} is key to this.

\begin{lemma}
\label{lem.coloring.of.tree.of.circles}
Let $\R$ be an $n$-regular rabbit.
Let $\Omega$ be a countable set together with a dense cyclic order $O_\Omega$ and, if $n$ is finite, let $[n] = \{1, \dots, n\}$ be equipped with its standard cyclic order.
For each circle $C$ of $\R$, equip $\mathrm{out}(C)$ with a countable dense cyclic order on the cut points of $C$ as described in \cref{rmk.rabbit.tree.separation.relation} and, if $n$ is finite, for each cut point $p$ of $\R$, equip $\mathrm{out}(p)$ with a cyclic order.
Then there is a legal coloring with colors $([n],\Omega)$ whose restriction to $\mathrm{out}(v)$ is an isomorphism of cyclic orders for each vertex $v$ (only for those that correspond to circles, if $n$ is infinite).
In particular, by \cref{lem.preserves.respects.cyclic.order} it is an isomorphism of the induced separation relations (\cref{def.separation.relation}).
\end{lemma}

\begin{proof}
We are going to repeatedly use the fact that the set of cut points of each circle $C$ is naturally endowed with two countable dense cyclic orders that induce the same separation relation.
On each circle, we fix any of the two.
This induces a unique separation relation on each $\mathrm{out}(C)$, as described in \cref{rmk.rabbit.tree.separation.relation}.

Fix a starting circle $C_0$ of $\R$.
The correspondence between the set $\mathrm{out}(C_0)$ of half-edges of $T_\R$ originating at $C_0$ and the set of cut points of $C_0$ induces the desired bijection $\mathrm{out}(C) \to \Omega$ that preserves the cyclic orders.

Proceeding by induction, suppose that $T$ is a subtree of $T_\R$ whose leaves are circles of $\R$ and that we have already built the map $k$ on the half-edges originating from $T$ (i.e., on the half-edges $(v,w) \in \mathrm{out}(v)$ for each $v \in V(T)$, even if $w$ does not belong to $T$) in such a way that
\begin{enumerate}
    \item for each circle $C$ of $T$, $k$ maps $\mathrm{out}(C)$ bijectively to $\Omega$ and preserves the cyclic orders;
    \item for each cut point $p$ of $T$, $k$ maps $\mathrm{out}(p)$ bijectively to $[n]$ and, if $n$ is finite, preserves the cyclic orders;
    \item for each vertex $v$ of $T$, $k$ is constant on $\mathrm{in}(v)$.
\end{enumerate}
Suppose that $p$ is a cut point that does not belong to $T$ and is adjacent in $T_\R$ to some circle $C$ of $T$.
We need to show that we can add $p$ to $T$ along with all circles that are adjacent to it and extend the coloring in such a way that the above conditions are satisfied.
This will eventually color every vertex of $T_\R$.
Given $p$ and $C$ as above, define $k(p,C)$ as $k(q,C)$ for any other cut point $q$ in $T$ that belongs to $C$ (in the first step, when $V(T) = \{C\}$, this choice is instead arbitrary).
Then complete $k(p,-)$ to a bijection such that $k(p,C')$ satisfies condition 2.
Next, for each circle $C' \neq C$ that is adjacent to $p$ define $k(C',p) = k(C,p)$, so that condition 3 holds at the vertex $p$.
For each $C'$, color $\mathrm{out}(C')$ in such a way that condition 1 holds.
The map $k$ extended in this way satisfies all conditions, so we are done.
\end{proof}

In the remainder of this section, which of the two cyclic orders we fix on each circle will not matter, as they induce the same separation relation.
The cyclic order on $\mathrm{out}(p)$ (for $p$ a cut point) too is not going to matter here.
These are only going to come into play in \cref{sec.lamination.automorphisms}.

\subsection{Universal groups of biregular trees}
\label{sub.universal.groups}

A graph automorphism $g$ permutes the half-edges around each vertex.
In particular, for any vertex $v$ we consider the bijection $g|_{\mathrm{out}(v)} \colon \mathrm{out}(v) \to \mathrm{out}(g(v))$.

\begin{definition}
\label{def.universal.group}
Let $T$ be an $(n,m)$-biregular tree with bipartition $\{A, B\}$ and fix a legal coloring $k$ of $T$.
Given $N \leq \Sym(K_n)$ and $M \leq \Sym(K_m)$, the associated \textbf{universal group} with prescribed local actions $\U_k (N,M)$ (or simply $\U (N,M)$ if $k$ is understood) is the group of those automorphisms $g$ of $T$ that preserve the bipartition and such that, for all $v \in V(T)$,
\begin{equation}
    k|_{\mathrm{out}(g(v))} \circ g|_{\mathrm{out}(v)} \circ k|_{\mathrm{out}(v)}^{-1} \in
    \begin{cases}
      N & \text{if $v \in A$},\\
      M & \text{if $v \in B$},
    \end{cases}
\end{equation}
i.e., the permutation of colors induced by $g$ at the vertex $v$ belongs to $N$ or $M$, depending on whether $v$ belongs to $A$ or $B$.
We endow $\U_k(N,M)$ with the pointwise convergence topology.
\end{definition}

We will denote the color permutation induced by $g$ around $v$ by
\[ \sigma_g(v) \coloneqq k|_{\mathrm{out}(g(v))} \circ g|_{\mathrm{out}(v)} \circ k|_{\mathrm{out}(v)}^{-1}. \]
This is called the \textbf{local action} of $g$ at $v$.

The paper \cite{biregular} studies the topological and transitive properties that $\U(N,M)$ inherits from $N$ and $M$.
It also shows that $\U(N,M)$ is topologically isomorphic to its permutation group of the vertices of one of the two parts of the bipartition, which Smith calls the \textit{box product} and denotes by $N \boxtimes M$.

\begin{remark}
\label{rmk.BM.and.biregular}
As noted in the third paragraph of Remark 4 of \cite{biregular}, when $n=2$ there is natural identification between the universal group $\U(\Sym([2]),M)$ and the Burger-Mozes group $\U(M)$ that acts on an $m$-regular tree.
More precisely, $\U(\Sym([2]),M)$ and $\U(M)$ are permutation isomorphic if we consider the action of the first on the set of vertices of degree $m$.
This identification is found simply by taking the barycentric subdivision of the $m$-regular tree.

In particular, it is not hard to see that the group $\U(\Sym([2]),\Aut(S))$ that we will describe in the next subsection for the basilica (which is the $2$-regular rabbit) is the same as the Burger-Mozes group found in \cite{Neretin}.
\end{remark}

\subsection{Uniqueness of the regular rabbits}

Existence of $n$-regular rabbits is proved in \cref{sec.julia.sets,,sec.limit.spaces}, when $n$ is finite.
On the other hand, an $\infty$-rabbit cannot be a Julia set nor a limit space (indeed, points of Julia sets and limit spaces have finite order).
However, $\infty$-rabbits exist as topological spaces.
Let us sketch a construction of $\R_\infty$.
Consider the biregular tree $T=T_{\infty,\infty}$ and fix an embedding of $T$ in the plane such that each edge is realized as a segment and for each vertex the set of directions of segments attached to it is dense in the projective line.
Let $A$ and $B$ be the sets of the bipartition.
For a finite subset $F$ of $B$, let $X_F$ be the topological space obtained by blowing up in the plane all vertices in $F$.
This gives an inverse system.
Let $X_\infty$ be the inverse limit.
Consider the equivalence relation that identifies all points lying on the same edge.
The quotient space is an $\infty$-regular rabbit.

The following proposition shows that there is a unique $n$-regular rabbit for each $n$ (up to homeomorphisms).
Given a rabbit $\R$ and its tree of circles $T$, below we denote by $\pi \colon \R \to \overline{T}$ the map that sends each end point to the corresponding point of $\partial T$, each cut point to itself (as a vertex of $T$) and each regular point to the unique circle to which it belongs (as a vertex of $T$).

\begin{proposition}
\label{prop.rabbit.surjectivity}
Let $\R_1$ and $\R_2$ be two $n$-regular rabbits and let $T$ be their tree of circles (which is the same up to isomorphisms by \cref{lem.tree.of.circles.is.tree}).
Let $k_1$ and $k_2$ be two legal colorings of $T$ as in \cref{lem.coloring.of.tree.of.circles}.
Then, for every automorphism $\varphi$ of $T$ whose local action on $\mathrm{out}(v)$ at each vertex $v$ is an isomorphism between the two separation relations, there is a homeomorphism $\psi \colon \R_1 \to \R_2$ that makes the following diagram commute.
\begin{center}
\begin{tikzcd}
    \R_1 \arrow[r, "\psi"] \arrow[d, "\pi_1"] & \R_2 \arrow[d, "\pi_2" ] \\
    T \arrow[r, "\varphi"] & T
\end{tikzcd}
\end{center}
\end{proposition}

\begin{proof}
Let $\R_1$ and $\R_2$ be two $n$-regular rabbits and let $T_1$ and $T_2$ be their trees of circles.
Fix two legal colorings satisfying \cref{lem.coloring.of.tree.of.circles} and let $\varphi \colon T_1 \to T_2$ be a color-preserving automorphism (the local action at any vertex is the identity).
Note that the maps $\pi_1$ and $\pi_2$ are injective on the sets of cut points of $\R_1$ and $\R_2$.
This defines a bijection $\psi = \pi^{-1}_2 \circ \varphi \circ \pi_1$ between the cut points of $\R_1$ and $\R_2$.
Since $\varphi$ is a tree isomorphism and the edge adjacency in $T_1$ and $T_2$ is given by the inclusion of cut points to circles, the map $\psi$ restricts to bijections between the cut points of each circle.
Moreover, since $\varphi$ is color-preserving, the circular order of cut points on each circle is preserved, and since the sets of cut points are dense in $\R_1$ and $\R_2$ (by condition 1 of \cref{def.rabbit} together with \cref{lem.coloring.of.tree.of.circles}), $\psi$ extends uniquely to a continuous map between a circle of $\R_1$ and its image.
Doing this for all circles of $\R_1$ extends $\psi$ to a bijection between $\R_1$ and $\R_2$ such that $\varphi \circ \pi_1 = \pi_2 \circ \psi$.
It now remains to show that $\psi$ is a homeomorphism.

Note that the connected components of $\R_i$ minus finitely many circles correspond precisely to the connected components of $T_i$ minus finitely many circles together with their cut points (seen as vertices of the tree).
Then, if $U$ is such a subset of $\R_2$, its preimage $\psi^{-1}(U)$ is a connected component of $\R_1$ minus finitely many circles, so in particular it is open.
By \cref{cor.rabbit.subbasis}, the collection of such subsets is a subbasis for the topology of a rabbit, so this shows that $\psi$ is continuous.
Ultimately, this shows that $\psi$ is a continuous bijection from a compact space to a Hausdorff space, so it is a homeomorphism.
\end{proof}

By \cref{lem.coloring.of.tree.of.circles}, we can always find legal colorings of any two $n$-regular rabbits satisfying the hypotheses of \cref{prop.rabbit.surjectivity}.
Hence we obtain the desired result.

\begin{theorem}
\label{thm.uniqueness.rabbit}
Any two $n$-regular rabbits are homeomorphic.
\end{theorem}

Thanks to this theorem, we can now refer to any $n$-regular rabbit as \textit{the} $n$-regular rabbit.
For example, \textit{the} basilica and \textit{the} Douady rabbit.

\subsection{Homeomorphism groups of regular rabbits}
\label{sub.rabbit.homeo.group.is.universal}

Here we show that the group of homeomorphisms of a regular rabbit is the universal group $\U(\Sym([n]),\Aut(S))$, where $S$ is the separation relation of the countable dense cyclic order.

\begin{proposition}
\label{prop.homeo.rabbit.acts.on.tree}
Given an $n$-regular rabbit $\R$, the group $\Homeo(\R)$ acts faithfully by graph automorphisms on $T_\R$ in such a way that defines an embedding
\[ \Pi \colon \Homeo(\R) \hookrightarrow \U_k (\Sym([n]),\Aut(S)), \]
where $k$ is the legal coloring of $T_\R$ provided by \cref{lem.coloring.of.tree.of.circles}
\end{proposition}

\begin{proof}
A homeomorphism $\varphi$ of $\R$ induces a permutation of the circles and of the cut points, which are the vertices of the tree of circles $T_\R$.
Recall that $\{p,C\}$ is an edge of $T_\R$ if and only if $p \in C$, so clearly the permutation induced by $\varphi$ preserves the edge adjacency of $T_\R$.
Thus, $\Homeo(\R)$ acts by graph automorphisms on $T_\R$.
The set of cut points is dense in $\R$ by condition 1 of \cref{def.rabbit} and \cref{cor.rabbit.cut.points}, so the action is faithful.

Now, consider the legal coloring $k$ of $T_\R$ provided by \cref{lem.coloring.of.tree.of.circles}.
The bipartition of $T_\R$ distinguishes between cut points and circles, so clearly the action of $\Homeo(\R)$ preserves it.
For each circle $C$ of $\R$, every homeomorphism $\varphi$ of $\R$ maps bijectively between the sets of cut points belonging to $C$ and those belonging to $\varphi(C)$ in a separation-preserving fashion.
This means that the permutation of colors induced by $\varphi$ around the vertex $C$ of $T_\R$ belongs to $\Aut(S)$, as needed.
\end{proof}

The surjectivity of this embedding follows immediately from \cref{prop.rabbit.surjectivity}, so we have the following fact.

\begin{proposition}
\label{prop.rabbit.embedding.is.surjective}
The embedding $\Pi$ of \cref{prop.homeo.rabbit.acts.on.tree} is surjective.
\end{proposition}

We are now only missing the bicontinuity of $\Pi$.
We endow $\Homeo(\R)$ with the uniform convergence topology.

\begin{proposition}
\label{prop.rabbit.embedding.is.continuous}
The embedding $\Pi$ of \cref{prop.homeo.rabbit.acts.on.tree} is a homeomorphism.
\end{proposition}

\begin{proof}
Since both are Polish groups, it suffices to prove that the abstract group isomorphism $\Pi \colon \Homeo(\R) \to \U(\Sym([n]),\Aut(S))$ is continuous (for example because of \cite[Corollary 2.3.4]{MR2455198}).
Suppose that $g_m$ is a sequence in $\Homeo(\R)$ that converges to the identity and consider an arbitrary circle $C_0$ of the rabbit $\R$.
We will show that, for large enough $m \in \mathbf{N}$, the homeomorphism $g_m$ fixes $C_0$.
Let $C_1$, $C_2$ and $C_3$ be circles of $\R$ that lie in three distinct connected components $U_1$, $U_2$ and $U_3$ of $\R \setminus C_0$.
Then, for $m$ large enough, say $m \geq N$, the circle $g_m (C_i)$ is included in $U_i$ for each $i = 1, 2, 3$.
In the tree of circles, consider the unique paths between $C_1$ and $C_2$, between $C_1$ and $C_3$ and between $C_2$ and $C_3$:
by \cref{rmk.arcs.and.paths}, they intersect precisely at $C_0$, so the vertex $C_0$ is fixed by $\Pi(g_m)$ for all $m \geq N$.
This argument shows that each vertex is fixed by $\Pi(g_m)$ for large enough $m$, so the sequence $\Pi(g_m)$ converges to the identity of $\U(\Sym([n]),\Aut(S))$ and thus $\Pi$ is continuous.
\end{proof}

Together, the previous propositions immediately imply the desired result.

\begin{theorem}
\label{thm.homeo.rabbits.are.auto.trees}
The group $\Homeo(\R_n)$ of homeomorphisms of an $n$-regular rabbit $\R_n$ is topologically isomorphic to $\U (\Sym([n]),\Aut(S))$.
\end{theorem}

\subsection{Consequences of the identification with universal groups}\label{sec.consequenses.identification.rabbit}

\cref{thm.homeo.rabbits.are.auto.trees} has the following consequence.

\begin{corollary}
\label{cor.rabbits.simple}
Let $\R_n$ be the $n$-regular rabbit.
\begin{itemize}
    \item If $n=2$, the group $\Homeo(\R_2) = \Homeo(\B)$ has a simple index-$2$ commutator subgroup $\Homeo(\B)'$.
    \item If $n \geq 3$, the group $\Homeo(\R_n)$ is simple.
\end{itemize}
\end{corollary}

\begin{figure}
    \centering
    \includegraphics[width=.7\textwidth]{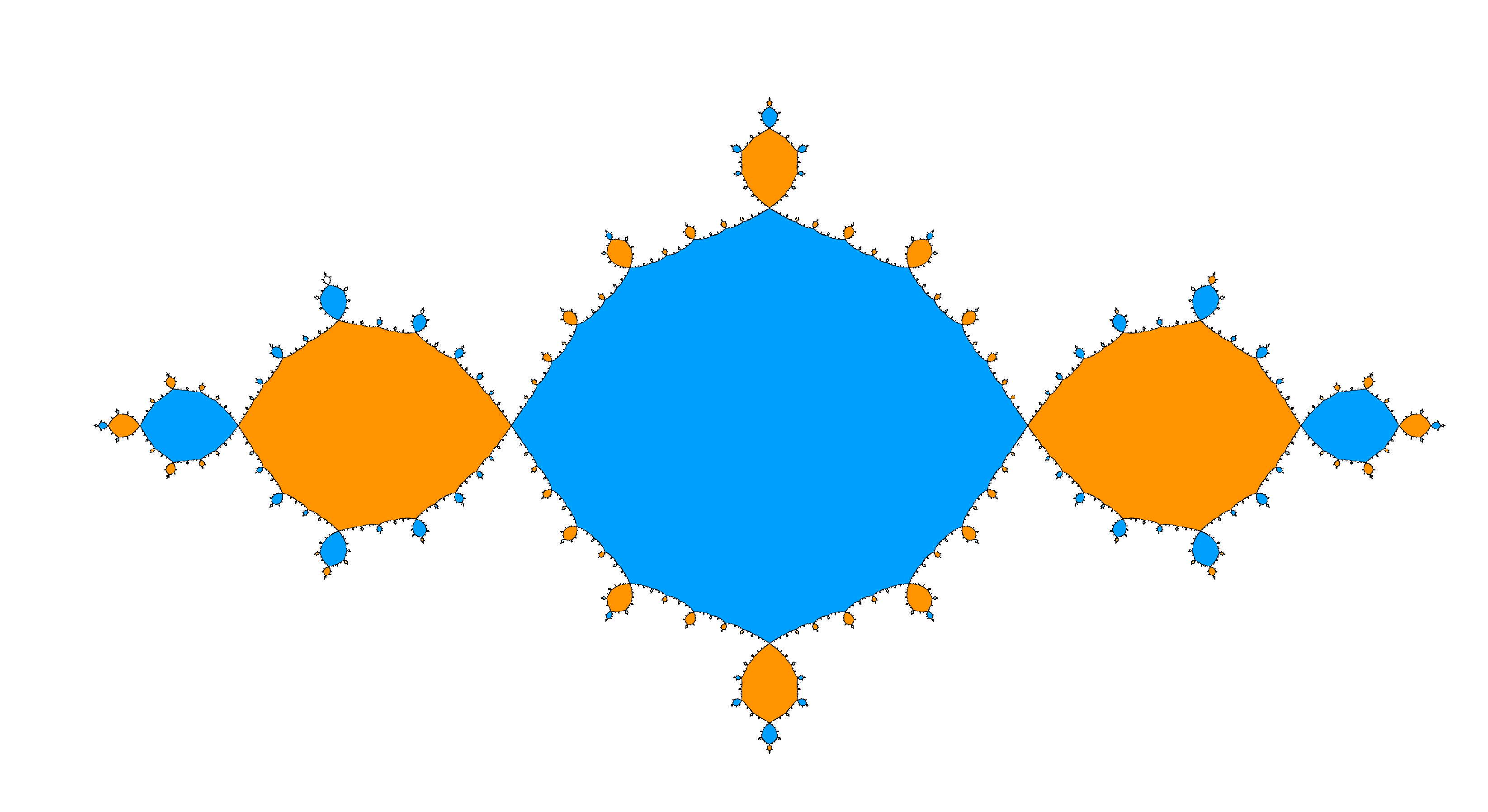}
    \caption{The alternating $2$-coloring of the circles of $\B$.}
    \label{fig.julia_basilica_colored}
\end{figure}

\begin{proof}
Theorem 23 of \cite{biregular} states that $\U(N,M)$ is simple if and only if $N$ or $M$ is transitive, when assuming that at least one of $N$ and $M$ is non-trivial and that they are both generated by their point stabilizers.
In our case, $M = \Aut(S)$ is transitive, non-trivial and generated by its point stabilizers (\cref{cor.aut.circular.order.generated.by.stabilizers}) and, as soon as $n \geq 3$, $\Sym([n])$ is generated by its point stabilizers too.
This shows that $\Homeo(\R_n)$ is simple as soon as $n \geq 3$.

Consider now the case $n=2$, which concerns the basilica $\R_2 = \B$.
Neretin noted that $\Homeo(\B)$ is the universal Burger-Mozes group $\U(\Aut(S))$ on a regular tree \cite{Neretin}, which has an index-$2$ simple subgroup generated by the edge stabilizers by \cite[Proposition 3.2.1]{BM00}.
The regular tree in question is a simplified tree of circles whose vertices are the circles of the basilica and whose edges connect adjacent circles.
Hence, the group generated by the pointwise stabilizers of edges corresponds to the group that preserves the alternating $2$-coloring of the circles depicted in \cref{fig.julia_basilica_colored}.
This is the kernel of the group homomorphism $f \colon \Homeo(\B) \to \Sym([2]) = \{0,1\}$ that has value $0$ or $1$ depending on whether the element preserves the coloring or not, which implies that $\Homeo(\B)' \leq \mathrm{Ker}(f)$.
Moreover, note that $\mathrm{Ker}(f) = \U(1,\Aut(S))$, so $\mathrm{Ker}(f)$ is simple by the aforementioned Theorem 23 of \cite{biregular}.
Clearly $\Homeo(\B)'$ is not trivial and it is a normal subgroup of $\mathrm{Ker}(f)$, so we can finally conclude that the index-$2$ simple subgroup $\mathrm{Ker}(f)$ of $\Homeo(\B)$ is none other than $\Homeo(\B)'$.
\end{proof}

What makes $\Homeo(\B)$ different from the other $\Homeo(\R_n)$ regarding their simplicity is that $\Sym([n])$ is cyclic if and only if $n=2$.
The upcoming \cref{thm.simple.commutators} relies precisely on this fact and shows how the cyclicity of the group of local actions comes into play.

The following is an additional easy consequence of \cref{thm.homeo.rabbits.are.auto.trees}.

\begin{corollary}
\label{cor.rabbit.groups.closed}
For every $n \geq 2$, the group $\Homeo(\R_n)$ embeds into $\Aut(T_{n,\infty})$ as a closed subgroup.
\end{corollary}

\begin{proof}
Theorem 6 (v) of \cite{biregular} states that $\U(N,M)$ is closed in $\Aut(T)$ if $N$ and $M$ are closed in $\Sym([n])$ and $\Sym([m])$, respectively.
Since $\Aut(S)$ is a closed subgroup of the infinite symmetric group (\cref{rmk.Aut.S.C.closed}), we are done.
\end{proof}

Further properties of the groups $\Homeo(\R_n)$ are explored later in \cref{sec.embedding,,sec.geometric.analytic.properties,,sec.universal.minimal.flows}.


\section{The airplane}
\label{sec.airplane}

In this section we discuss the topological characterization of the airplane fractal, which we gave in \cref{def.airplane} (it is shown in \cref{sec.julia.sets} that a Julia set satisfies this topological characterization).
Then we show that its homeomorphism group is isomorphic to a kaleidoscopic group.
Introduced and studied in \cite{Kaleidoscopic}, kaleidoscopic groups are groups of homeomorphisms of a dendrite with prescribed local actions at each branch point.

\subsection{Topological characterization of the airplane}

Recall the abstract definition of airplane that we gave \cref{def.airplane}.

The following Lemma essentially states that, unlike rabbits, \textit{traversing} (\cref{def.traverse}) and \textit{meeting} (i.e., having non-empty intersection) are the same for arcs in an airplane, possibly except for the end points of the arcs.

\begin{lemma}
\label{lem.airplane.meet.means.traverse}
In an airplane $\A$, whenever the interior $\mathring{A}$ of an arc $A$ meets a circle $C$, it traverses it (in the sense of \cref{def.traverse}, i.e., $\mathring{A} \cap C$ includes a non-degenerate arc).
\end{lemma}

\begin{proof}
By contradiction, assume that the intersection $\mathring{A} \cap C$ between the interior of an arc $A$ and a circle $C$ does not include non-degenerate arcs.
Then there is a subarc $A_0$ such that $\mathring{A_0}$ only meets $C$ in a single point, so we can simply assume that $\mathring{A} \cap C = \{p\}$ up to passing to subarcs.

Let us denote by $x$ and $y$ the two end points of $A$.
By condition 4 of \cref{def.airplane}, the order of $p$ is either $1$ or $2$.
In both cases, since $\A \setminus \{p\}$ features at most two connected components, either the two subarcs $[x,p)$ and $(p,y]$ of $A$ lie in the same connected component of $\A \setminus \{p\}$, or one of them lies in the connected component that includes $C \setminus \{p\}$.

In the first case, a connected component $X$ of $\A \setminus \{p\}$ includes both subarcs $[x,p)$ and $(p,y]$ of $A$.
In particular it contains $x$ and $y$ and is path-connected (because it is a connected and open in $\A$), so $X$ includes an arc $A^*$ between $x$ and $y$.
Thus, the subset $A^* \cup [x,p) \cup [y,p) \cup \{p\}$ of $X \cup \{p\}$ includes a circle that meets $C$ in $p$, which contradicts condition 3 of \cref{def.airplane}.

In the second case, a connected component $X$ of $\A \setminus \{p\}$ includes both $C \setminus \{p\}$ and a subarc $[x,p)$ or $(p,y]$ of $A$.
If we suppose for example that $[x,p) \subseteq X$, then $X$ includes an arc $[x,z]$ between $x$ and some $z \in C \setminus \{p\}$.
Let $[z,p]$ be any arc between $z$ and $p$ that is included in the circle $C$.
The subset $[x,z] \cup [z,p] \cup (p,x]$ includes a circle in $X \cup \{p\}$ that meets $C$ in $p$ and is distinct from $C$, which contradicts condition 3 of \cref{def.airplane}.
\end{proof}

The following fact will allow us to build paths of circles inside an airplane.

\begin{lemma}
\label{lem.airplane.arcs}
Given two points of an airplane, any two arcs between them meet the same circles.
\end{lemma}

\begin{proof}
Suppose that $\A$ is an airplane and let $x \neq y$ be points of $\A$.
We will show that, given an arc $A$ between $x$ and $y$, the set of circles that it meets is
\[ C(x,y) = \{ C \text{ circle } \mid x \in C \text{ or } y \in C \} \cup \{ C \text{ circle } \mid C \text{ separates } x \text{ and } y \}. \]
Since this set does not depend on the arc $A$ itself, this claim immediately implies the statement.

If $C$ is a circle that contains either $x$ or $y$ then clearly $A$ meets $C$.
If $C$ is a circle that separates $x$ and $y$, then every arc $A$ between $x$ and $y$ must meet $C$.
Hence, each circle that belongs to $C(x,y)$ is met by every arc between $x$ and $y$.

Conversely, let $C$ be a circle that is met by an arc $A$ between $x$ and $y$.
Assume that $C$ does not contain $x$ nor $y$, otherwise $C \in C(x,y)$ by definition and there is nothing to prove.
In particular, $A$ is not entirely included in $C$.
Suppose by contradiction that $C \notin C(x,y)$.
Then $C$ does not separate $x$ and $y$, so $\A \setminus C$ includes an arc $A^*$ between $x$ and $y$.
Thus, the subset $A^* \cup A$ of $\A$ includes a circle that meets $C$ but is not $C$, which contradicts condition 3 of \cref{def.airplane}.
\end{proof}

In an airplane, a \textbf{wing} is a connected component of the complement of a circle.
While wings alone do not form a subbasis for the topology of $\A$, a larger collection of similar sets do.
This is expressed in the following statement, which is an immediate corollary of \cref{prop.subbasis}.

\begin{corollary}
Given an airplane $\A$, the collection of all connected components of $\A$ minus finitely many of its circles and cut point that lie on circles is a subbasis for the topology of $\A$.
\end{corollary}

\subsection{The dendrite of circles of an airplane}
\label{sub.dendrite.of.circles}

Recall that a \textbf{dendrite} is a Peano continuum such that every two distinct points are joined by a unique arc.
A \textbf{branch point} of a dendrite is a point of order at least $3$.
We write $\Br(D)$ for the set of all branch points of a dendrite $D$.
A \textbf{branch} around $b \in \Br(D)$ is a connected component of $D \setminus \{b\}$.

The set of branch points of a dendrite comes equipped with a \textbf{betweenness relation}: $p$ is between $x$ and $y$ when $p$ belongs to the unique arc joining $x$ and $y$.
Every homeomorphism of $D$ induces a bijection of $\Br(D)$ that preserves the betweenness relation.
The converse is also true by \cite[Proposition 2.4]{Kaleidoscopic}, as long as every arc of $D$ contains branch points (which holds for all dendrites discussed in this paper), in which case $\Homeo(D)$ can be identified with the group of betweenness-preserving bijections of $\Br(D)$.

Given a subset $S \subseteq \mathbf{N}_{\geq 3} \cup \{\infty\}$, the \textbf{Wa\.zewski dendrite} associated to $S$ is a dendrite such the order of any branch point belongs to $S$ and that, for each $s \in S$, any arc of the dendrite contains points of order $s$.
Such dendrite is usually denoted by $D_S$ or simply by $D_n$ when $S = \{n\} \subseteq \mathbf{N}_{\geq 3} \cup \{\infty\}$.
By \cite[Theorem 6.2]{SelfHomeomorphicSpaces}, for each distinct $S$ there is a unique Wa\.zewski dendrite $D_S$ up to homeomorphism, so we can refer to $D_S$ as \textit{the} Wa\.zewski dendrite of order $S$.
Mimicking the construction of the tree of circles for rabbits, we now associate a copy of $D_\infty$ to the airplane (this construction was previously outlined in \cite[Section 7.3]{dendriterearrangement}).

\begin{definition}
\label{def.dendrite.of.circles}
Let $\A$ be an airplane.
The \textbf{circular equivalence relation} on $\A$ is defined by saying that two points are equivalent when they lie on a common circle.
The quotient topological space $Q$ of an airplane under the circular equivalence relation is called the \textbf{dendrite of circles}.
\end{definition}

The fact that $Q$ is really a dendrite is proved right below.

\begin{proposition}
\label{prop.airplane.quotient}
If $\A$ is an airplane, the quotient $Q$ of $\A$ under the circular equivalence relation is homeomorphic to the universal Wa\.zewski dendrite $D_\infty$.
Moreover, the quotient map induces a bijection between the set of wings of $\A$ and that of branches of $Q$.
\end{proposition}

\begin{proof}
We need to show that the quotient space $Q$ is a compact, metrizable, connected, and locally connected topological space where every two distinct points are joined by a unique arc and such that its branch points have infinite order and meet every arc.

Let us first show that $Q$ is Hausdorff.
Consider distinct $p_1$ and $p_2$ in $Q$.
Each of their preimages is either a circle or a single point of $\A$.
Whatever the case, by Conditions 1 and 2 of \cref{def.airplane} there exists a circle $C$ in $\A$ that separates them.
Then the two branches $B_1$ and $B_2$ of $Q$ at $C$ that contain $p_1$ and $p_2$, respectively, are open (indeed, their preimages are wings, which are open) and each only contains one of the two points, as needed.

Recall that a compact Hausdorff space is metrizable if and only if it is second countable.
Since $Q$ is a continuous image of a compact space, it is compact.
Second countability of $Q$ descends from the one of $\A$:
wings are open and saturated for the relation (i.e., they are preimages of subsets of $Q$), and using \cref{prop.subbasis} it is easy to show that their images form a countable subbase for the topology of $Q$.
Thus, $Q$ is metrizable.

Being a continuous image of a path-connected and locally path-connected space, $Q$ itself is path-connected and locally path-connected.
We claim that there is a unique arc between any two points of $Q$.
Let us say that a point $c\in Q$ is a \textbf{circle point} if its preimage is a circle (i.e., if it corresponds to a non-trivial equivalence class).
By condition 2 of \cref{def.airplane}, any two distinct points of $Q$ are separated by a circle point.
More precisely, a circle point $c$ separates $x,y$ if and only its preimage $C$ separates the preimages of $x$ and $y$ for the quotient map in $\A$.
Note that an arc $I$ between any two distinct points $x,y \in Q$ is the closure in $Q$ of the set of circle points between $x$ and $y$ by \cref{lem.airplane.arcs}.
This implies that there is a unique arc between any two points of $Q$, so $Q$ is a dendrite.

For each cut point $x$ on a circle $C$ of $\A$, the complement $C \setminus \{x\}$ is included in one connected component of $\A \setminus \{x\}$.
Hence, if $x$ and $y$ are two cut points on the same circle and $U_x$ and $U_y$ are the connected components of $\A\setminus\{x\}$ and $\A\setminus\{y\}$ that do not contain $x$ and $y$, then $U_x \cap U_y = \emptyset$.
Since there is an arc between $x$ and $y$ on $C$, there would be another arc going from $x$ to some $z\in U_x\cap U_y$ concatenated with a path from $z$ to $y$ and this would meet other circles than only $C$, which contradicts condition 3 of \cref{def.airplane}.
This implies that for each point of the dense set of cut points in $C$, we have connected saturated open subsets that are distinct and thus any circle point is a branch point of the dendrite $Q$.
Moreover, by condition 2 of \cref{def.airplane} the set of circle points meets every arc of $Q$ and the order of a circle point is infinite by condition 1.

In order to conclude that $Q$ is homeomorphic to $D_\infty$, it remains to show there is no branch point in $Q$ that is not a circle point.
Assume that $x \in \A$ does not belong to a circle, meaning that $x$ is the unique preimage of its image.
Then the one or two connected components of $\A \setminus \{x\}$ (they cannot be more than two by condition 4 of \cref{def.airplane}) are open and saturated.
Thus the image of $x$ has order $1$ or $2$, so it is not a branch point, concluding our proof.
\end{proof}

From the construction of the quotient map $\pi$, we can distinguish between four types of points in an airplane.

\begin{definition}
\label{def.points.airplane}
A point $x$ in an airplane is said to be:
\begin{itemize}
    \item an \textbf{end point} if $\pi(x)$ is an end point of $Q$ (equivalently, $x$ has order $1$ and does not lie on a circle),
    \item a \textbf{regular cut point} if $\pi(x)$ is a regular point of $Q$ (equivalently, $x$ does not lie on a circle and has order $2$),
    \item a \textbf{circular cut point} if $x$ is a cut point that lies on a circle.
    \item a \textbf{circular non-cut point} if $x$ lies on a circle but is not a cut point.
\end{itemize}
\end{definition}

\begin{remark}
One can see that, for airplanes, our circular equivalence relation coincides with the equivalence relation defined by Bowditch for general separable continuum.
In \cite[Theorem 5.23]{Bowditch} Bowditch showed that, if $M$ is a separable continuum, then $M / \sim$ is a dendrite, where $\sim$ is defined by setting $x \not\sim y$ if and only if there is a set of cut points of $M$ separating $x$ and $y$ which is order-isomorphic to the rational numbers in the natural linear order.
However, our circular equivalence relation is arguably more intuitively formulated in the context of the airplane.

Note that, on the other hand, Bowditch's quotient of a rabbit is just a singleton, so unfortunately the tree of circles of a rabbit and the dendrite of circles of a rabbit cannot both be formalized by this equivalence relation.
Finding a unique construction that realizes both the trees of circles for rabbits (\cref{def.tree.of.circles}) and the dendrite of circles for the airplane (\cref{def.dendrite.of.circles}) might allow to generalize this construction to Julia sets such as the one depicted in \cref{fig.quasi.rabbit}.
\end{remark}

\subsection{Kaleidoscopic coloring of the dendrite of circles}

By condition 1 of \cref{def.airplane}, the set of cut points on a fixed circle is dense in the circle, so the set of cut points on $C$ is equipped naturally with a countable dense cyclic order.
By condition 4, each cut point $x$ on a circle $C$ corresponds to a unique wing of $\A \setminus C$.
Together these two facts imply that, for every circle $C$, there is a countable dense cyclic order on the set of wings at $C$.

The bijection between wings of $\A$ and branches of $Q$ from \cref{prop.airplane.quotient} thus induces a countable dense cyclic order on the set of branches around each branch point of $Q$.
We claim that this corresponds to a coloring of $Q$, using the notions defined in \cite{Kaleidoscopic}.
Before showing this in \cref{lem.kaleidoscopic.coloring.order}, we briefly recall what kaleidoscopic colorings are and we introduce some useful notations.

Given a dendrite $D$ and a branch point $b$, we denote by $\hat{b}$ the set of branches around $b$ (i.e., the set of connected components of $D \setminus \{b\}$).
For two distinct branch points $b$ and $b'$ of a dendrite, we denote by $c_b(b')$ the unique branch around $b$ that contains $b'$.

\begin{definition}[Definitions 3.1 and 3.2 of \cite{Kaleidoscopic}]
\label{def.kaleidoscopic.coloring}
Given a set $\Omega$, a \textbf{kaleidoscopic coloring} of a Wa\.zewski dendrite $D_n$ ($n \in \mathbf{N}_{\geq 3} \cup \{\infty\}$) is a map
\[ k \colon \bigsqcup_{b\in\Br(D_n)} \hat{b} \to \Omega \]
such that
\begin{enumerate}
    \item for every branch point $b$, the restriction $k_b \colon \hat{b} \to \Omega$ of $k$ to the branches around $b$ is a bijection;
    \item for all distinct branch points $b_1, b_2$ and colors $i,j \in \Omega$, there exists a branch point $b$ separating $b_1$ and $b_2$ in distinct branches $B_1 = c_b(b_1), B_2 = c_b(b_2) \in \hat{b}$ such that $k_b(B_1) = i$ and $k_b(B_2) = j$.
\end{enumerate}
\end{definition}

Proposition 3.3 of \cite{Kaleidoscopic} shows that kaleidoscopic colorings exist (and are, in fact, common).
However, we will need specific kaleidoscopic colorings that represent the cyclic orders on circles.

\begin{lemma}
\label{lem.kaleidoscopic.coloring.order}
Let $(\Omega,O_\Omega)$ be a countable dense cyclic order.
Consider a collection $(O_b)_{b\in \Br(D_\infty)}$ of countable dense cyclic orders $O_b$ on $\hat{b}$ for each $b\in\Br(D_\infty)$.
There exists a kaleidoscopic coloring 
\[ k\colon \bigsqcup_{b\in\Br(D_\infty)} \hat{b} \to \Omega \]
such that, for any $b\in \Br(D_\infty)$ and $x,y,z\in \hat{b}$, 
\begin{equation}
\label{eq.kaleidoscopic.coloring.order}
    O_b(x,y,z)\iff O_\Omega(k(x),k(y),k(z)).
\end{equation}
\end{lemma}

\begin{proof}
The set of colorings satisfying \cref{eq.kaleidoscopic.coloring.order} is closed in the space of all colorings and thus is a Baire space.
Let $\mathbb{X}$ be the space of colorings satisfying \cref{eq.kaleidoscopic.coloring.order}.
For any two distinct branch points $b,b'$, the kaleidoscopic condition is open and dense, as in Proposition 3.3 of \cite{Kaleidoscopic}.
So the collection of kaleidoscopic colorings satisfying \cref{eq.kaleidoscopic.coloring.order} is a dense $G_\delta$ set in $\mathbb{X}$ and thus non-empty by the Baire property.
\end{proof}

\subsection{Kaleidoscopic groups}

Consider the Wa\.zewski dendrite $D_n$ with a kaleidoscopic coloring $k$ with set of colors $[n]$, where $n \in \mathbf{N}_{\geq 3} \cup \{\infty\}$.
The \textbf{local action} of a homeomorphism $g$ of $D_n$ at a branch point $b$ is the element $\sigma_g$ of $\Sym([n])$ defined as
\[ \sigma_g(b) \coloneqq k_{g(b)} \circ g \circ k_b^{-1}. \]
The set of elements with prescribed actions around each branch point forms a group, whose definition we recall right below.

\begin{definition}[Definition 3.8 of \cite{Kaleidoscopic}]
\label{def.kaleidoscopic.group}
Given a permutation group $\Gamma \leq \Sym([n])$, the \textbf{kaleidoscopic group} with local action $\Gamma$ is
\[ \K_k(\Gamma) = \{ g \in \Homeo(D_n) \mid \forall b \in \Br(D_n), \sigma_g(b) \in \Gamma \}, \]
which we simply denote by $\K(\Gamma)$ if the kaleidoscopic coloring $k$ is understood.
We endow $\K_k(\Gamma)$ with the permutation topology under its action on the set of branch points of $D_n$.
\end{definition}

For more details on the topic of kaleidoscopic groups, we refer to the paper \cite{Kaleidoscopic} that defines them, which also investigates their algebraic, homological and topological properties.

\subsection{Uniqueness of the airplane}

The following proposition shows that there is a unique airplane (up to homeomorphisms) and is going to be useful later for \cref{prop.airplane.embedding.is.surjective}.

\begin{proposition}
\label{prop.airplane.surjectivity}
Let $\A_1$ and $\A_2$ be two airplanes with quotient maps $\pi_1$ and $\pi_2$ to the Wa\.zewski dendrite $D_\infty$.
Let $c_1$ and $c_2$ be two kaleidoscopic colorings of $D_\infty$ satisfying \cref{eq.kaleidoscopic.coloring.order} relatively to the dense cyclic orders induced by $\pi_1$ and $\pi_2$, respectively.
Then, for every homeomorphism $\varphi$ of $D_\infty$ whose local action at each branch point is an isomorphism between the two separation relations, there is a homeomorphism $\psi$ that makes the following diagram commute.
\begin{center}
\begin{tikzcd}
    \A_1 \arrow[r, "\psi"] \arrow[d, "\pi_1"] & \A_2 \arrow[d, "\pi_2" ] \\
    D_\infty \arrow[r, "\varphi"] & D_\infty
\end{tikzcd}
\end{center}
\end{proposition}

This proof follows the general outline of that of \cref{prop.rabbit.surjectivity}.

\begin{proof}
Let $\A_1$ and $\A_2$ be two airplanes and let $Q_1$ and $Q_2$ be their quotients with respect to the circular equivalence relation, which are homeomorphic to $D_\infty$ by \cref{prop.airplane.quotient}.
Fix two kaleidoscopic colorings satisfying \cref{eq.kaleidoscopic.coloring.order}, which exist by \cref{lem.kaleidoscopic.coloring.order}.

Let $\varphi\colon Q_1\to Q_2$ be a color-preserving homeomorphism (the local action at any branch point is the identity).
The quotient maps $\pi_i\colon\A_i\to Q_i$ are injective on points that are not on a circle, so this defines a bijection $\psi=\pi_2^{-1}\circ \varphi\circ\pi_1$ between points that are not on a circle in $\A_1$ and points that are not on a circle in $\A_2$.
We wish to extend $\psi$ to a homeomorphism $\A_1 \to \A_2$.
To do this, we rely on the fact that each $\pi_i$ induces a bijection between the set of wings of $\A$ and that of branches of $Q$ (\cref{prop.airplane.quotient}).

For each circle $C$ of an airplane, each cut point $x \in C$ corresponds to a unique wing $U_x$ not including $C$.
For $C \subseteq \A_1$, the image $\pi_1(U_x)$ is a branch of $Q_1$ around the image $b_1 = \pi_1(C)$.
The homeomorphism $\varphi$ maps such a branch of $Q_1$ to a branch of $Q_2$ around $b_2=\varphi(b_1)$ which is the image of a unique wing $U_y$ associated to a circular cut point on the circle $\pi_2^{-1}(\varphi\circ\pi_1(C))$, so we can extend $\psi$ to the set of all cut points that lie on circles.
Since $\varphi$ induces a bijection that preserves the separation relation between branches around a branch point $b$ to branches around the branch point $\varphi(b)$, the map $\psi$ restricted to cut points on a given circle $C$ also preserves the separation relation and thus extends uniquely to a continuous map between $C$ and the circle that contains the images of cut points in $C$.
Doing that for all circles $C\subset \A_1$ extends $\psi$ to a bijection between $\A_1$ and $\A_2$ such that $\varphi \circ \pi_1 = \pi_2 \circ \psi$.
It now remains to show that $\psi$ is a homeomorphism.

Let $x \in \A_1$ be a point that does not lie on a circle.
Then $\psi(x)$ does not lie on a circle.
By conditions 1 and 2 of \cref{def.airplane}, the collection of wings that contain $\psi(x)$ satisfies the hypothesis of  \cref{prop.subbasis} and is thus a subbasis of neighborhoods at $\psi(x)$.
For any wing $W$ that contains $\psi(x)$, its preimage $\psi^{-1}(W)$ is a wing that contains $x$, which is open.
Thus, $\psi$ is continuous at points $x$ that do not lie on circles.

Assume that $x$ instead lies on some circle $C$ and suppose $x_n \to x$.
It suffices to consider 2 cases:
$x_n \in C$ for all $n$ or $x_n \notin C$ for all $n$.
In the first case, $\psi(x_n)\to\psi(x)$ because $\psi$ is continuous in the restriction to each circle.
In the second case, since $\varphi(\pi_1(x_n)) \to \varphi(\pi_1(x))$ we know that eventually $\psi(x_n)$ lies in any wing that contains $\psi(x)$, so no subsequence of $\psi(x_n)$ converges to a point that is not on $\psi(C)$.
For any open interval $I \subset \psi(C)$ that contains $\psi(x)$, the preimage $U_I$ of all wings attached to $\psi(C)$ at a point in $I$ is an open subset that contains $x$, so $x_n$ eventually belongs to it.
In particular, any subsequence of $\psi(x_n)$ converges to a point in $I$, which implies that $\psi(x_n)$ converges to $\psi(x)$.
Ultimately, this shows that $\psi$ is a continuous bijection from a compact space to a Hausdorff space, so it is a homeomorphism.
\end{proof}

Since by \cref{lem.kaleidoscopic.coloring.order} the hypotheses of \cref{prop.airplane.surjectivity} are met by any two airplanes, we obtain the desired result.

\begin{theorem}
\label{thm.uniqueness.airplane}
Any two airplanes are homeomorphic.
\end{theorem}

Thanks to this theorem, we can now refer to any airplane as \textit{the} airplane.

\subsection{Homeomorphism group of the airplane}
\label{sub.airplane.homeo.group.is.kaleidoscopic}

We now show that the group of homeomorphisms of an airplane is topologically isomorphic to the kaleidoscopic group $\K(\Aut(S))$, where $S$ is the separation relation of the countable dense cyclic order.

\begin{proposition}
\label{prop.homeo.airplane.acts.on.dendrite}
Given an airplane $\A$, the group $\Homeo(\A)$ acts faithfully by homeomorphisms on $Q$ in such a way that defines an embedding
\[ \Pi \colon \Homeo(\A) \hookrightarrow \K_k (\Aut(S)), \]
where $k$ is the kaleidoscopic coloring of $Q$ provided by \cref{lem.kaleidoscopic.coloring.order}
\end{proposition}

\begin{proof}
Every homeomorphism $\varphi$ of $\A$ induces a permutation of the circles of $\A$, which are the branch points of $Q$.
This permutation clearly preserves the betweenness relation, so it induces a unique homeomorphism of $Q$ (see \cite[Proposition 2.4]{Kaleidoscopic}).
Moreover, since $\varphi$ preserves the separation relation of each circle, the homeomorphism of $Q$ that it induces has actions that are locally in $\Aut(S)$, meaning that the homeomorphism belongs to $\K(\Aut(S))$.

Finally, this action of $\Homeo(\A)$ on $Q$ is faithful because each $\varphi \in \Homeo(\A)$ is completely determined by its action on the circular cut points of $\A$ (\cref{def.points.airplane}), as they densely populate $\A$ by conditions 1 and 2 of \cref{def.airplane}.
\end{proof}

The surjectivity of this embedding follows immediately from \cref{prop.airplane.surjectivity}, so we have the following fact.

\begin{proposition}
\label{prop.airplane.embedding.is.surjective}
The embedding $\Pi$ of \cref{prop.homeo.airplane.acts.on.dendrite} is surjective.
\end{proposition}

We are now only missing the bicontinuity of $\Pi$.
We endow $\Homeo(\A)$ with the uniform convergence topology.

\begin{proposition}
\label{prop.airplane.embedding.is.continuous}
The embedding $\Pi$ of \cref{prop.homeo.airplane.acts.on.dendrite} is a homeomorphism.
\end{proposition}

\begin{proof}
Both groups are Polish, so it suffices to prove that $\Pi \colon \Homeo(\A) \to \K(\Aut(S))$ is continuous (for example because of \cite[Corollary 2.3.4]{MR2455198}).
If $g_m$ is a sequence in $\Homeo(\A)$ converging to the identity, fix an arbitrary circle $C_0$ of the airplane.
We will show that, for large enough $m \in \mathbf{N}$, the homeomorphism $g_m$ fixes $C_0$.
Let $C_1$, $C_2$ and $C_3$ be circles of $\A$ that lie in three distinct wings at $C_0$, say $W_1$, $W_2$ and $W_3$, respectively.
Then, for $m$ large enough, say $m \geq N$, the circle $g_m(C_i)$ is included in $W_i$ for each $i = 1,2,3$.
Now, the circles $C_0$, $C_1$, $C_2$ and $C_3$ correspond to branch points in the dendrite of circles $Q$.
Consider the unique arcs in $Q$ between $C_1$ and $C_2$, between $C_1$ and $C_3$ and between $C_3$:
they intersect precisely at $C_0$, so the branch point $C_0$ of $Q$ is fixed by $\Pi(g_m)$ for all $m \neq N$.
This argument shows that each branch point of $Q$ is fixed by $\Pi(g_m)$ for large enough $m$, so the sequence $\Pi(g_m)$ converges to the identity of $\K(\Aut(S))$ and thus $\Pi$ is continuous.
\end{proof}

Together, the previous propositions immediately imply the desired result.

\begin{theorem}
\label{thm.homeo.airplane.is.kaleidoscopic}
The group $\Homeo(\A)$ of homeomorphisms of the airplane is topologically isomorphic to $\K(\Aut(S))$.
\end{theorem}

\subsection{Consequences of the identification with a kaleidoscopic group}

\cref{thm.homeo.airplane.is.kaleidoscopic} has the following consequences.

\begin{corollary}
\label{cor.airplane.simple}
The group $\Homeo(\A)$ is simple.
\end{corollary}

\begin{proof}
Thanks to the fact that $\Homeo(\A) = \K(\Aut(S))$, the statement follows immediately from \cite[Theorem 1.1]{Kaleidoscopic}.
\end{proof}

\begin{corollary}
\label{cor.airplane.group.closed.oligomorphic}
$\Homeo(\A)$ embeds as an oligomorphic closed subgroup of the group $\Sym(\Br(Q))$ of permutations of the circles of $\A$.
\end{corollary}

\begin{proof}
$\Aut(S)$ is oligomorphic by \cref{lem.aut.circular.order.and.relation.properties}, so $\Homeo(\A) = \K(\Aut(S))$ is oligomorphic too by \cite[Theorem 1.4]{Kaleidoscopic}.

Theorem 1.9 of \cite{Kaleidoscopic} states that $\K(N)$ is closed in $\Sym(\Br(Q))$ if and only if $N$ is closed in $\Sym([n])$.
Since $\Aut(S)$ is a closed subgroup of $\Sym([\infty])$ (\cref{rmk.Aut.S.C.closed}), we are done.
\end{proof}

Further properties of the group $\Homeo(\A)$ are explored later in \cref{sec.embedding,,sec.geometric.analytic.properties,,sec.universal.minimal.flows}.


\section{Automorphism groups of the Julia set laminations}
\label{sec.lamination.automorphisms}

In this section we study the subgroups of $\Homeo(\R_n)$ and of $\Homeo(\A)$ that preserve an orientation of the fractals induced by their laminations.
Throughout this section, $n$ will always be finite.

\subsection{Laminations of quadratic Julia sets}

In \cref{sec.julia.sets}, we prove that $n$-rabbits and the airplane are homeomorphic to the Julia set of some quadratic polynomial $f_c(z)=z^2+c$, where $c\in\mathbf{C}$ depends on the space considered.
Identifying our spaces with a Julia set $J_c$, we obtain a continuous surjective map $\varphi\colon \mathbf{S}^1\to J_c$ called the \textbf{Carathéodory loop} of $J_c$ (see for example \cite[Exposé II]{douady-hubbard-I}).
Let us fix one of the two cyclic orders on $\mathbf{S}^1$.

The Carathéodory loop defines a \textbf{laminational equivalence relation} on $\mathbf{S}^1$:
$x\sim y\iff\varphi(x)=\varphi(y)$, which enjoys the following properties:
\begin{enumerate}
    \item each class is finite;
    \item the graph of $\sim$ is closed;
    \item convex hulls in the unit open disk $\mathbf{D}$ of distinct classes are disjoint. 
\end{enumerate}

This induces a \textbf{lamination} $\mathcal{L}J_c$ on the unit disk, that is a collection of chords in the open disk $\mathbf{D}$ called \textbf{leaves} such that the intersection of two distinct leaves is empty and their union is closed in $\mathbf{D}$.
More precisely, the leaves of the lamination $\mathcal{L}J_c$ are the sides of the convex hulls of classes not reduced to a point, that are polygons (possibly reduced to a single leaf).
\cref{fig.laminations} portrays two laminations as examples.

\begin{remark}
\label{rmk.hyperbolic components}
For distinct values of $c\in\mathbf{C}$, the Julia sets of the polynomial maps $z\mapsto z^2+c$ may be homeomorphic but the laminations can be different.
For example, the Douady rabbit and its image under complex conjugation (often called the \textit{corabbit}) have different laminations.
Nonetheless, if $c$ and $c'$ belong to the same hyperbolic component of the Mandelbrot set (or are the root of this component), then the equivalence relations and thus the laminations are the same \cite[Exposé XVIII, Proposition 1]{douady-hubbard-II}.
\end{remark}

We denote by $\Aut(\mathcal{L}J_c)$ the group of \textbf{automorphisms} of the lamination $\mathcal{L}J_c$, that is the subgroup of $\Homeo^+(\mathbf{S}^1)$ that preserves the laminational equivalence relation $\sim$ (i.e., $g(x)\sim g(y)\iff x\sim y$ for all $x,y\in\mathbf{\mathbf{S}^1}$).

Let us recall that a point $c\in\mathbf{C}$ belongs to the \textbf{Mandelbrot set} if the Julia set of $f_c$ is connected, which is equivalent to the fact that the orbit of the origin is bounded.
The polynomial $f_c$ is \textbf{(sub)hyperbolic} if it is (sub)expanding on its Julia set $J_c$.
The filled-in Julia set $K_c$ is the set of points $z\in\mathbf{C}$ with bounded orbits.
The boundary $\partial K_c$ is $J_c$ and the connected components of $\interior{K_c}$ are topological disks which are in bijection with the circles of $J_c$ (their boundaries).
When $f_c$ is hyperbolic, both $K_c$ and $J_c$ are connected and locally connected.
See \cite[Exposé III]{douady-hubbard-I} for more details.

\subsubsection{Orientations and projections on circles and cut points}

Let $c\in\mathcal{M}$ such that $f_c$ is hyperbolic.
Note that, by \cref{rmk.hyperbolic components}, this does not limit our choice of Julia sets.

For each cut point $p$ of $J_c$, we define an orientation on the set of connected components of the complement $J_c \setminus \{p\}$ induced by a fixed orientation $O_{\mathbf{S}^1}$ of $\mathbf{S}^1$:
for components $X$, $Y$ and $Z$, we set $O_\varphi^p(X,Y,Z) \coloneqq +$ if the preimages $\varphi^{-1}(X)$, $\varphi^{-1}(Y)$ and $\varphi^{-1}(Z)$, which are pairwise disjoint open arcs, are positively oriented in $\mathbf{S}^1$;
else we set $O_\varphi^p(X,Y,Z) \coloneqq -$.

The following fact descends immediately from the definition of $O_\varphi^p$, where $P_p(q)$ denotes the connected component of $J_c\setminus\{p\}$ that contains $q\neq p$.

\begin{lemma}
\label{lem.cut.points.projection}
For a fixed cut point $p$ and for all $x, y, z \in \mathbf{S}^1$ with $P_p (\varphi(x))$, $P_p(\varphi(y))$ and $P_p(\varphi(z))$ distinct, the following holds:
\[ O_\varphi^p(P_p (\varphi(x)), P_p(\varphi(y)), P_p(\varphi(z))) = O_{\mathbf{S}^1}(x,y,z). \]
\end{lemma}

Similarly, for each circle $C$ of $J_c$, the orientation $O_{\mathbf{S}^1}$ of $\mathbf{S}^1$ induces a cyclic order $O_\varphi^C$ on $C$:
if $\varphi(x)$, $\varphi(y)$ and $\varphi(z)$ belong to $C$, we set $O_\varphi^C(\varphi(x),\varphi(y),\varphi(z))$ as $+$ or $-$ depending on the orientation of the points $x$, $y$ and $z$ on $\mathbf{S}^1$.
Note that each cyclic order $O_\varphi^C$ on a circle $C$ of $J_c$ corresponds to one of the two cyclic orders that define the separation relation $S$ on $C$ which we used in \cref{sub.rabbit.homeo.group.is.universal,,sub.airplane.homeo.group.is.kaleidoscopic} to determine the local actions $\Aut(S)$ of the homeomorphism groups of the rabbits and the airplane.

For each circle $C$ of $J_c$ and each point $p$ of $J_c$ there is a unique point that is the projection of $p$ onto $C$ \cite[Exposé II, Proposition 4]{douady-hubbard-I}. We denote it by $P_C(p)$.
In the next lemma we use the classical notions of \textit{external rays} and \textit{legal arcs};
see \cite{Douady-Hubbard-English} for definitions (page 10 for external rays and page 19 for legal arcs, which are called \textit{allowable arcs} in the reference).

\begin{lemma}
\label{lem.circle.projection}
For a fixed circle $C$ and for all $x, y, z \in \mathbf{S}^1$ with $P_C(\varphi(x))$, $P_C(\varphi(y))$ and $P_C(\varphi(z))$ distinct, the following holds:
\[ O_\varphi^C(P_C(\varphi(x)), P_C(\varphi(y)), P_C(\varphi(z))) = O_{\mathbf{S}^1}(x,y,z). \]
\end{lemma}

\begin{proof}
Let $E_x$, $E_y$ and $E_z$ be the external rays of $K_c$ with arguments $x$, $y$ and $z$, respectively.
Consider legal arcs $[\varphi(x),P_C(\varphi(x))]$, $[\varphi(y),P_C(\varphi(y))]$ and $[\varphi(z),P_C(\varphi(z))]$ in $K_c$.
Since $P_C(\varphi(x))$, $P_C(\varphi(y))$ and $P_C(\varphi(z))$ are distinct points of $C$, these arcs lie in distinct components of $K_C\setminus D$, where $D$ is the disk in $K_c$ such that $C=\partial D \subseteq J_c$.
So $E_x\cup[\varphi(x),P_C(\varphi(x))]$,  $E_y\cup[\varphi(y),P_C(\varphi(y))]$ and $E_z\cup[\varphi(z),P_C(\varphi(z))]$ are pairwise disjoint and the statement about the cyclic order follows from the definition of $O_\varphi^C$.
\end{proof}

Note that the projections $P_C$ map to points of $J_c$ whereas the projections $P_p$ map to subsets of $J_c$.
This choice is justified by the fact that we canonically identify each $\mathrm{out}(C)$ with a set of points of $J_c$ and, for rabbits, each $\mathrm{out}(p)$ with a set of subsets of $J_c$.

\subsubsection{Laminations of the rabbit Julia sets}

Let $T_{n,\infty}$ be the tree of circles associated to the $n$-regular rabbit $\R_n$ (see \cref{sec.tree.circles.rabbit}).
Recall that this tree is biregular:
vertices of degree $n$ correspond to cut points and vertices of infinite degree correspond to circles in $\R_n$.
Every cut point $p$ belongs to $n$ circles, each corresponding to one connected component of $\R_n$.

Let us equip the set of colors $[n]$ with its natural cyclic order, which we denote by $O_n \colon [n]^3 \to \{+,-\}$.
We denote by $\Cyc(n)$ the group of isomorphisms of this cyclic order, which is a cyclic group of order $n$.
\cref{lem.coloring.of.tree.of.circles} allows us to find a legal coloring $k$ of $T_n$ with the following properties:
\[ O_\varphi^p(X,Y,Z) = O_n(k(X), k(Y), k(Z)), \forall X, Y, Z \text{ components at a cut point } p, \]
\[ O_\varphi^C(p,q,r)=O_\Omega(k(p),k(q),k(r)), \forall p, q, r \text{ cut points on a circle } C, \]
where, in the first equation, we are using the natural identification of the half-edges at vertices corresponding to cut points with the connected components of the complements of the cut points (described in \cref{rmk.rabbit.tree.separation.relation}) and, in the second equation, the identification between half-edges at a vertex corresponding to a circle $C$ in $\R_n$ and cut points in $C$.

Using this coloring $k$, we define the subgroup
\[ \Homeo^+(\R_n) \coloneqq \Pi^{-1} (\U(\Cyc(n),\Aut(O))) \leq \Homeo(\R_n), \]
where $\Pi$ is the topological group isomorphism from \cref{sub.rabbit.homeo.group.is.universal}.

\subsubsection{Lamination of the airplane Julia set}

Recall that the arrangement of the circles of the airplane $\A$ is described by the dendrite of circles $Q$, which is homeomorphic to the universal Wa\.zewski dendrite $D_\infty$ (see \cref{sub.dendrite.of.circles}).
Each branch point of $Q$ corresponds bijectively to a circle of $\A$ and each branch of $Q$ corresponds to a wing of $\A$ (\cref{prop.airplane.quotient}).

Let us fix the quotient map $\pi \colon \A \to Q$.
Consider the isomorphism of topological groups $\Pi \colon \Homeo(\A) \to \K(\Aut(S))$ defined in \cref{sub.airplane.homeo.group.is.kaleidoscopic} and consider the subgroup $\Aut(O)$ of $\Aut(S)$.
We define the subgroup
\[ \Homeo^+(\A) \coloneqq \Pi^{-1} (\K(\Aut(O))) \leq \Homeo(\A). \]

\subsection{Identification of the lamination automorphism groups}

Given $x$, $y$ and $z$ distinct vertices of a tree or points of a dendrite, their \textbf{center} is the unique point in the intersection $[x,y] \cap [y,z] \cap [z,x]$, where $[p,q]$ denotes the unique path in the tree or arc in the dendrite that joins $p$ and $q$.

\begin{theorem}
\label{thm.auto.laminations}
Let $J$ be an $n$-rabbit or the airplane.
Endowed with the compact-open topologies, $\Homeo^+(J)$ and $\Aut(\mathcal{L}J)$ are isomorphic as Polish groups.
\end{theorem}

\begin{proof}
Let $R$ denote the preimage via $\varphi$ of the set of those points of $J$ each of which belongs to a unique circle (for rabbits, this is all regular points; for the airplane, it is the circular non-cut points).
This is a dense subset of $\mathbf{S}^1$.
Consider the action of $\Homeo^+(J)$ on $R$ defined by the equation
\[ \varphi(g \cdot x) = g(\varphi(x)), \]
for all $g \in \Homeo^+(J)$ and all $x \in R$.
We claim that this action preserves the cyclic order induced from $\mathbf{S}^1$ on its subset $R$.

To show this, let $x$, $y$ and $z$ be distinct elements of $R$.
Denote by $C_x$, $C_y$ and $C_z$ the unique circles to which $\varphi(x)$, $\varphi(y)$ and $\varphi(z)$ belong, respectively.
Note that these are vertices of the tree of circles in the case of rabbits and branch points of the dendrite of circles in the case of the airplane.
Let $v$ be the center of $C_x$, $C_y$ and $C_z$ in the tree or dendrite of circles.

Assume first that $v=C$ is a circle of $J$ (which is always the case when $J=\A$).
Then the projections $P_C(\varphi(x))$, $P_C(\varphi(y))$ and $P_C(\varphi(z))$ are distinct and we can apply \cref{lem.circle.projection} to obtain the following equalities:
\begin{align*}
O_{\mathbf{S}^1}(x,y,z) &=
O_\varphi^C(P_C\circ\varphi(x),P_C\circ\varphi(y),P_C\circ\varphi(z)) =\\
&=O_\varphi^{g(C)}(g(P_C\circ\varphi(x)),g(P_C\circ\varphi(y)),g(P_C\circ\varphi(z))) =\\
&=O_\varphi^{g(C)}(P_{g(C)}(g(\varphi(x))), P_{g(C)}(g(\varphi(y))), P_{g(C)}(g(\varphi(z)))) =\\
&=O_\varphi^{g(C)}(P_{g(C)}\circ \varphi(g\cdot x), P_{g(C)}\circ \varphi(g\cdot y), P_{g(C)}\circ \varphi(g\cdot z)) =\\
&=O_{\mathbf{S}^1}(g\cdot x, g\cdot y, g\cdot z),
\end{align*}
where the second equality holds because the local action at $C$ is $\Aut(O)$.

In the case of rabbits, $v$ can also be a cut point, say $v=p$.
Then the projections $P_p(\varphi(x))$, $P_p(\varphi(y))$ and $P_p(\varphi(z))$ must be distinct and we can apply \cref{lem.cut.points.projection} to obtain analogous equalities:
\begin{align*}
O_{\mathbf{S}^1}(x,y,z) &=
O_\varphi^p(P_p\circ\varphi(x),P_p\circ\varphi(y),P_p\circ\varphi(z)) =\\
&=O_\varphi^{g(p)}(g(P_p\circ\varphi(x)), g(P_p\circ\varphi(y)), g(P_p\circ\varphi(z))) =\\
&=O_\varphi^{g(p)}(P_{g(p)}(g\circ\varphi(x)), P_{g(p)}(g\circ\varphi(y)), P_{g(p)}(g\circ\varphi(z))) =\\
&=O_\varphi^{g(p)}(P_{g(p)}\circ \varphi(g\cdot x), P_{g(p)}\circ \varphi(g\cdot y), P_{g(p)}\circ \varphi(g\cdot z)) =\\
&=O_{\mathbf{S}^1}(g\cdot x, g\cdot y, g\cdot z),
\end{align*}
where the second equality holds because the local action at $v$ is $\Cyc(n)$.

Ultimately, the action of $\Homeo^+(J)$ preserves the cyclic order on the dense subset $R$ of $\mathbf{S}^1$.
By \cite[Lemma 3.5]{KaleidoscopicOnCircle}, it thus extends to a unique element of $\Homeo^+(\mathbf{S}^1)$.
By construction, it preserves the laminational equivalence relation $\mathcal{L}J$, so we have an embedding of $\Homeo^+(J)$ into $\Aut(\mathcal{L}J)$.
Conversely, any element of $\Homeo^+(\mathbf{S}^1)$ that preserves $\mathcal{L}J$ induces a homeomorphism of the quotient space $\mathbf{S}^1 / \mathcal{L}J = J$.
By construction of the cyclic orders $O_\varphi^p$ and $O_\varphi^C$, such homeomorphism has local actions in  $\Cyc(n)$ or $\Aut(O)$.
The fact that this correspondence is an isomorphism of topological groups follows from the fact that the Polish topologies coincide with the pointwise convergence topologies on $R$ and $\varphi(R)$.
Indeed, the compact-open topology is always stronger.
For the converse, it suffices to show it on all of $\mathbf{S}^1$.
Here, since $R$ is dense in $\mathbf{S}^1$, one can one can cover $\mathbf{S}^1$ with arbitrarily small intervals whose end points lie in $R$.
\end{proof}

\begin{remark}
By \cref{lem.aut.circular.order.and.relation.properties} the group $\Aut(O)$ is an index-$2$ normal subgroup of $\Aut(S)$.
This contrasts with the fact that $\Homeo^+(\R_n)$ (respectively $\Homeo^+(\A)$) has infinite index:
note that, for any infinite path (arc) of the tree (dendrite) of circles, for each circle lying on it, one can find an element of $\Homeo(\R_n)$ ($\Homeo(\A)$) that reverses the orientation of that circle but not of the other circles of the path (arc);
all such elements belong to distinct cosets.
And $\Homeo^+(\R_n)$ (respectively $\Homeo^+(\A)$) is not normal in $\Homeo(\R_n)$ (respectively $\Homeo(\A)$, the orientation depends on the embedding of the Julia set in $\C$.
\end{remark}

\begin{corollary}
\label{cor.lamination.groups.closed.in.circle.group}
$\Homeo^+(\R_n)$ and $\Homeo^+(\A)$ embed into $\Homeo^+(\mathbf{S}^1)$ as closed subgroups.
\end{corollary}

\begin{proof}
Let $J$ be a rabbit or airplane Julia set and $\rho \colon \Homeo^+(J) \hookrightarrow \Homeo^+(\mathbf{S}^1)$ be the topological group embedding provided by \cref{thm.auto.laminations}.
Note that the image of $\rho$ is the stabilizer of $\mathcal{L}J$ for the diagonal action of $\Homeo^+(\mathbf{S}^1)$ on $\mathbf{S}^1 \times \mathbf{S}^1$.
Since $\mathcal{L}J$ is closed in $\mathbf{S}^1\times\mathbf{S}^1$, its stabilizer is closed in $\Homeo(\mathbf{S}^1\times\mathbf{S}^1)$, so the image of $\rho$ is closed in $\Homeo(\mathbf{S}^1\times\mathbf{S}^1)$ and thus in $\Homeo^+(\mathbf{S}^1)$.
\end{proof}

Since $\Aut(O)$ is a closed subgroup of $\Sym([\infty])$ (\cref{rmk.Aut.S.C.closed}) the proof of \cref{cor.rabbit.groups.closed} applies as is to $\Homeo^+(\R_n)$, so we have the following.

\begin{corollary}
\label{cor.rabbit.lamination.group.closed}
For all $n \geq 3$, the group $\Homeo^+(\R_n)$ embeds into $\Aut(T_{n,\infty})$ as a closed subgroup.
\end{corollary}

Moreover, since $\Aut(O)$ is also oligomorphic (\cref{lem.aut.circular.order.and.relation.properties}) the proof of \cref{cor.airplane.group.closed.oligomorphic} applies to $\Homeo^+(\A)$, so we have the following.

\begin{corollary}
\label{cor.airplane.lamination.group.closed.oligomorphic}
$\Homeo^+(\A)$ embeds as an oligomorphic closed subgroup of the group $\Sym(\Br(Q))$ of permutations of the circles of $\A$.
\end{corollary}

\subsection{Simplicity of commutator subgroups}\label{sec.simplicity.commutator.subgroups}

The following statement shows that, unlike $\Homeo(\R_n)$ (\cref{cor.rabbits.simple}), the group $\Homeo^+(\R_n)$ is never simple itself, although it is virtually simple.

\begin{figure}
    \centering
    \includegraphics[width=.7\textwidth]{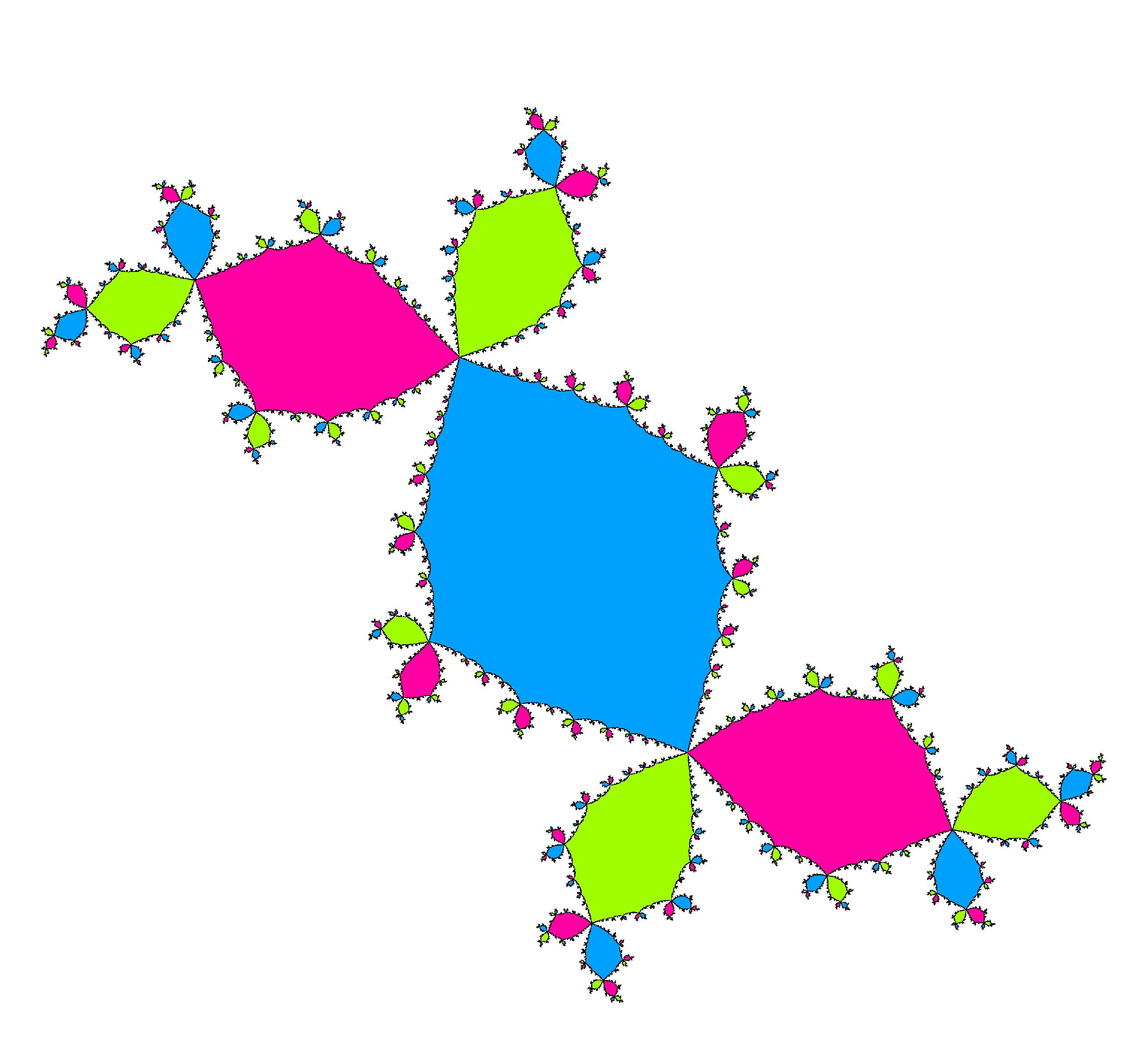}
    \caption{The $3$-coloring of the circles of $\R_3$.}
    \label{fig.julia_rabbit_colored}
\end{figure}

\begin{theorem}
\label{thm.simple.commutators}
For each $n \geq 2$, the commutator subgroup $\Homeo^+(\R_n)'$ is an index-$n$ simple subgroup of $\Homeo^+(\R_n)$ and it consists of those elements that preserve the $n$-coloring of $\R_n$ (depicted for example in \cref{fig.julia_rabbit_colored} for $n=3$).
On the other hand, $\Homeo^+(\A)$ is simple.
\end{theorem}

\begin{proof}
First note that $\Homeo^+(\A)$ is simple because every kaleidoscopic group is such, by \cite[Theorem 1.1]{Kaleidoscopic}.

For the rabbits, consider the map $f \colon \Homeo^+(\R_n) \to \Cyc(n)$ defined as $f(g) \coloneqq \sigma_g(v)$ for any vertex $v$ of degree $n$.
Of course we need to show that this map is well defined, i.e., that $\sigma_g(v) = \sigma_g(w)$ for each pair of vertices $v$ and $w$ of degree $n$ and for all $g \in \Homeo^+(\R_n)$.
It suffices to show that $\sigma_g(v) = \sigma_g(w)$ when $v$ and $w$ have distance $2$, i.e., when there is a vertex $x$ of infinite degree that is adjacent to both $v$ and $w$.
In this situation, consider the half-edges $(v,x)$ and $(w,x)$ and let $c$ denote its color, which must be the same by definition of legal colorings (\cref{def.legal.coloring}).
If $\sigma_g(v)$ maps $c$ to $c'$, then $x$ is mapped to a vertex such that $k|_{\mathrm{in}(g(x))}$ is constantly $c'$.
In particular, since $w$ is mapped to a vertex that is adjacent to $g(x)$, we have that $k(g(w),g(x)) = c'$, meaning that $\sigma_g(w)$ maps $c$ to $c'$.
Since $\sigma_g(v)$ and $\sigma_g(w)$ belong to the cyclic permutation group $\Cyc(n)$, whose action on $[n]$ is simply transitive, this is enough to conclude that $\sigma_g(v) = \sigma_g(w)$, so the map $f$ is well defined.

The map $f$ is also a group homomorphism, since $f(g \circ h) = \sigma_{g \circ h}(v) = \sigma_g(h(v)) \circ \sigma_h(v) = f(g) \circ f(h)$, and it is clearly surjective.
Let us consider the kernel $\mathrm{Ker}(f)$, which is an index-$n$ subgroup of $\Homeo^+(\R_n)$.
Note that $\mathrm{Ker}(f)$ coincides with $\U(1,\Aut(O))$ by its very definition.
Then it is a simple group by \cite[Theorem 23]{biregular}, as both the trivial group and $\Aut(O)$ are generated by their point stabilizers (\cref{cor.aut.circular.order.generated.by.stabilizers}) and $\Aut(O)$ is transitive (\cref{lem.aut.circular.order.and.relation.properties}).

Since the codomain of $f$ is abelian, we have that $\Homeo^+(\R_n)' \leq \mathrm{Ker}(f)$, so $\Homeo^+(\R_n)'$ is a normal subgroup of $\mathrm{Ker}(f)$.
Now, $\Homeo^+(\R_n)$ is clearly not abelian, so $\Homeo^+(\R_n)'$ is not trivial.
Thus $\Homeo^+(\R_n)' = \mathrm{Ker}(f)$ because $\mathrm{Ker}(f)$ is simple.
Ultimately, $\Homeo^+(\R_n)' = \mathrm{Ker}(f)$ is the index-$n$ simple subgroup of $\Homeo^+(\R_n)$ that preserves the $n$-coloring of $\R_n$, as needed.
\end{proof}

\begin{remark}
\label{rmk.rearrangement.groups}
The basilica and airplane rearrangement groups $T_B$ and $T_A$ from \cite{BF15,Belk-Forrest,airplane} are countable dense subgroups of the groups $\Homeo^+(\B)$ and $\Homeo^+(\A)$, respectively.
For the basilica this was shown in \cite{Neretin} and for the airplane a proof is sketched in \cite[Subsection 7.3]{dendriterearrangement}.

Note that the behavior of $\Homeo^+(\B)$ described in \cref{thm.simple.commutators} is shared by $T_B$, since its commutator subgroup $[T_B, T_B]$ is a simple index-$2$ normal subgroup consisting of those elements of $T_B$ that preserve the $2$-coloring depicted in \cref{fig.julia_basilica_colored} \cite[Section 8]{BF15}.
Instead, unlike $\Homeo^+(\A)$, the airplane rearrangement group $T_A$ has infinite-index simple commutator subgroup \cite[Corollary 7.4 and Thoerem 7.8]{airplane}.
\end{remark}

\begin{question}
In light of \cref{rmk.rearrangement.groups,,thm.simple.commutators,,cor.rabbits.simple}, it is worth asking the following questions.
\begin{enumerate}
    \item A rearrangement group of the basilica that does not preserve the orientation of the circles can be built by switching orientation of the edge attached to the initial vertex of the replacement graph of the basilica replacement system.
    The same arguments employed in \cite{Neretin,dendriterearrangement} show that it is dense in $\Homeo(\B)$.
    Does it have a simple index-$2$ commutator subgroup?
    \item Rearrangement groups of regular rabbits can be built as in \cite[Example 2.3]{Belk-Forrest} (they do not embed into $\Homeo^+(\R_n)$) and, as explained above for the basilica, non-orientation-preserving variants groups can be defined too.
    Is any of these groups simple?
    Are they dense in $\Homeo(\R_n)$?
    \item A non-orientation-preserving variant of $T_A$ can be built too, using the same modifications.
    Is it simple?
    Is it dense in $\Homeo(\A)$?
\end{enumerate}
\end{question}


\section{Embedding and non-embedding results}
\label{sec.embedding}

In this section we show a general embedding result about universal groups of regular trees $T_n$ and kaleidoscopic groups of Wa\.zewski dendrites $D_n$, which applies to $\Homeo(\B)$ and $\Homeo(\A)$.
Then we prove further embedding and non-embedding results about the homeomorphism and lamination automorphism groups of rabbits and the airplane.

\subsection{Patchwork}

Patching together partially defined homomorphisms is useful for both universal and kaleidoscopic groups.
Such tools were developed in \cite[Section 4]{Kaleidoscopic} for dendrites and partially in \cite[Lemma 17]{biregular} for biregular trees.
Here we develop the relevant results for biregular trees and we recall them from \cite{Kaleidoscopic} for dendrites.

Since most of the statements are essentially the same in both settings, we will let $T$ be a biregular tree and $D$ be a Wa\.zewski dendrite.
Given $p$ and $q$ vertices of a tree or points of a dendrite, we will denote by $[p,q]$ either the unique path in the tree (as a set of vertices) or the unique arc in the dendrite (as a subset of the dendrite) between $p$ and $q$.
Recall that, given $x$, $y$ and $z$ distinct vertices of $T$ or points of $D$, their \textbf{center} is the unique point in the intersection $[x,y] \cap [y,z] \cap [z,x]$.

\begin{definition}
\label{def.center.closed}
A subset $F$ of $X$ is \textbf{center-closed} if it contains the center of any of three distinct points in it.
\end{definition}

If $F$ is a non-empty subset of $V(T)$, we denote by $[F]$ the minimal subtree of $T$ whose set of vertices includes $F$.

\begin{definition}
\label{def.partial.morphism}
We define two notions of partial homomorphism.
\begin{itemize}
    \item Given subsets $F$ and $F'$ of $V(T)$, a \textbf{partial tree homomorphism} $F \to F'$ is a bijection that is a restriction of a tree isomorphism $[F] \to [F']$.
    \item Given center-closed subsets $F$ and $F'$ of $D$, a \textbf{partial dendrite homomorphism} $F \to F'$ is a bijection that preserves the betweenness relation and maps branch points to branch points, end points to end points and regular points to regular points.
\end{itemize} 
\end{definition}

It is easy to see that any partial tree homomorphism $F \to F'$ is the restriction of a unique tree isomorphism $[F] \to [F']$.

If $F$ is a subset of $V(T)$ or $D$, we denote by $\Omega_F$ the set of connected components of $T \setminus F$ or $D \setminus F$.
Note that of each element of $\Omega_F$ is either a subtree of $T$ or its closure is a subdendrite of $D$.

\begin{remark}
\label{rmk.boundary.omega.F}
If $F$ is a center-closed subset of a dendrite $D$, then each element of $\Omega_F$ has at most two points of $F$ on its topological boundary.
Indeed, if $\partial A \supseteq \{x,y,z\}$ for some $A \in \Omega_F$, then the center of $x$, $y$ and $z$ would belong to $F$ and so $A$ could not be a connected component of $D \setminus F$.

In the case of a tree $T$, if $A \in \Omega_F$ then we denote by $\partial A$ the set of vertices of $T$ that are at distance $1$ from $A$.
By the same argument as for dendrites, if $F$ is center-closed then each $\partial A$ consists either of one or two vertices of $F$.
\end{remark}

The following lemma is our main tool for constructing elements of $\U(N,M)$ and $\K(N)$ from partially defined homomorphisms.
It combines \cite[Lemma 4.3]{Kaleidoscopic} (which concerns dendrites) with a useful yet simple tool for trees.

\begin{lemma}
\label{lem.patchwork}
Suppose that $F$ and $F'$ are subsets of $V(T)$ or $D$ and let $f \colon F \to F'$ be a partial homomorphism.
In the case of a dendrite, suppose also that $F$ and $F'$ are closed, center-closed and that, for each two cut points $x$ and $y$ of $D$, $[x,y] \cap F$ and $[x,y] \cap F'$ are finite.
Assume the following:
\begin{enumerate}
    \item for each $A \in \Omega_F$ there is a tree isomorphism or a homeomorphism $h_A \colon A \cup \partial A \to B \cup \partial B$ such that $B \in \Omega_{F'}$ and $h_A|_{\partial A} = f|_{\partial A}$;
    \item every $B \in \Omega_{F'}$ is $h_A(A)$ for some unique $A \in \Omega_F$.
\end{enumerate}
Define a map $h \colon V(T) \to V(T)$ or $D \to D$ by setting $h(x) = h_A(x)$ when $x \in A$ for some $A \in \Omega_F$ and $h(x) = f(x)$ for all $x \in F$.
Then $h$ is an automorphism of $T$ or a homeomorphism of $D$ that agrees with $f$ on $F$.
\end{lemma}

Note that the additional hypotheses required on $F$ and $F'$ in the case of a dendrite always hold when such sets are finite and center-closed.

\begin{proof}
For dendrites, this is Lemma 4.3 from \cite{Kaleidoscopic} when $F$ and $F'$ are finite.
When they are infinite,
we will show that $h$ is continuous by proving that it maps connected components of the complement of a point to connected components of the complement of a point, thus preserving a subbasis for the topology of the dendrite.
This is equivalent to showing that $h$ maps arcs to arcs.
Observe that each $A \in \Omega_F$ or $B \in \Omega_{F'}$ has non-empty topological boundary.
Indeed, $A$ is open because $D$ is locally connected and $F$ is closed, so if $\partial A$ were empty then $A$ would be clopen, which contradicts the fact that $D$ is connected (and the same for any $B$ in $\Omega_{F'}$).
Then $\partial A$ (respectively $\partial B$) consists exactly of one or two points of $F$ (respectively $F'$) by \cref{rmk.boundary.omega.F}.
Now, each arc $[x,y]$ of $D$ with $x$ and $y$ cut points splits into finitely many closed subarcs with end points in $F \cup \{x,y\}$, so $h$ maps it to an arc, because it is a homeomorphism $h_A$ on the interior of each subarc and because two adjacent subarcs meet in a point of $F$ on which $h$ has the same value as $f$.
If $x$ is not a cut point but $y$ is (the case when neither are cut points is similar), then $(x,y]$ is an increasing union of arcs $[x_i,y_i]$ whose end points are cut points, so the previous argument applies to each $[x_i,y_i]$, thus on $(x,y]$, and ultimately to $[x,y]$ because if $x \in F$ then $h(x)=f(x)$ and if $x \notin F$ then $h(x)=h_A(x)$ for some $A \in \Omega_F$.

Let us now consider the case of a tree.
In this case, first note that $h$ is a well-defined bijection because $\Omega_F \cup \{F\}$ is a partition of $V(T)$.
Then we only need to check that the bijection $h$ of the vertices of $T$ is truly a tree automorphism, i.e., that $\{v,w\} \in E(T)$ if and only if $\{h(v),h(w)\} \in E(T)$.
When both $v$ and $w$ belong to $F$, this follows from the fact that $h|_F=f|_F$ and that $f$ descends from an isomorphism $[F] \to [F']$.
On the other hand, when either or both vertices belong to the complement of $F$, since they are at distance $1$ they must belong to a unique common $A \cup \partial A$, so the conclusion follows from the fact that each $h_A$ is a tree automorphism defined on $A \cup F$.
\end{proof}

\begin{remark}
\label{rmk.patchwork.coloring}
The map $h$ produced by \cref{lem.patchwork} belongs to $\U(N,M)$ or $\K(N)$ as soon as the local actions at the elements of $F$ induced by $f$ together with the correspondence $\Omega_F \to \Omega_{F'}$ belong to $N \sqcup M$ or $N$ and if the maps $h_A$ have local actions in $N \sqcup M$ or $N$.
\end{remark}

For dendrites, the maps $h_A$ needed in \cref{lem.patchwork} can be often built using \cite[Corollary 4.5]{Kaleidoscopic}, which we recall right below.
Recall that $c_x(y)$ denotes the unique connected component at $x$ that contains $y$.
Given two points $x$ and $y$, we denote by $C_{x,y}$ the intersection between $c_x(y)$ and $c_y(x)$.

\begin{proposition}[Corollary 4.5 from \cite{Kaleidoscopic}]
\label{prop.partial.dendrite.homeomorphism}
Consider two distinct pairs of points $a_0, a_1$ and $b_0, b_1$ of $D$.
Then there is a homeomorphism
\[ g \colon C_{a_0,a_1} \cup \{a_0,a_1\} \to C_{b_0,b_1} \cup \{b_0,b_1\} \]
such that $g(a_i)=b_i$ for $i=1,2$ and that, for all $p \in C_{a_0,a_1} \cap \Br$ and $q \in C_{a_0,a_1} \cup \{a_0,a_1\}$ distinct, $k_p(q) = k_{g(p)}(g(q))$ (i.e., the local action of $g$ at $p$ is trivial).
\end{proposition}

\subsection{Embedding trees into dendrites}

As done in \cite{biregular}, given a half-edge $a=(v,w)$, we denote by $\overline{a}$ the opposite half-edge $(w,v)$.
Moreover, denote by $HT(T_n)$ the set of \textbf{half-trees} of $T_n$, which are the connected components obtained by the removal of an edge.
Every edge $e = \{v,w\}$ identifies a half-tree for each of each two half-edges:
if $a=(v,w)$, we denote by $T_a$ the half-tree that contains the origin $v$ of $a$ among the two half-trees obtained by removing the edge $\{v,w\}$ (thus, $T_{\overline{a}}$ is the half-tree that contains $w$ instead).
Let $\tau$ denote the bijection $H(T_n) \to HT(T_n),\, a \mapsto T_{\overline{a}}$.

We denote by $\overline{T_n}$ the union of the geometric realization of $T_n$ and the boundary $\partial T_n$ (which is the set of semi-infinite paths from an arbitrary fixed vertex).
We equip $\overline{T_n}$ with the observer's topology, which is the topology generated by the subbasis consisting of all half-trees.
Note that $\overline{T_n}$ is a dendrite (though not a Wa\.zewski dendrite).
When $n$ is finite, $\overline{T_n}$ coincides with the end compactification of $T_n$.

A topological embedding $\Phi \colon \overline{T_n} \to D_n$ induces an injective map $\hat{\Phi} \colon V(T_n) \to \Br(D_n)$.
Let $N \leq \Sym([n])$ and assume that $T_n$ and $D_n$ are equipped with a legal coloring $k_T$ and a kaleidoscopic coloring $k_D$, respectively, both with $[n]$ as the set of colors.
If $S \subseteq D_n$ is fully included in some connected component at $b \in \Br(D_n)$, we denote that connected component by $c_b(S)$.

\begin{definition}
\label{def.color.preserving.embedding}
A topological embedding $\Phi \colon \overline{T_n} \to D_n$ is \textbf{color-preserving} if, for all $v \in V(T_n)$ and all $a \in \mathrm{out}(v)$,
\[ k_T (a) = k_D \circ c_{\Phi(v)} \circ \Phi \circ \tau (a), \]
i.e., it maps each half-tree of $T_n$ inside a branch of $D_n$ of the same color.
\end{definition}

\begin{lemma}
\label{lem.exist.nice.embeddings}
For every $n \geq 3$, there exist color-preserving topological embeddings $\Phi \colon \overline{T_n} \to D_n$ such that the closure of $\bigcup_{\{v,w\}\in E(T_n)} \overline{C_{\Phi(v),\Phi(w)}}$ is $D_n$.
\end{lemma}

\begin{proof}
Throughout this proof, let $\partial D_n$ denote the set of end points of $D_n$.
Denote by $B_{v_0}(m)$ the radius-$m$ ball of $T_n$ centered in an arbitrary vertex $v_0$.
We start by letting $\Phi_0$ be any map $B_{v_0}(0) = \{v_0\} \to \Br(D_n)$.
Now suppose that we are given an injective map $\Phi_m \colon B_{v_0}(m) \to \Br(D_n)$ such that
\begin{enumerate}
    \item $u \in [v,w]$ in $T_n$ if and only if $\Phi_m(u) \in [\Phi_m(v),\Phi_m(w)]$ in $D_n$;
    \item for all adjacent $v,w \in B_{v_0}(m)$, one has $k_T (v,w) = k_D \left( c_{\Phi_m(v)} \left( \Phi_m \left(w\right) \right) \right)$.
\end{enumerate}
The map $\Phi_0$ trivially satisfies these conditions.
Now let us show that we can extend $\Phi_m$ to some $\Phi_{m+1} \colon B_{v_0}(m+1) \to \Br(D_n)$ with these properties.

For all $v \in B_{v_0}(m+1) \setminus B_{v_0}(m)$, let $w$ be the unique vertex in $B_{v_0}(m)$ that is adjacent to it.
Consider the branches at $\Phi_m(w)$ and denote by $B$ the only one whose color is $k_T(w,v)$.
We will choose $\Phi_{m+1}(v)$ to be a branch point included in $B$:
with these choices we already know that $\Phi_{m+1}$ satisfies condition (1) and that condition (2) holds for half-edges originating from vertices in $B_{v_0}(m)$.
By definition of kaleidoscopic coloring (\cref{def.kaleidoscopic.coloring}), there exist branch points in $B$ such that $k_D\left(c_b(\Phi_m(w))\right) = k_T(v,w)$.
Let $\Phi_{m+1}(v)$ be any such point.
Then condition (2) is satisfied for half-edges originating from vertices in $B_{v_0}(m+1)$ too.
Doing this for all $v \in B_{v_0}(m+1) \setminus B_{v_0}(m)$ produces a map $\Phi_{m+1} \colon B_{v_0}(m+1) \to \Br(D_n)$ that extends $\Phi_m$ and satisfies conditions (1) and (2).

Let $\Phi_* \colon V(T_n) \to \Br(D_n)$ be defined as $\Phi_m(v)$ for $m$ such that $v \in B_{v_0}(m)$.
For each edge $\{v,w\}$ of $T_n$, extend $\Phi_*$ on the corresponding arc $[v,w]$ of $\overline{T_n}$ to be any homeomorphism from $[v,w]$ to the unique arc $[\Phi_*(v),\Phi_*(w)]$ of $D_n$.
Thanks to condition (1) on each $\Phi_m$, this defines an embedding $\overline{T_n} \setminus \partial (T_n) \to D_n \setminus \partial D_n$, which thus extends to an embedding $\Phi \colon \overline{T_n} \to D_n$.
Thanks to condition (2) on the maps $\Phi_m$, the embedding $\Phi$ is color-preserving in the sense of \cref{def.color.preserving.embedding}.

Finally, consider the set $\Omega_{\Phi(V(T_n))}$ of connected components of $D_n \setminus \Phi(V(T_n))$.
For each edge $\{v,w\}$ of $T_n$, let us denote by $C_{v,w}$ the unique element of $\Omega_{\Phi(V(T_n))}$ whose boundary contains $\Phi(v)$ and $\Phi(w)$ (which is $C_{\Phi(v),\Phi(w)}=c_{\Phi(v)}(\Phi(w)) \cap c_{\Phi(w)}(\Phi(v))$ by the notation established right before \cref{prop.partial.dendrite.homeomorphism}).
Observe that the subspace
\[ D_n' \coloneqq \bigcup_{\{v,w\}\in E(T_n)} \left(C_{v,w}\cup\{\Phi(v),\Phi(w)\}\right) \cup \Phi(\partial T_n) \subseteq D_n \]
is a subdendrite of $D_n$ that is homeomorphic to $D_n$ itself.
This follows from \cite[Lemma 2.14]{DendriteGroups}, together with the fact that this set is the closure of the open connected subset $\bigcup_{\{v,w\}\in E(T_n)} \left(C_{v,w}\cup\{\Phi(v),\Phi(w)\}\right)$.
Let us equip it with the coloring induced from that of $D_n$, which is readily seen to be a kaleidoscopic coloring of $D_n'$.
Note that $D_n'$ includes the image of $\Phi$ and that $\Phi$ maps the end points of $\overline{T_n}$ to end points of $D_n'$.
Then, if we replace the codomain $D_n$ of $\Phi$ with $D_n'$, we have that $\Phi$ is color-preserving and satisfies the property in the statement, as needed.
\end{proof}

\subsection{Embedding universal groups into kaleidoscopic groups}

We now prove a general fact that relates universal groups to kaleidoscopic groups.

\begin{theorem}
\label{thm.universal.embeds.into.kaleidoscopic}
A color-preserving topological embedding $\Phi \colon \overline{T_n} \hookrightarrow D_n$ such that the closure of $\bigcup_{\{v,w\}\in E(T_n)} \overline{C_{\Phi(v),\Phi(w)}}$ is $D_n$ induces a topological group embedding $\varphi \colon \U(N) \to \K(N)$ such that $\Phi \circ g |_{V(T_n)} = \varphi(g) \circ \Phi|_{V(T_n)}$ for all $g \in \U(N)$.
\end{theorem}

\begin{proof}
Given $g \in \U(N)$, we want to build $\varphi(g) \in \K(N)$.
We start by considering $f \coloneqq \Phi \circ g \circ \Phi^{-1}$ on $\Phi \left( V(T_n) \cup \partial T_n \right)$.
Note that $f$ is a partial dendrite homomorphism (\cref{def.partial.morphism}).
We will now extend $f$ ``canonically'' everywhere else using \cref{lem.patchwork} with $F=F'=\Phi \left( V(T_n) \cup \partial T_n \right)$, as explained below.

By hypothesis, the closure of $\bigcup_{\{v,w\}\in E(T_n)} \overline{C_{\Phi(v),\Phi(w)}}$ is the entirety of $D_n$, so $\Phi$ must map end points of $\overline{T_n}$ to end points of $D_n$.
Thus the set of accumulation points of $\Phi \left( V(T_n) \cup \partial T_n \right)$ is precisely $\Phi(\partial T_n)$ and so $[x,y] \cap \Phi\left( V(T_n) \cup \partial T_n \right)$ is finite for all cut points $x$ and $y$ of $D_n$.
Moreover, a connected component of $D_n\setminus\Phi \left( V(T_n) \cup \partial T_n \right)$ coincides with a connected component of $D_n\setminus\Phi \left(V(T_n)\right)$.
Finally, the set $\Phi\left( V(T_n) \cup \partial T_n \right)$ is closed and center-closed, so it satisfies all the hypotheses of \cref{lem.patchwork}.

Consider the set $\Omega_{\Phi(V(T_n))}$ of connected components of $D_n \setminus \Phi(V(T_n))$, which we will denote by $\Omega$ for the sake of brevity.
Recall that the boundary of each element of $\Omega$ in $D_n$ consists of at most two elements of $\Phi(V(T_n))$ by \cref{rmk.boundary.omega.F}.
Actually, since $\Phi$ is such that the closure of $\bigcup_{\{v,w\}\in E(T_n)} \overline{C_{\Phi(v),\Phi(w)}}$ is $D_n$, the boundary of each element $A$ of $\Omega$ must have two points which belong to $\Phi(V(T_n))$, otherwise $A$ would be left out of such set.
For the sake of brevity, given an edge $\{v,w\}$ of $T_n$, let us write $C_{v,w}$ for $C_{\Phi(v),\Phi(w)}$, which is the element of $\Omega$ whose boundary is $\{\Phi(v),\Phi(w)\}$.
Let $D_{v,w}$ denote the closure of $C_{v,w}$, which is $C_{v,w} \cup \{\Phi(v),\Phi(w)\}$.
Since $D_{v,w}$ is homeomorphic to $D_n$, for all $\{v,w\} \in E(T_n)$ we can apply \cref{prop.partial.dendrite.homeomorphism} to find color-preserving homeomorphisms $\psi_{(v,w)} \colon D_{v,w} \to D_n$ such that
\[ \psi_{(v_1,w_1)}(\Phi(v_1)) = \psi_{(v_2,w_2)}(\Phi(v_2)) \text{ and } \psi_{(v_1,w_1)}(\Phi(w_1)) = \psi_{(v_2,w_2)}(\Phi(w_2)) \]
for every pair of edges $\{v_1,w_1\}$ and $\{v_2,w_2\}$ of $T_n$.
We do this for both half-edges of $\{v,w\} \in E(T_n)$ in order to find both $\psi_{(v,w)}$ and $\psi_{(w,v)}$.

For all $A \in \Omega$, say $A=C_{v,w}$ for some $\{v,w\} \in E(T_n)$, we let $h_A \coloneqq \psi_{\left(g(v),g(w)\right)}^{-1} \circ \psi_{(v,w)} \colon D_{v,w} \to D_{g(v),g(w)}$.
With these maps $h_A$, we can apply \cref{lem.patchwork}, which produces an element $h \in \Homeo(D_n)$.
Let us now show that each $h$ actually belongs to $\K(N)$.

First note that, by construction, the local action at each branch point outside $\Phi(V(T_n))$ belongs to $N$.
Then, thanks to \cref{rmk.patchwork.coloring}, we only need to show that the local action of $h$ at $\Phi(v)$ belongs to $N$ for each $v \in V(T_n)$.
To do so we will prove that $\sigma_h\left(\Phi(v)\right) = \sigma_g(v)$.
Let $b_0 \coloneqq \Phi(v)$ and $b_1 \coloneqq \Phi(g(v)) = h(b_0)$, which are branch points of $D_n$.
For $A \in \Omega$, say $A = C_{v,w}$, the corresponding color $k_D\left(c_{b_0}(A)\right)$ is mapped to
\begin{gather*}
k_D \circ c_{b_1} \circ h_A \circ c_{b_0}^{-1} \circ (k_D|_{\widehat{b_0}})^{-1} \left(k_D\left(c_{b_0}(A)\right)\right) =
\\
k_D \circ c_{b_1} \circ h_{C_{v,w}} (C_{v,w}) = k_D \left( c_{b_1} \left( C_{g(v),g(w)} \right) \right).
\end{gather*}
Since $\Phi$ is color-preserving (\cref{def.color.preserving.embedding}) and the restriction $c_{b_0}|_{\Phi(HT(T_n))}$ is a bijection $\Phi(HT(T_n)) \to \widehat{b_0}$, we have
\begin{gather*}
k_T = k_D \circ c_{b_1} \circ \Phi \circ \tau,
\\
\left(k_T|_{\mathrm{out}(v)}\right)^{-1} = \tau^{-1} \circ \Phi^{-1} \circ \left(c_{b_0}|_{\Phi(HT(T_n)))}\right)^{-1} \circ (k_D|_{\widehat{b_0}})^{-1}.
\end{gather*}
Using these two equations, we find that
\begin{gather*}
\sigma_g(v) \left(k_D\left(c_{b_0}(A)\right)\right) = k_T \circ g \circ \left(k_T|_{\mathrm{out}(v)}\right)^{-1} \left(k_D\left(c_{b_0}(A)\right)\right) =
\\
k_D \circ c_{b_1} \circ \Phi \circ \tau \circ g \circ \tau^{-1} \circ \Phi^{-1} \circ \left(c_{b_0}|_{\Phi(HT(T_n)))}\right)^{-1} \circ (k_D|_{\widehat{b_0}})^{-1} \left(k_D\left(c_{b_0}(A)\right)\right) =
\\
k_D \circ c_{b_1} \circ \Phi \circ \tau \circ g \circ \tau^{-1} \circ \Phi^{-1} \circ \left(c_{b_0}|_{\Phi(HT(T_n)))}\right)^{-1} \circ c_{b_0} (C_{v,w}) =
\\
k_D \circ c_{b_1} \circ \Phi \circ \tau \circ g \circ \tau^{-1} \circ \Phi^{-1} (\Phi(T_{\overline{(v,w)}})) = k_D \circ c_{b_1} \circ \Phi \circ \tau \circ g (v,w) =
\\k_D \circ c_{b_1} \circ \Phi (T_{\overline{(g(v),g(w))}}) = k_D \left( c_{b_1} \left(C_{g(v),g(w)}\right) \right).
\end{gather*}
Thus, ultimately all of the local actions of $h$ belong to $N$, so $h \in \K(N)$.
We let $\varphi(g) = h$ for $h$ found in this manner, which defines a map $\varphi \colon \U(N) \to \K(N)$.
It remains to show that it is a topological group embedding.

Consider two elements $g_1, g_2 \in \U(N)$.
For any $v \in V(T_n)$, it is clear that $\varphi(g_1) \circ \varphi(g_2) (\Phi(v)) = \varphi(g_1 \circ g_2) (\Phi(v))$.
The set $\Omega$ is a partition of $D_n \setminus \Phi(V(T_n))$, so it suffices to see that $\varphi(g_1) \circ \varphi(g_2) = \varphi(g_1 \circ g_2)$ on every $C_{v,w} \in \Omega$.
This is an immediate computation using the fact that $\varphi(g)|_{D_{v,w}} = \psi_{g(v),g(w)}^{-1} \circ \psi_{v,w}$, so $\varphi$ is a group homomorphism.

Now, $\varphi$ is clearly injective, as $\varphi(g)$ being trivial implies that it fixes each element of $\Phi(V(T_n))$, so any $g \in \mathrm{Ker}(\varphi)$ must fix every vertex of $T_n$.

Continuity is readily proved as follows.
If $(g_m)_{m \in \mathbf{N}}$ is a convergent sequence of $\U(N)$ then for all $v \in V(T_n)$ there is a $K(v) \in \mathbf{N}$ such that $g_i(v)=g_j(v)$ for all $i,j \geq K(v)$.
We want to show that the same holds for $\left(\varphi(g_m)\right)_{m \in \mathbf{N}}$ for all branch points.
If $b \in \Br(D_n) \setminus \Phi(V(T_n))$, then $b \in D_{v,w}$ for a unique edge $\{v,w\} \in E(T_n)$.
For all $i,j \geq \mathrm{max}(K(v),K(w))$, the elements $g_i$ and $g_j$ agree on $\{v,w\}$, so $\varphi(g_i)|_{D_{v,w}} = \varphi(g_j)|_{D_{v,w}}$ by construction.
The fact that $\varphi(g)|_{\Phi(V(T_n))} = g \circ \Phi$ immediately shows that the same holds for every $b \in \Phi(V(T_n))$.
\end{proof}

\begin{remark}
It is likely that statements analogous to \cref{lem.exist.nice.embeddings,,thm.universal.embeds.into.kaleidoscopic,,cor.universal.embeds.into.kaleidoscopic} also hold for universal groups $\U(N,M)$ on biregular trees $T_{n,m}$ and generalized kaleidoscopic groups $\K(N,M)$ on Wa\.zewski dendrites $D_{n,m}$.
However, the groups $\K(N,M)$ have not been formally defined in the literature.
They can be built by using two disjoint sets of colors $[n]$ and $[m]$ and fixing two local groups $N$ and $M$ acting on $[n]$ and $[m]$, respectively, as noted in \cite[Remark 3.10]{Kaleidoscopic}.
\end{remark}

Our main application of \cref{lem.exist.nice.embeddings,,thm.universal.embeds.into.kaleidoscopic} is the following fact, which we will apply to $\Homeo(\B)$ and $\Homeo(\A)$ in the next subsection.

\begin{corollary}
\label{cor.universal.embeds.into.kaleidoscopic}
For every $n \geq 3$ at most countably infinite and for all $N \leq \Sym([n])$, there are topological group embeddings
\[ \U(N) \hookrightarrow \K(N). \]
\end{corollary}

It is worth mentioning that, when we consider the case $N=\Sym([n])$, \cref{cor.universal.embeds.into.kaleidoscopic} yields the following fact.

\begin{corollary}
For every $n \geq 3$ at most countably infinite, there are topological group embeddings
\[ \Aut(T_n) \hookrightarrow \Homeo(D_n). \]
\end{corollary}

\begin{remark}
For any tree $T$ whose degrees are at most countable, if $n$ is the supremum of its degrees, then $T$ embeds into $T_n$.
Arguments similar to those that we will soon use to prove \cref{lem.embeddings} show that this induces an embedding $\Aut(T) \hookrightarrow \Aut(T_n)$, so ultimately the previous corollary also shows that $\Aut(T)$ embeds into $\Homeo(D_n)$.
\end{remark}

\subsection{Embeddings and non-embeddings among the homeomorphism groups}

Recall that universal and kaleidoscopic groups (\cref{def.universal.group,,def.kaleidoscopic.group}) are defined by permutation groups $N$ (and $M$).
Before considering the homeomorphism groups of the rabbits and the airplane, let us show that embeddings of permutation groups induce embeddings of universal and kaleidoscopic groups.

\begin{lemma}
\label{lem.embeddings}
Consider two permutation groups $N_1$ and $N_2$ of $K_N^{(1)}$ and $K_N^{(2)}$, respectively.
Assume that there exist an injective map $\chi \colon K_N^{(1)} \hookrightarrow K_N^{(2)}$ and a group embedding $\iota \colon N_1 \hookrightarrow N_2$ such that
\[ \iota(g) \cdot x = g \cdot \chi(x) \]
for all $x \in K_N^{(1)}$ and $g \in N_1$.
Then $\U(N_1)$ embeds into $\U(N_2)$ and $\K(N_1)$ embeds into $\K(N_2)$, both as topological groups.

Moreover, if we have the same situation for two additional permutation groups $M_1$ and $M_2$ of $K_M^{(1)}$ and $K_M^{(2)}$, respectively, then $\U(N_1, M_1)$ embeds into $\U(N_2, M_2)$ as a topological group.
\end{lemma}

\begin{proof}
We will prove this for regular trees first and then for Wa\.zewski dendrites, since the case of biregular trees is essentially identical to that of regular trees, only notationally heavier.

Let $n_1 = |K_N^{(1)}|$ and $n_2 = |K_N^{(2)}|$ and denote by $T_{n_1}$ and $T_{n_2}$ the $n_1$-regular tree and the $n_2$-regular tree, respectively.
Assuming the hypotheses of the statement, $\chi$ induces embeddings $\Chi \colon T_{n_1} \hookrightarrow T_{n_2}$ that is color-preserving, i.e., $k_2(v,w) = k_1(\Chi(v),\Chi(w))$ for all $v,w\in V(T_{n_1})$, where $k_i$'s are the colorings of $T_{n_1}$ and $T_{n_2}$, respectively.
Each such embedding is completely determined by the image of a sole vertex.
Let us fix one such embedding $\Chi$.
The group $\U(N_1)$ acts faithfully on $T_{n_2}$ as follows.
On $\Chi(T_{n_1})$, the action is $g \cdot \Chi(x) \coloneqq \Chi(g \cdot x)$ (where the second action is that of $\U(N_1)$ on $T_{n_1}$).
If $A$ is a connected component of $T_{n_2} \setminus \Chi(T_{n_1})$ that is adjacent to $v \in \Chi(T_{n_1})$, then $g$ maps it to the unique component adjacent to $g \cdot v$ with color $\sigma_{\iota(g)}(v)\left(k_2(\tau_2^{-1}(A))\right)$, where $\tau_2$ is the bijection $H(T_{n_2}) \to HT(T_{n_2})$.
The action of $g$ on $A$ is the unique color-preserving action.
This truly defines an element of $\U(N_2)$ by \cref{lem.patchwork,rmk.patchwork.coloring}, so it determines an embedding $\U(N_1) \hookrightarrow \U(N_2)$.

In the case of dendrites and kaleidoscopic groups, let us see that again $\chi$ induces color-preserving embeddings $\Chi \colon D_{n_1} \hookrightarrow D_{n_2}$.
We will first produce an injective map $\tau\colon\Br(D_{n_1})\to\Br(D_{n_2})$ that preserves the betweenness relation and maps any branch $B$ of $D_{n_1}$ to a branch of $D_{n_2}$ whose color is the image via $\chi$ of the color of $B$.
Let $F_n$ be an increasing exhausting sequence of finite subsets of $\Br(D_{n_1})$ such that $F_0$ contains just a point $b_0$ and $F_{n+1}\setminus F_n$ is a point in $[F_n]$ or a point whose projection is in $[F_n]$.
The map $\tau$ is defined inductively as follows.
Let $\tau(b_0)$ be any branch point of $D_{n_2}$.
Assume that $\tau$ has been defined on $F_n$, is a partial dendrite homomorphism (\cref{def.partial.morphism}) and, for any $b \neq b'$ in $F_n$, is such that $k_2(c_{\tau(b)}(\tau(b'))) = \chi(k_1(c_b(b')))$.
Now let $b\in F_{n+1}\setminus F_n$ and let $p$ be its projection on $[F_n]$.
If $p\in F_n$ then we choose a point $b_2$ in the branch  $B$ of $D_{n_2}$ at $\tau(p)$ such that $k_2(B)=\chi(k_1(c_p(b)))$ and $k_2(c_{b_2}(\tau(p)))=\chi(k_1(c_b(p))$ and we define $\tau(b)=b_2$.
If $p\notin F_n$ then $p=b\in [F_n]$ and there are $b_1,b_1'\in F_n$ such that $b\in [b_1,b_1']$ and the interior of $[b_1,b_1']$ contains no points of $F_n$.
By the kaleidoscopic property, there is $b_2\in [\tau(b_1), \tau(b'_1)]$ such that $k_2(c_{b_2}(\tau(b_1)))=\chi(k_1(c_b(b_1)))$ and $k_2(c_{b_2}(\tau(b'_1)))=\chi(k_1(c_b(b'_1)))$.
We set $\tau(b)=b_2$.
This defines $\tau$ on $\Br(D_{n_1})$ with the required properties, which thus extends continuously to an embedding $\Chi\colon D_{n_1}\to D_{n_2}$ by \cite[Proposition 2.4]{Kaleidoscopic}.

Now, the action of $\K(N_1)$ on $D_{n_2}$ is defined as before on $\Chi(D_{n_1})$ (i.e., $g \cdot \Chi(x) \coloneqq \Chi(g \cdot x)$) and, if $B$ is a connected component of $D_{n_2} \setminus \Chi(D_{n_1})$ at the branch point $b \in \Chi(D_{n_1})$, then $g$ maps it color-preservingly to the unique component $B'$ adjacent to $g \cdot b$ with color $\sigma_{\iota(g)}(b)\left(k_2(B)\right)$.
However, unlike in the case of trees, for Wa\.zewski dendrites there are infinitely many color-preserving homeomorphisms $B \to B'$.
As we did in the proof of \cref{thm.universal.embeds.into.kaleidoscopic}, we can make a canonical choice using \cref{prop.partial.dendrite.homeomorphism}:
for each connected component $B$ of $D_{n_2} \setminus \Chi(D_{n_1})$ with $\{b\}$ as its boundary, we fix a color-preserving homeomorphism $\psi_B \colon \overline{B} \to D_{n_2}$ (where $\overline{B} = B \cup \{b\}$ is the closure of $B$ in $D_{n_2}$) in such a way that these homeomorphisms all map $\psi_B(b)$ to the same end point of $D_{n_2}$, then we define $g|_B$ as $\psi_{B'}^{-1} \circ \psi_B$ (where, as before, $B'$ is the unique component adjacent to $g \cdot b$ with color $\sigma_{\iota(g)}(b)\left(k_2(B)\right)$) and we use \cref{lem.patchwork,rmk.patchwork.coloring} to find conclude as before.

In both cases, the fact that the continuity of the group embeddings follows from the fact that the preimage of the stabilizer of a finite set of vertices (branch points) is the stabilizer of a finite set of vertices (branch points).
\end{proof}

\cref{cor.rabbits.simple} allows us to distinguish between the basilica and the other rabbit homeomorphism groups.
The upcoming \cref{prop.rabbit.embeddings} shows that we can distinguish almost every two other homeomorphism groups of rabbits, as it implies that the homeomorphism groups of two distinct rabbits $\R_n$ and $\R_m$ are not isomorphic as soon as $n, m \geq 4$.

\begin{proposition}
\label{prop.rabbit.embeddings}
For every $2 \leq n < m$, there are topological group embeddings
\[ \Homeo(\R_n) \hookrightarrow \Homeo(\R_m). \]
Conversely, $\Homeo(\R_m)$ does not embed into $\Homeo(\R_n)$ for all $4 \leq n < m$.
\end{proposition}

\begin{proof}
As soon as $n \leq m$, the group $\Sym([n])$ embeds into $\Sym([m])$ as the pointwise stabilizer of $[m]\setminus[n]$ points.
With this embedding, the conditions of \cref{lem.embeddings} are satisfied, so the first part of the statement is proved.

For the second part, let $m > n \geq 4$.
Since $\Homeo(\R_m)$ includes copies of the alternating group $\Alt(m)$, it suffices to show that $\Alt(m)$ does not embed into $\U(\Sym([n]),\Aut(S))$.
We do this by contradiction, so let us assume that $\Alt(m)$ acts faithfully on $T_{n,\infty}$ with local actions in $\Aut(S)$ at the vertices of infinite degree.

Note that $m \geq 5$, so $\Alt(m)$ is simple.
Then it is not hard to see that every faithful action of $\Alt(m)$ on a tree must have a fixed vertex $v$ on whose adjacent edges the action is the natural action of $\Alt(m)$ (see \cite[Claim 6.3]{dendriterearrangement} for details).
Since $n < m$, the vertex $v$ of $T_{n,\infty}$ must be of infinite degree, so the local action around $v$ must induce an embedding of $\Alt(m)$ into $\Aut(S)$.
This is impossible:
if such an embedding existed, $\Alt(m)$ would either embed into $\Aut(O)$ or have an index-$2$ subgroup that embeds into $\Aut(O)$, both of which cannot happen because finite subgroups of $\Aut(O)$ must be cyclic.
Thus, $\Alt(m)$ does not embed into $\U(\Sym([n]),\Aut(S)) \simeq \Homeo(\R_n)$, as needed.
\end{proof}

\cref{cor.rabbits.simple} and \cref{prop.rabbit.embeddings} imply that, assuming that $\{n,m\} \neq \{3,4\}$, $\Homeo(\R_n) \simeq \Homeo(\R_m)$ if and only if $n=m$.
For the case $n=3, m=4$, the strategy of \cref{prop.rabbit.embeddings} does not apply because $\Alt(4)$ is not simple.
Indeed, \cite{MR586434} describes a faithful action of $\Sym([4])$ on a $3$-regular tree, so we cannot rule out that $\Aut(4)$ embeds into $\Homeo(\R_3)$.
Hence, the following questions remain open.

\begin{question}
Are $\Homeo(\R_4)$ and $\Homeo(\R_3)$ isomorphic topological groups?
If they are not, does $\Homeo(\R_4)$ embed into $\Homeo(\R_3)$?
\end{question}

Note that \cref{prop.rabbit.embeddings} does not state that there do not exist embeddings of $\Homeo(\R_3)$ into $\Homeo(\B)$.
And indeed there are such embeddings:

\begin{corollary}
\label{cor.3.rabbit.in.basilica.in.airplane}
There are topological group embeddings
\[ \Homeo(\R_3) \hookrightarrow \Homeo(\B) \hookrightarrow \Homeo(\A). \]
\end{corollary}

\begin{proof}
For the first embedding, recall that $\Sym([3])$ is the dihedral group with $6$ elements and thus acts continuously on the circle.
This defines an embedding of $\Sym([3])$ into $\Aut(S)$.
(Note that this fails for $\Sym([m])$ as soon as $m \geq 4$.)
This embedding satisfies the conditions of \cref{lem.embeddings}, so we have that $\Homeo(\R_3) \simeq \U(\Sym([3]),\Aut(S)) \leq \U(\Aut(S),\Aut(S))$, which in turn embeds into $\U(\Aut(S)) \simeq \Homeo(\B)$, so we are done.

Since $\Homeo(\A) \simeq \K(\Aut(S))$ and $\Homeo(\B) \simeq \U(\Aut(S))$, the second embedding is given by direct application of \cref{cor.universal.embeds.into.kaleidoscopic}.
\end{proof}

Embeddings between groups acting on the basilica and the airplane is not new behavior:
as shown in \cite[Section 10]{airplane}, the basilica rearrangement group $T_B$ embeds into the airplane rearrangement group $T_A$ and the overall idea is ultimately not dissimilar to that of \cref{cor.3.rabbit.in.basilica.in.airplane}.

\subsection{Embeddings and non-embeddings among the lamination automorphism groups}

Let us consider the groups $\Homeo^+(\R_n)$ and $\Homeo^+(\A)$ from \cref{sec.lamination.automorphisms}.
As a consequence of \cref{thm.simple.commutators} we immediately have distinction.

\begin{theorem}
\label{thm.pairwise.distinction.+}
The groups $\Homeo^+(\R_n)$ and $\Homeo^+(\R_m)$ are abstractly isomorphic if and only if $n=m$.
\end{theorem}

It is also easy to see that $\Homeo^+(\R_n)$ embeds into $\Homeo^+(\R_m)$ as soon as $n|m$. 
Indeed, if $n|m$ then $[n]$ embeds in $[m]$ via $\iota \colon x\mapsto mx/n$ and $\Cyc(n)$ embeds as the setwise stabilizer of the image of $\iota$ in $\Cyc([m])$, which satisfy the conditions of \cref{lem.embeddings} and so, thanks to \cref{thm.auto.laminations}, one has that $\Homeo^+(\R_n) \simeq \U(\Cyc(n),\Aut(O))$ embeds into $\U(\Cyc(m),\Aut(O)) \simeq \Homeo^+(\R_m)$.

Moreover (and perhaps more surprisingly), these groups embed into $\Homeo^+(\B)$.
Recalling again that $\Homeo^+(\R_n) \simeq \U(\Cyc(n),\Aut(O))$ by  \cref{thm.auto.laminations}, this can be seen as follows:
\[ 
\U(\Cyc(n),\Aut(O)) \hookrightarrow \U(\Aut(O),\Aut(O))
\hookrightarrow \U(\Aut(O)) = \U(\Cyc(2),\Aut(O)) \]
where the first embedding is induced by any faithful action of $\Cyc(n)$ that preserves the cyclic order on the countable dense cyclic order $O$ (this defines an embedding that satisfies the conditions of \cref{lem.embeddings}), the second by observing that the biregular tree on which $\U(\Aut(O),\Aut(O))$ acts is actually the regular tree with infinite degree and the last identity by taking the barycentric subdisivion of such a tree.

As for the embedding between $\Homeo^+(\R_n)$ and $\Homeo^+(\A)$, we simply need to gather previously obtained results.
As we have just seen in the previous paragraph, $\Homeo(\R_n)^+$ embeds into $\U(\Aut(O))$.
By \cref{cor.universal.embeds.into.kaleidoscopic} we have that $\U(\Aut(O)) \leq \K(\Aut(O))$ and by \cref{thm.auto.laminations} we have that $\K(\Aut(O)) \simeq \Homeo(\A)^+$.
Hence, ultimately each $\Homeo^+(\R_n)$ embeds into $\Homeo^+(\A)$.

Here is a collection of what we have discussed in the previous paragraphs.

\begin{corollary}
As soon as $n|m$, there are topological group embeddings
\[ \Homeo^+(\R_n) \hookrightarrow \Homeo^+(\R_m) \hookrightarrow \Homeo^+(\B) \hookrightarrow \Homeo^+(\A). \]
\end{corollary}


\section{Some geometric and analytic properties of the groups}
\label{sec.geometric.analytic.properties}

For topological groups that are not necessarily locally compact, one can generalize compact subsets in different directions.
As introduced by Rosendal in \cite{MR4327092}, a subset $A$ of a topological group $G$ is \textbf{coarsely bounded} if it has finite diameter with respect to every continuous left invariant écart (also known as pseudo-metric) on $G$.
A group is \textbf{locally bounded} if it has a coarsely bounded identity neighborhood.
For $\sigma$-compact locally compact groups, coarsely bounded subsets coincides with relatively compact subsets.
In another direction, one can define \textbf{Roelcke precompact subsets}, which are precompact subsets for the Roelcke uniformity (the meet of the left and right uniformities).
Such subsets $A$ in a topological group are characterized by the following property:
for any identity neighborhood $V$, there is a finite subset $F$ such that $A\subseteq VFV$.
Roelcke precompact subsets are always coarsely bounded.
For further details, see \cite[\S 3]{MR3707816}.

\subsection{Quasi-isometry types of the groups}

Rosendal extended geometric group theory to non-locally compact groups in \cite{MR4327092}.
In particular, he found under what assumptions a Polish group has a well-defined quasi-isometry type, which we determine for the groups under study in this paper.
Let us start with homeomorphism and lamination automorphism groups of the airplane.

\begin{corollary}
\label{cor.airplane.QItype}
The groups $\Homeo(\A)$ and $\Homeo^+(\A)$ are are Roelcke precompact and thus coarsely bounded.
\end{corollary}

\begin{proof}
Thanks to \cref{cor.airplane.group.closed.oligomorphic,,cor.airplane.lamination.group.closed.oligomorphic}, both $\Homeo(\A)$ and $\Homeo^+(\A)$ are oligomorphic.
In particular they are Roelcke precompact.
See \cite[\S 3.1]{MR4327092} for details.
\end{proof}

For regular rabbits, we rely on \S 6 of \cite{MR4327092} about automorphism groups of countable structures.
Our countable structure is the regular tree of infinite degree $T_\infty$ or the bi-regular tree $T_{n,\infty}$ with degrees $(n,\infty)$, where $n\geq3$.

\begin{proposition}
\label{prop.rabbit.QItype}
The groups $\Homeo(\R_n)$ and $\Homeo^+(\R_n)$ are locally Roelcke precompact, in particular locally bounded, and quasi-isometric to $T_\infty$.
\end{proposition}

\begin{proof}
Let $T$ be the tree of circles associated to a regular rabbit.
We prove that the stabilizer of a vertex is Roelcke precompact (even if it is not oligomorphic because it has an orbit for each fixed distance from the fixed vertex) and thus yields the desired result about local Roelcke precompactness.
Let $v$ be a vertex and let $G$ be its stabilizer (in either $\Homeo(\R_n)$ or $\Homeo^+(\R_n)$).
Let $V_n$ be the subtree given by the ball of radius $n$ around $v$ in $T$.
Let $G_n$ be the group of automorphisms of $V_n$ fixing $v$ and with local action in $\Aut(S)$ (or $\Aut(O)$) for all vertices with infinite degree of $V_n$.
We consider the pointwise convergence topology on $G_n$ for its action on vertices of $V_n$.
Let $\Pi_n\colon G\to G_n$ be the restriction map to $V_n$ and similarly, for $m\geq n\in\mathbf{N}$, let $\Pi_{m,n}\colon G_m\to G_n$ be the restriction map to $V_n$.
These maps are continuous group homomorphisms.
Since any element of $G_n$ can be extended to an element of $G$, we see that $G$ is topologically isomorphic to the inverse limit $\varprojlim G_n$ and Roelcke precompactness follows from \cite[Proposition 2.2 (iii)]{MR2929072} as soon as we know that all $G_n$ are Roelcke precompact.
Let us prove that each $G_n$ is actually oligomorphic.

Let $F$ be some finite set of vertices of $V_n$.
Up to enlarging $F$, we may assume that it is center-closed (\cref{def.center.closed}) and contains $v$.
Let us say that two such finite center-closed sets $F$ and $F'$ have the same configuration if there is a partial tree homomorphism (\cref{def.partial.morphism}) between $F$ and $F'$ fixing $v$ that induces isomorphisms for the separation relation (for the cyclic orders in the case of $\Homeo^+(\R_n))$ on neighbors of each vertex of infinite degree of $F$ and that induces isomorphism of cyclic orders on neighbors of finite degree of $F$ in the case of $\Homeo^+(\R_n)$.
For each finite cardinal, there are only finitely many configurations and $G_n$ acts transitively on finite sets with the same configuration, since any isomorphism testifying that two sets $F$ and $F'$ have the same configuration can be extended to an element of $G_n$ by \cref{lem.patchwork,rmk.patchwork.coloring}.
So $G_n$ is an oligomorphic group.

Now let us prove that $\Homeo(\R_n)$ and $\Homeo^+(\R_n)$ are quasi-isometric to $T_\infty$.
Since these groups act transitively on the set of (unoriented) edges of $T$, we can rely on the argument given in \cite[Example 6.34]{MR4327092}.
The fact that $T$ is quasi-isometric to $T_\infty$ follows from the following lemma.
\end{proof}

\begin{lemma}
For all $k\geq2$, the biregular tree $T_{k,\infty}$ is quasi-isometric to $T_\infty$.
\end{lemma}

\begin{proof}
Let us fix infinite-degree vertices $u_0$ in $T_{k,\infty}$ and $v_0$ in $T_\infty$ and $U_n, V_n$ the respective $n$-neighborhood of these points.
We define by induction  surjective $(1,2)$-quasi-isometric maps $f_n\colon V_n\to U_n$ such that $f_{n+1}$ restricted to $V_n$ is $f_n$. The desired quasi-isometry $f$ will be defined to coincide with $f_n$ on $V_n$.

For $n=0$, we set $f_0(v_0)=u_0$.
We extend $f_0$ to $f_1$ by choosing any bijection between the infinite countable sets $V_{1}\setminus\{v_0\}$ and $U_1\setminus\{u_0\}$.
Assume that $f_n$ has been defined for some $n\geq1$.
We define $f_{n+1}$ to be $f_n$ on $V_n$. If $n$ is odd, for each vertex $v$ in $V_{n}\setminus V_{n-1}$, we choose a surjective map between neighbors of $v$ not in $V_n$ and the $k-1$ neighbors of $f(v_n)$ not in $U_n$.
This defines $f_{n+1}$ on the 1-neighborhood of $v$ and we do the same for all other vertices in $V_n\setminus V_{n-1}$. 

Assume now that $n$ is even.
Let us fix $v\in V_n\setminus V_{n-1}$.
Let $N_v$ be the collection of all neighbors of vertices in $f_n^{-1}(f_n(v))$ except for the one in $V_{n-1}$ and $N_{f(v)}$ the collection of all neighbors of $f(v)$ except the one in $U_{n-1}$.
Since $N_v$ and $N_{f(v)}$ are both infinite countable, we can choose $f_{n+1}$ so that its restriction on $N_v$ coincides with some bijection $N_v \to N_{f(v)}$.
We do the same for all distinct sets of the form $N_v$ to define $f_{n+1}$. 

We have thus defined a surjective map $f\colon T_\infty\to T_{k,\infty}$ mapping adjacent vertices to adjacent vertices.
In particular, it is 1-Lipschitz and, for any vertex $v$, $d(f(v),f(v_0))=d(v,v_0)$.
Moreover, by construction $f$ is injective on the subset of vertices in $T_\infty$ at odd distance from $v_0$.
Let $x,y$ be vertices in $T_\infty$ and $v_1=x,v_2,\dots,v_{n+1}=y$ the shortest path between them.
Let $v_l$ be the unique point at shortest distance from $v_0$ on this path.
If $v_l$ is at even distance from $v_0$ then $f(v_{l-1})$ and $f(v_{l+1})$ are distinct and $f$ is then injective on this path.
Otherwise, if $v_l$ is at odd distance from $v_0$ and $f(v_{l-1})=f(v_{l+1})$, then $f(v_{l-2})$ and $f(v_{l+2})$ are distinct, as are the other points on the path.
So $d(x,y)-2\leq d(f(x),f(y))\leq d(x,y)$, which means that $f$ is a quasi-isometry.
\end{proof}

\subsection{Property (T) and the Haagerup property}\label{sec.haagerup.property.(T)}
 
Recall that a topological group $G$ has \textbf{strong Property (T)} if $G$ has a Kazhdan pair $(F,\varepsilon)$ with $F$ finite (see \cite{MR2415834} for details about Property (T)).
In the opposite direction, the Haagerup property has been extended from locally compact groups to general topological group $G$ in \cite{MR3707816}.
A topological group $G$ has the \textbf{Haagerup property} if $G$ has a continuous action by isometries on a Hilbert space $\mathcal{H}$ such that there is an orbit map that is a coarse embedding from $G$ to $\mathcal{H}$ (i.e., the preimages of bounded subsets are coarsely bounded).
These properties are opposite in the sense that a topological group with both property (T) and the Haagerup property is coarsely bounded.

\begin{remark}
Thanks to the GNS construction, the Haagerup property is equivalent to the existence of a continuous coarsely proper function conditionally of negative type.
See \cite[\S 2.10]{MR2415834} for explanations.
\end{remark}

\begin{theorem}
\label{thm.property.T}
The groups $\Homeo(\R_n)$ and $\Homeo^+(\R_n)$ have the Haagerup property for all $n \geq 2$, whereas $\Homeo(\A)$ and $\Homeo^+(\A)$ both have strong Property (T).
\end{theorem} 

\begin{proof} 
For the airplane, thanks to the identification of $\Homeo(\A)$ and $\Homeo^+(\A)$ with oligomorphic kaleidoscopic groups (\cref{cor.airplane.group.closed.oligomorphic,,cor.airplane.lamination.group.closed.oligomorphic}), strong property (T) follows immediately from \cite[Theorem 1.1]{MR3417738}.

Now consider the rabbits.
Let $d$ be the distance function on $T_\infty$ or $T_{n,\infty}$ and $v$ some fixed vertex in this tree, the map $f\colon g\mapsto d(g(v),v)$ is a continuous, metrically proper and conditionally of negative type function \cite[C.2.2.(iii)]{MR2415834} and thus gives rise to a continuous proper affine isometric on some Hilbert space (that is the Hilbert space of square-integrable functions on the countable set of oriented edges).
\end{proof}

By \cref{thm.property.T}, the group $\Homeo^+(\A)$ is a new example of an uncountable group with property (T) acting non-elementarily on the circle, such as that studied in \cite{KaleidoscopicOnCircle}.


\section{Universal minimal flows in the case of the airplane}
\label{sec.universal.minimal.flows}

Recall that a topological group $G^*$ is \textbf{extremely amenable} if all of its continuous actions on compact spaces have a fixed point, or equivalently its universal minimal flow $\mathcal{M}(G^*)$ is a singleton.
Moreover, a subgroup $G^*$ of a group topological $G$ is \textbf{coprecompact} if the completion of its quotient $G/G^*$ is compact, or equivalently if there is a finite subset $F \subseteq G$ such that $G = VFH$.
In particular, a permutation group $G^* \leq \Sym(\Omega)$ is coprecompact if and only if it is oligomorphic.

We denote by $\widehat{X}$ the completion of a uniform space $X$.
In particular, if $H$ is a closed subgroup of a topological group $G$, then $G/H$ is endowed with the uniformity coming from the right uniformity on $G$ and we denote its completion by $\widehat{G/H}$.

Assume $G^*$ is a closed extremely amenable coprecompact subgroup of $G$.
Then any minimal subflow of the completion $\widehat{G / G^*}$ of the quotient $G/G^*$ is $\mathcal{M}(G)$.
In particular, if $\widehat{G / G^*}$ itself is minimal then it is $\mathcal{M}(G)$.

Universal minimal flows of kaleidoscopic groups have been studied in \cite{Kaleidoscopic2}.
Using and expanding from ideas from that paper, we can obtain a more precise description of the universal minimal flow of $\Homeo(\A)$ and $\Homeo^+(\A)$.

Consider the dendrite $D_\infty$, identified with the dendrite of circles of $\A$ (\cref{sub.dendrite.of.circles}).
Fix an end point $\xi$ of $D_\infty$.

\begin{proposition}
The stabilizer subgroup $\K(\Aut(O))_\xi$ is an extremely amenable oligomorphic group.
\end{proposition}

\begin{proof}
The stabilizer subgroup in $\Aut(O)$ of a point is isomorphic to $\Aut(\mathbf{Q},<)$ by \cref{lem.aut.circular.order.and.relation.properties} and thus is extremely amenable (see \cite{MR1608494}).
Since $\Aut(O)$ is transitive, extreme amenability of $\K(\Aut(O))_\xi$ follows from \cite[Corollary 8.3(i)]{Kaleidoscopic2}.
The fact that this group is oligomorphic follows from Proposition 4.6 and Corollary 5.10 in \cite{Kaleidoscopic2}.
\end{proof}

\begin{theorem}
\label{thm.universal.minimal.flow.airplane.completion}
The universal minimal flow of $\K(\Aut(O))$ is $$\widehat{\K(\Aut(O))/\K(\Aut(O) )_\xi}.$$
\end{theorem}

\begin{proof}
Since $\K(\Aut(O) )_\xi$ is oligomorphic for the action on the branch points of $D_\infty$, it is coprecompact in $\K(\Aut(O))$, so the completion of the quotient space coincides with the Samuel compactification.
Thus, the desired result follows from the fact that $\K(\Aut(O) )_\xi$ is extremely amenable and from applying Theorem 1.4 in \cite{Kaleidoscopic2} with $\Gamma=\Aut(O)$, $\Delta=\Gamma_c\simeq \Aut(\mathbf{Q},<)$, $\K^*(\Delta)=\K(\Aut(O) )_\xi$.
\end{proof}

We can find a more concrete description of the universal minimal flow of the group $\Homeo^+(\A)\simeq\K(\Aut(O))$ and thus, as will be shown soon in \cref{cor.c.ordered}, give an example of a Polish group that is intrinsically c-ordered (as asked in \cite[Question 4.9]{Glasner-Megrelishvili-2018}).
The universal minimal flow has a continuous invariant cyclic order, but admits a minimal action on a compact space without invariant cyclic order, as we will see in \cref{prop.no.invariant.cyclic.order}

Let us recall a few definitions from \cite{Glasner-Megrelishvili-2018,Glasner-Megrelishvili-2021}.
Let $\Omega$ be some cyclically ordered set and $A\subset\Omega$.
Then there exists a cyclically ordered set $\Split(X,A)$ with an order-preserving map $\nu\colon\Split(X,A)\to X $ such that any point $x\in X\setminus A$ has a unique preimage and any point $a\in A$ has two preimages $a^-, a^+$ (see \cite[Lemma 2.11]{Glasner-Megrelishvili-2018}).
The cyclic order on $\Split(X,A)$ is such that $[a^-,a^+,x]$ for any $x\neq a^-,a^+$ and $[x,y,z]\iff [\nu(x),\nu(y),\nu(z)]$ if $x,y,z\in \Split(X,A)$ have distinct images in $X$.
If $X$ is compact (for the topology induced by the cyclic order, see \cite[\S2]{Glasner-Megrelishvili-2018}), then so is $\Split(X,A)$.
If some topological group $G$ acts continuously on $X$ preserving the cyclic order and $A\subset X$ is $G$-invariant, then $G$ acts continuously and order-preservingly on $\Split(X,A)$ as well.

For $\Homeo^+(\A)$ and its action on $\mathbf{S}^1$ constructed in \cref{sec.lamination.automorphisms}, we consider the set $P$ of preimages of cut points of $\A$ that belong to circles.
This is a countable set and the Carathéodory loop $\varphi$ is $2$-to-$1$ on this subset.

\begin{theorem}
\label{thm.universal.minimal.flow.airplane}
The Universal minimal flow of $\Homeo^+(\A)$ is $\Split(\mathbf{S}^1,P)$.
\end{theorem}

\begin{proof}
As before, let us fix an end point $\xi$ of the airplane.
Throughout this proof, we denote $\Homeo^+(\A)$ by $G$ and $\Homeo^+(\A)_\xi$ by $H$.
Let $E$ be the set of points $y\in\Split(\mathbf{S}^1,P)$ such that $\varphi(\nu(y))$ is an end point of the airplane.
Since $\varphi$ is bijective when restricted to the set $E$, there exists a unique $x\in E$ such that $\varphi(\nu(x))=\xi$.
Thanks to \cref{thm.universal.minimal.flow.airplane.completion} and the fact that $G$ is $\K(\Aut(O))$, it suffices to prove that the orbit map 
\[ \begin{matrix}
G/H&\to&\Split(\mathbf{S}^1,P)\\
gH&\mapsto&gx
\end{matrix} \]
is a bi-uniformly continuous injection with dense image.

It is injective with image $E$ (which is dense in $\Split(\mathbf{S}^1,P)$) because the action of $G$ on the set of end points is transitive (this is the case for every kaleidoscopic group on $D_\infty$, including $\K(1)$, see \cite[Proposition 5.5]{Kaleidoscopic}).

For a cycle $(b_1,\dots,b_m)$ in $P$ (i.e., the map $i\mapsto b_i$ is a cyclic order-preserving map from $[n]$ to $P$), we define $V_{b_1,\dots,b_m}=\{(y_1,y_2)\in E^2 \mid \exists i,\ y_1,y_2\in(b_i^+,b_{i+1}^-)\}$.
The collection of all $V_{b_1,\dots,b_m}$ is a basis for the uniform structure on $E$.
The map $gH\mapsto gx$ is uniformly continuous because the action of $G$ on $\Split(\mathbf{S}^1,P)$ is continuous.
For any $y\in E$, we choose $g_y\in G$ such that $g_yx=y$.
We want to show that $y\to g_y H$ is uniformly continuous.

By definition of the Polish topology on $\K(\Aut(O))$, a base for the uniform structure on $G/H$ is 
\[ V_F=\{(gH,hH),\ g\in U_FhH\} \]
where $U_F$ is the pointwise stabilizer in $G$ of a finite set $F$ of branch points of $D_\infty$.
Let $F_0$ be some finite center-closed set of branch points in $D_\infty$ containing at least two points.
Let $F_1$ be the (finite) subset of $\mathcal{A}$ consisting of all projections of a circle $C_1$ to a circle $C_2$ where $C_1,$ $C_2$ are distinct preimages of elements of $F$ under the quotient map $\A\to D_\infty$.
Finally, let $F\subset \Split(\mathbf{S}^1,P)$ be the subset of all preimages of $F_1$ under $\varphi\circ\nu$.
Let us order cyclically the points of $F$: $(b_1,\dots,b_m)$.
For any pair $(y_1,y_2)\in V_{b_1,\dots,b_m}$, there is $u\in U_F$ such that $gy_1=y_2$ (by the same construction as in the proof of \cite[Theorem 1.14]{Duchesne-2023}).
Thus, if $(y_1,y_2)\in V_{b_1,\dots,b_m}$ then $(g_{y_1}H,g_{y_2}H)\in V_F$.
This proves that the inverse of the orbit map is uniformly continuous.
Since $E$ is dense in the compact space  $\Split(\mathbf{S}^1,P)$, we have that $G$ and $\Split(\mathbf{S}^1,P)$ are homeomorphic as $G$-spaces.
\end{proof}

For a topological group $G$, a flow $X$ is strongly proximal if any probability measure contains a Dirac mass in the closure of its orbit (for the induced action of $G$ on $\Prob(X)$).
Among strongly proximal minimal $G$-flows, there is a universal one (called the universal minimal strongly proximal $G$-flow or the Furstenberg boundary of $G$).
The following corollary shows that the Furstenberg boundary of $\Homeo^+(\A)$ is $\Split(\mathbf{S}^1,P)$.

\begin{corollary}
Any minimal $\Homeo^+(\A)$-flow is strongly proximal.
\end{corollary}

\begin{proof}
It suffices to prove that the action of on $\Split(\mathbf{S}^1,P)$ is strongly proximal. The argument is the same as in \cite[Corollary 5.1.]{KaleidoscopicOnCircle}
\end{proof}

Let us recall that a topological group $G$ is \textbf{intrinsically c-ordered} if its universal minimal flow has an invariant cyclic order that induces the topology.
This property implies in particular that any minimal flow of $G$ is tame, i.e. the envelopping semigroup is separable and Fréchet, so it does not contain a topological copy of the Stone-Cech compacitifion of the integers $\beta\mathbf{N}$ (see \cite{Glasner-Megrelishvili-2018} for this fact and \cite{Glasner,MoreTameDynSys} for the relevance of tameness).
An immediate consequence of \cref{thm.universal.minimal.flow.airplane} is the following.

\begin{corollary}
\label{cor.c.ordered}
The Polish group $\Homeo^+(\A)$ is intrinsically c-ordered.
\end{corollary}

We have the following sequence of $\Homeo^+(\A)$-factors: \[ \Split(\mathbf{S}^1,A)\to\mathbf{S}^1\to\A\to D_\infty. \]
The first two ones are cyclically ordered.

\begin{proposition}
\label{prop.no.invariant.cyclic.order}
Both the dendrite $D_\infty$ and the airplane $\A$ have no cyclic order that is invariant under $\Homeo^+(\A)$.
\end{proposition}

\begin{proof}
Let $X$ be either $D_\infty$ or $\A$.
If $X=D_\infty$, let $p_1$, $p_2$ and $p_3$ be branch points such that $p_2$ is between the other two.
If $X=\A$, let $p_1$, $p_2$ and $p_3$ be regular cut points (\cref{def.points.airplane}) such that $p_2$ is between the other two.
In both cases, an application of \cref{lem.patchwork,rmk.patchwork.coloring} shows that there exists a $g \in \Homeo^+(\A)$ that fixes $p_2$ and switches $p_1$ and $p_3$, so both $D_\infty$ and $\A$ admit no cyclic order that is invariant under $\Homeo^+(\A)$.
\end{proof}

\begin{proposition}
\label{prop.UMF.homeo.airplane.is.metrizable}
The universal minimal flow of $\Homeo(\A)$ is metrizable.
\end{proposition}

\begin{proof}
Let us use the identification $\Homeo(\A)\simeq\K(\Aut(S))$ from \cref{thm.homeo.airplane.is.kaleidoscopic}.
Since $\Aut(S)$ acts transitively and oligomorphically on $(\Omega,S)$ by \cref{lem.aut.circular.order.and.relation.properties} and since $\mathcal{M}(\Aut(S))$ is metrizable by \cref{prop.univ.min.flow.aut(s)}, the result follows immediately from \cite[Theorem 1.1]{Kaleidoscopic2}.
\end{proof}

\begin{remark}
It is natural to ask what happens for the basilica and regular rabbits with at least two ears.
For the airplane, we relied on the work \cite{Kaleidoscopic2} on the universal minimal flow of kaleidoscopic groups;
however, to the best of our knowledge, no such work exists for universal groups.
Let us explain what could be a strategy for homeomorphism groups of rabbits.
Let $\xi\in\R_n$ be some end point, it corresponds to an end (which we also denote by $\xi$) in the tree $T_{n,\infty}$ of circles of $\R_n$.
Its stabilizer $\Homeo^+(\R_n)_\xi$ has a continuous surjective homomorphism to $\mathbf{Z}$ coming from the Busemann cocycle associated to $\xi$.
In particular, $\Homeo^+(\R_n)_\xi$ is not extremely amenable and the universal minimal flow of $\Homeo^+(\R_n)$ is not the completion of $\Homeo^+(\R_n)/\Homeo^+(\R_n)_\xi$ for the left uniform structure.
Instead, the universal minimal flows of $\Homeo^+(\R_n)$ and $\Homeo(\R_n)$ can likely be obtained using by some suspension process as in \cite[\S7]{Duchesne-2023}.
\end{remark}


\begin{appendices}

\section{The countable dense cyclic order and its separation relation}
\label{sec.order.and.separation}

In this section we gather some useful facts about the dense cyclic order $O$ and its separation relation $S$ that we use throughout the paper.

\subsection{Cyclic orders and their separation relations}

\begin{definition}[Cyclic order]
\label{def.cyclic.order}
Let $\Omega$ be a set.
A \textbf{cyclic order} on $\Omega$ is a ternary relation $[\ ,\ ,\ ]$ on $X$ that enjoys the following properties:
\begin{enumerate}
    \item \textbf{Cyclicity:} If $[a, b, c]$ then $[b, c, a]$,
    \item\textbf{Asymmetry:} If $[a, b, c]$ then not $[c, b, a]$,
    \item\textbf{Transitivity:} If $[a, b, c]$ and $[a, c, d]$ then $[a, b, d]$
    \item\textbf{Connectedness:} If $a, b,c$ are distinct, then either $[a, b, c]$ or $[c, b, a]$.
\end{enumerate}
\end{definition}

Given a cyclic order, one can define the inverse cyclic order $[x,y,z]^{-1}\iff[z,y,x]$.
It is readily seen that this is a well-defined cyclic order.

Let us denote by $\Omega^{(3)}$ the set of all triples $(x,y,z) \in \Omega^3$ such that $x,y,z$ are distinct.
Every cyclic order defines an \textbf{orientation map}
\[ O \colon \Omega^{(3)} \to \{+,-\}, \]
that is $O(x,y,z)=+ \iff [x,y,z]$.
Since a cyclic order is entirely determined by its orientation map, we take the liberty of using the symbol $O$ to also refer to the cyclic order.

Following Vailati \cite[\S204, p. 215]{Russel}, we introduce the following quaternary relation (see Coxeter \cite[Chapter 3]{Coxeter}, for a more recent reference).

\begin{definition}[Separation relation]
\label{def.separation.relation}
Let $\Omega$ be a set equipped with a cyclic order $[\ ,\ ,\ ]$.
The \textbf{separation relation} associated to it is the relation $S$ such that $S(a,b,c,d) \iff [a,b,c]\wedge [c,d,a]$ or $[a,d,c]\wedge [c,b,a]$.
\end{definition}

It is easy to verify that a cyclic order and its inverse induce the same separation relation.

\begin{example}
Let $\Omega=\mathbf{S}^1$ be the unit circle in $\mathbf{C}$ with the following cyclic order.
For $z_0,z_1,z_2\in \mathbf{S}^1$, choose two preimages $x_1,x_2\in\mathbf{R}$ of $z_1,z_2$ via the exponential map $x\mapsto \exp(ix)$ that lies in the same connected component of $\mathbf{R}\setminus\exp^{-1}(z_0)$ and write $[z_0,z_1,z_2]\iff x_1<x_2$ for the usual order on $\mathbf{R}$.
\end{example}

\begin{example}
The separation relation associated to the cyclic order on $\mathbf{S}^1$ can be described in a topological way:
$S(a,b,c,d)$ holds if and only if $b$ and $d$ lie in different connected components of $\mathbf{S}^1\setminus\{a,c\}$.
That is, $\{b,d\}$ is separated by $\{a,c\}$ and vice versa.
\end{example}

Linear orders, cyclic orders, betweenness and separations relations are closely related, as explained in \cite{Huntington}.
In particular, given a cyclic order on a set $\Omega$ and $x\in\Omega$, we can set $y<_xz$ if $(x,y,z)$ is cyclically ordered.
The binary relation $<_x$ is a linear order on $\Omega$.
Thus, given a subset $F$ of a cyclically ordered set $\Omega$, we call \textbf{cyclic successor} of $x\in F$ the minimum of the linear order $<_x$ on $F\setminus\{x\}$ (if it exists).
Conversely if $<$ is a linear order on some set $\Omega$, we can define a cyclic order via $[x,y,z]\iff x<y<z\vee y<z<x \vee z<x<y$.
In particular, we call \textbf{standard} the cyclic order on $[n]=\{1,\dots,n\}$ induced by the  natural linear order $<$.

Let $\Omega$ be a set with a cyclic order $[\ ,\ ,\ ]$.
We say that a bijection $g$ of $\Omega$ \textbf{preserves} the cyclic order if $[x,y,z]\iff [g(x),g(y),g(z)]$ and that it \textbf{reverses} the order if $[x,y,z]\iff [g(x),g(y),g(z)]^{-1}$.
If $S$ is a separation relation, we say that a bijection $g$ of $\Omega$ \textbf{preserves} $S$ if $S(x,y,z,w)\iff S(g(x),g(y),g(z),g(w))$.

\begin{lemma}
\label{lem.preserves.respects.cyclic.order}
Let $\Omega$ be a set with a cyclic order and let $S$ be the associated separation relation.
A bijection $g$ of $\Omega$ preserves $S$ if and only if it preserves or reverses the cyclic order.
\end{lemma}

\begin{proof}
Since the cyclic order and its inverse cyclic order induce the same separation relation, it is clear that a bijection that preserves or reserves the cyclic order preserves the separation relation.

Conversely, assume that $(u,v,w)$ satisfies $[u,v,w]$.
Then, for any triple $a,b,c$,
\begin{align*}
    [a,b,c]&\iff &(S(a,u,v,w)\wedge S(a,b,v,w)\wedge S(a,b,c,w))\\
    &&\vee(S(b,u,v,w)\wedge S(b,c,v,w)\wedge S(b,c,a,w))\\
    &&\vee(S(c,u,v,w)\wedge S(c,a,v,w)\wedge S(c,a,b,w))\\
    &&\vee(S(a,u,v,w)\wedge S(a,b,v,w)\wedge S(a,b,w,c))\\
    &&\vee(S(b,u,v,w)\wedge S(b,c,v,w)\wedge S(b,c,w,a))\\
    &&\vee(S(c,u,v,w)\wedge S(c,a,v,w)\wedge S(c,a,w,b))
\end{align*}
Let us write $R(a,b,c,u,v,w)$ for the relation with $6$ variables that appears on the right-hand side of the equivalence above.
Let us fix $u,v,w \in \Omega$ such that $[u,v,w]$ and let us assume that $g$ preserves $S$.
If $[g(u),g(v),g(w)]$ then $[a,b,c]$ implies $R(a,b,c,u,v,w)$ and, since $g$ preserves $S$, $R(g(a),g(b),g(c),g(u),g(v),g(w))$ and thus $[g(a),g(b),g(c)]$.
Similarly if $[g(u),g(v),g(w)]^{-1}$, for $a,b,c$ such that $[a,b,c]$, $R(a,b,c,u,v,w)$ holds and thus $[g(a),g(b),g(c)]^{-1}$.
This means that $g$ preserves or reverses the cyclic order.
\end{proof}

\subsection{Countable dense cyclic orders}
 
A cyclic order is \textbf{dense} if for all distinct $x,z\in \Omega$, there is $y \in \Omega$ such that $[x,y,z]$.
There is a unique countable dense cyclic order $O$ up to isomorphism (this follows from the analogous statement for linear orders).
Such a cyclic order $O$ has been studied in \cite{Truss}.
The Dedekind completion of $O$ is $\mathbf{S}^1$ with its usual topology.

\begin{example}
A realization of this unique dense countable cyclic order can be obtained by choosing any dense countable subset of the circle $\mathbf{S}^1$.
For example, if we identify $\mathbf{S}^1$ with $\mathbf{R}/\mathbf{Z}$ then $\mathbf{Q}/\mathbf{Z}$ and $\mathbf{Z}[1/2]/\mathbf{Z}$ are standard choices.
\end{example}

Throughout the rest of this appendix and all other sections, we denote by $O$ this countable dense cyclic order and by $S$ the associated separation relation.
We denote by $\Aut(S)$ the group of all bijections of $\Omega$ that preserve this separation relation $S$ and by $\Aut(O)$ the subgroup of $\Aut(S)$ of bijections of $\Omega$ that preserve this cyclic order $O$.
Let us gather a few facts about these two groups.

\begin{remark}
\label{rmk.Aut.S.C.closed}
The structures $(\Omega,O)$ and $(\Omega,S)$ are Fraïssé limits of all finite cyclic orders and all finite separation relations respectively.
In particular, it follows that $\Aut(S)$ and $\Aut(O)$ are closed subgroups of $\Sym([\infty])$ for the pointwise convergence topology.
\end{remark}

\begin{lemma}
\label{lem.aut.circular.order.and.relation.properties}
The groups $\Aut(S)$ and $\Aut(O)$ enjoy the following properties.
\begin{enumerate}
    \item $\Aut(S)$ is the semi-direct product $\Aut(O)\rtimes\mathbf{Z}/2\mathbf{Z}$.
    \item $\Aut(O)$ is 2-transitive.
    \item $\Aut(S)$ is 3-transitive.
    \item $\forall x \in X$, the stabilizer $\Aut(O)_x$ is isomorphic to $\Aut(\mathbf{Q},<)$.
    \item $\forall x \in X$, the stabilizer $\Aut(S)_x$ is isomorphic to $\Aut(\mathbf{Q},B)$ where $B$ is the betweenness relation associated to $<$.
    \item $\Aut(S)$ and $\Aut(O)$ are oligormorphic groups.
    \item $\Aut(O)$ is a simple group.
\end{enumerate}
\end{lemma}

\begin{proof}\

1. Any bijection preserving the separation relation either preserves or reverses the cyclic order, so $\Aut(O)$ has index 2 in $\Aut(S)$.
If $\sigma\in\Aut(S)\setminus\Aut(O)$ is an involution (for example the map induced by $x\mapsto 1-x$ in the model $\Omega=\mathbf{Q}/\mathbf{Z}$), then it is readily verified that $\Aut(S) = \Aut(O)\rtimes<\sigma>$.

2. We use $\mathbf{Q}/\mathbf{Z}$ as a model for $\Omega$. Let $(a,b)$ and $(a',b')$ be two pairs of distinct points in $\mathbf{Q}/\mathbf{Z}$.
Since the left action of $(\mathbf{Q},+)$ on $\mathbf{Q}/\mathbf{Z}$ is transitive and preserves the cyclic order, we may assume that $a=a'=0$. Now the following piecewise affine map, where $\mathbf{Q}/\mathbf{Z}$ is identified with $\mathbf{Q}\cap[0,1)$, does the job:
\[ \begin{array}{rcl}
\mathbf{Q}/\mathbf{Z}&\to&\mathbf{Q}/\mathbf{Z}\\
x\in[0,b)&\mapsto& \frac{b'x}{b}\\
x\in[b,1)&\mapsto& b'+\frac{(1-b')(x-b)}{1-b}\\
\end{array} \]

3. Let $(a,b,c)$ and $(a',b',c')$ be two triples of distinct points in $\mathbf{Q}/\mathbf{Z}$.
By the previous argument, we may assume that $a=a'=0$ and $b=b'=1/2$.
If $c$ and $c'$ belong to the same interval $(0,1/2)$ or $(1/2,1)$, a piecewise affine map as above maps $c$ to $c'$ and preserves the cyclic order.
Otherwise the symmetry at $1/2$ swaps these two intervals preserving the separation relation and we are in the previous situation.

4. Once $x\in\Omega$ is fixed, the linear order $y<_xz\iff[x,y,z]$ on $\Omega\setminus\{z\}$ is invariant and isomorphic to $(\mathbf{Q},<)$.

5. As above, the linear order $<_x$ is preserved or reversed.

6. The group of orientation-preserving piecewise affine maps (over $\mathbf{Q}$) in $\Aut(O)$ acts transitively on cycles $(x_1,\dots,x_n)$, i.e., images of cyclic order-preserving maps $[n]\to\mathbf{Q}/\mathbf{Z}$.
So two $n$-tuples $(x_1,\dots,x_n)$ and $(y_1,\dots,y_n)$ are in the same orbit if and only if $[x_i,x_j,x_k]\iff[y_i,y_j,y_k]$.
Since there are $(n-1)!$ cyclic orders on $[n]$, the group $\Aut(O)$ (and thus $\Aut(S)$) acts oligomorphically on $\Omega$.

7. This is \cite[Proposition 4.2.8]{Macpherson}.
\end{proof}

\begin{remark}
Mathematicians with algebraic taste may think of $\mathbf{Q}/\mathbf{Z}$ as the projective line $\mathrm{P}^1\mathbf{Q}$ over $\mathbf{Q}$.
Under this point of view, the group of piecewise affine maps used above is the group $\mathrm{PL}(\mathrm{P}^1\mathbf{Q})$ of piecewise projective homeomorphisms of $\mathrm{P}^1\mathbf{Q}$ (or $\mathbf{Q}$ if $\infty$ is fixed).
This is a countable group studied in \cite{Monod-2013} as a discrete group.
With the pointwise convergence topology, $\mathrm{PL}^+(\mathrm{P}^1\mathbf{Q})$ is extremely amenable. This is a consequence of \cite[Main Theorem]{MR1608494}.
\end{remark}

From the fact that they are oligomorphic, we know that these groups have strong property (T) and are coarsely bounded.

\begin{corollary} 
\label{cor.aut.circular.order.generated.by.stabilizers}
Both groups $\Aut(S)$ and $\Aut(O)$ are generated by their point stabilizers.
\end{corollary}

\begin{proof}
By points 2 and 3 of \cref{lem.aut.circular.order.and.relation.properties}, $\Aut(S)$ and $\Aut(O)$ are 2-transitive, so they are generated by their point stabilizers.
\end{proof}

We conclude this appendix by emphasizing that the universal minimal flow of $\Aut(O)$ has been shown in \cite{Glasner-Megrelishvili-2021} to be $\Split(\mathbf{S}^1,\Omega)$ where $\Omega$ is any countable dense subset of $\mathbf{S}^1$.
Since $\Aut(S)\simeq\Aut(O)\rtimes\mathbf{Z}/2\mathbf{Z}$, we fall in a situation studied in \cite[Appendix 3]{AlexanderS2012} where an induced action of $\Aut(S)$ on $\Split(\mathbf{S}^1,\Omega)\times \mathbf{Z}/2\mathbf{Z}$ is defined via a lifting $\mathbf{Z}/2\mathbf{Z}\to\Aut(S)$.
In this situation, it can be shown that any such action is isomorphic to the diagonal action of $\Aut(S)$ on $\Split(\mathbf{S}^1,\Omega)\times\Aut(S)/\Aut(O)$, so we have the following fact.

\begin{proposition}
\label{prop.univ.min.flow.aut(s)}
The universal minimal flow of $\Aut(S)$ is $\Split(\mathbf{S}^1,\Omega)\times \mathbf{Z}/2\mathbf{Z}$ with the induced action  from the action of $\Aut(O)$ on $\Split(\mathbf{S}^1,\Omega)$. 
\end{proposition}


\section{Rabbits and airplane as Julia sets}
\label{sec.julia.sets}

In this section we show that the famous basilica, Douady rabbit and airplane Julia sets satisfy our \cref{def.rabbit} and \cref{def.airplane}.

We rely on \cite{douady-hubbard-I,douady-hubbard-II} (an English translation is available on John Hubbard's webpage \cite{Douady-Hubbard-English}) and refer to this for definitions and classical facts about the Mandelbrot, Julia and Fatou sets.
We will merely introduce the notations we use.

We denote by $f_c$ the polynomial map $z\mapsto z^2+c$, where $c\in\mathbf{C}$.
Its Julia set is denoted by $J_{c}$ and $K_c$ denotes the filled-in Julia set.
When $c$ belongs to the Mandelbrot set, $J_c$ and $K_c$ are connected.
Moreover, if $f_c$ is subhyperbolic then $J_c$ is locally connected \cite[Proposition 4, Exposé III]{douady-hubbard-I}, so in this case $J_c$ is a Peano continuum and there is a continuous surjective map $\varphi\colon \mathbf{R}/\mathbf{Z}\to J_c=\partial K_c$, called the Carathéodory loop, that we also employed in \cref{sec.lamination.automorphisms}.
Dynamically, the doubling map $2\mapsto 2x$ is semi-conjugate to the action of $f_c$ on $J_c$:
$$\forall x\in\mathbf{R}/\mathbf{Z},\ \varphi(2x)=f_c(\varphi(x)).$$
For an $x\in J_c$, its preimages (often identified with the representatives in $[0,1)$) are called the external arguments of $x$.

By \cite[Proposition 3, Exposé II]{douady-hubbard-I}, the boundary of any component of $K_{c}$ is a circle.
Conversely, If $C$ is a circle in $J_{c}$ then the bounded component $B$ of $\mathbf{C}\setminus C$ is in $K_{c}$ (the complement of $K_{c}$ is connected) and it is open, so $B\subset U$ for some component $U$ of $\mathring{K_{c}}$ and thus $B=U$.
In conclusion, circles in $K_{c}$ are exactly boundaries of components of $\mathring{K_{c}}$.

We will only consider points $c$ in the Mandelbrot set with external arguments of the form $\frac{p}{2^k-1}$.
These points are roots of hyperbolic components by \cite[Exposé VIII, Théorème 1]{douady-hubbard-II}.

\subsection{The rabbit Julia sets}

\begin{proposition}
Let $\theta_n=\frac{1}{2^n-1}$ for $n\in\mathbf{N}_{\geq2}$ and $c_n$ the point of the Mandelbrot set with external argument $\theta_n$.
The Julia set $J_{c_n}$ is a regular $n$-rabbit.
\end{proposition}

\begin{remark}
\label{rmk.basilica.rabbit}
For $n=2,3$, these Julia sets are known as the basilica and the Douady rabbit, respectively.
They are often defined by choosing $c$ to be the center (where we use the root) of the corresponding hyperbolic component of the Julia set but this gives homeomorphic Julia set thanks to \cite[Exposé XVIII, Proposition 1]{douady-hubbard-II}.
\end{remark}

\begin{proof}
Let us fix $n\geq2$, let $\theta_n=\frac{1}{2^n-1}\in\mathbf{Q}/\mathbf{Z}$ and let $c_{n}$ be the point in the Mandelbrot set with external argument $\theta_n$. The external argument $\theta_n$ is $n$-periodic for the doubling map.
There is a unique other external ray landing at $c_{n}$ with same period \cite[Proposition 2]{Douady}, so it is of the form $\frac{k}{2^n-1}$.
The external argument of this ray is $\theta_n'=\frac{2}{2^n-1}$ (for example because of \cite{Douady-86}).
The point $\alpha_1$ with external argument $\theta_n$ in the Julia set $J_{c_{n}}$ associated to $c_{n}$, also has $\theta_n'$ as external argument.
This implies that $\alpha_1$ is fixed ($2\theta_n=\theta_{n}'$) and thus the external arguments of rays landing at $\alpha_1$ are exactly $\frac{1}{2^n-1},\dots,\frac{2^{n-1}}{2^n-1}$.

Let $U_1$ be the Fatou component of $c_1$ and $U_0$ that of the critical point $0$.
For a component $U$ of $\mathring{K_c}$, its image $f(U)$ is also some component of $\mathring{K_c}$ (this follows from the fact that $f_c$ is an open map) and, if $U$ does not contain $0$ (i.e., $U\neq U_0$), then $f_c$ restricts to a homeomorphism from $U$ to $f_c(U)$.

Since $\alpha_1$ is fixed and $f_c(U_0)=U_1$, we have that $\alpha_1\in \overline{U_0} \cap \overline{U_1}$.
The point $\alpha_1$ and its preimages thus belong to multiple circles.
The set collecting such points is dense, so point 1 of \cref{def.rabbit} is proved.

The projection to any component is a cut point.
Thus, if $x\in\partial U$ and $y\in\partial V$, where $U$ and $V$ are distinct components of $\mathring{K_{c_n}}$, then the projection of $x$ to $V$ is a cut point that separates $x$ from $y$ \cite[Exposé I, Proposition 4]{douady-hubbard-I}.
This proves point 2 of \cref{def.rabbit}

For $x\in K_{c_n}$, we denote by $\ell(x)$ the number of components $V$ of $\mathring{K_{c_n}}$ such that $V\cap[x,\alpha_1]\neq\emptyset$, where $[x,\alpha_1]$ is the legal arc between $x$ and $\alpha_1$ (legal arcs were introduced in \cite[Exposé II, \S6]{douady-hubbard-I} as \textit{arcs réglementaires}).
By construction of legal arcs, this number is constant on each component $U$ of $\mathring{K_c}$, so we denote by $\ell(U)$ this common number (which, a priori, could be infinite).
We claim that, for any component $U$, $|\ell(f_c(U))-\ell(U)|\leq 1$.
If the legal arc between $x\in U$ and $\alpha_1$ does not contain $0$ then there is equality, since $f_c$ maps components to components and, by \cite[Exposé IV, Lemme 1]{douady-hubbard-I}, it maps legal arcs to legal arcs injectively if the legal arc does not contain the critical point in its interior.
In case $[x,\alpha_1]$ contains $0$, we have that $[x,\alpha_1]=[x,0]\cup[0,\alpha_1]$.
If $[f_c(x),f_c(0)]=[f_c(x),c]$ contains $\alpha_1$, then $\ell(f_c(U))=\ell(U)-1$ because $[f_c(x),c]=[f_c(x),\alpha_1]\cup[\alpha_1,c]$ and $c$ belongs to the component $U_1$ which contains $\alpha_1$ in its boundary.
If instead $[f_c(x),c]$ does not contain $\alpha_1$, then $\ell(f_c(U))=\ell(U)$.

We know that $U_0$ contains a point of the attractive cycle, so for any component $U$, there is $n$ such that $f_c^n(U)=U_0$.
In particular, $\ell(U)\leq n+1$.
So, if $x$ lies on some circle of $J_c$, i.e., $x\in\partial U$ for some component $U$ of $\mathring{K_c}$, then $\ell(x)=\ell(U)$ or $\ell(x)=\ell(U)-1$, depending on whether $x$ is the projection of $\alpha_1$ on $U$ or not.
If $x$ and $y$ belong to circles, then the collection of circles that separate $x$ from $y$ is included in the collection of circles that separate $x$ from $\alpha_1$ and $y$ from $\alpha_1$. This is easily seen using the fact that the legal hull of $x$, $y$ and $\alpha_1$ (introduced in \cite[Exposé II, \S7]{douady-hubbard-I} as \textit{enveloppe réglementaire}) is a topological tree.
This collection is finite, since $\ell(x),\ell(y)<\infty$.
This proves point 3 of \cref{def.rabbit}.

By \cref{lem.cut.points.order}, any cut point of $J_{c_n}$ belongs to at least two circles.
Any point that belongs to at least two circles is eventually mapped to $\alpha_1$ (indeed, any component is eventually mapped to the cycle $U_0, f_c(U_0), \dots, f_c^n(U_0)$ where the $U_i$'s are the images of $U_0$ by $f_{c_n}$).
Finally, any cut point has the same order as $\alpha_1$, which is $n$, thus $J_{c_n}$ is an $n$-regular rabbit. 
\end{proof}

\subsection{The airplane Julia set}
 
\begin{proposition}
\label{sub.airplane.julia.set}
Let $\theta=3/7$ and $c$ be the point in the Mandelbrot set with external argument $\theta$.
The Julia set $J_c$ is an airplane.
\end{proposition}

\begin{remark}
As in \cref{rmk.basilica.rabbit}, this Julia set is known as the airplane.
\end{remark}

\begin{proof}
The point $c$ is the root of a hyperbolic component.
The external arguments $1/7$ and $2/7$ land at the same point of the Mandelbrot set, as do $5/7$ and $6/7$ by symmetry, so there is only one possibility: $3/7$ must land at the same point of the Mandelbrot set as $4/7$.

Let $\alpha_1$ be the point in $J_c$ with the external rays $3/7$ and $4/7$.
This is the projection of $0$ on the component $U_1$ of $\mathring{K_c}$ that contains $c$.
The point $\alpha_1$ belongs to the circle $\partial U_1$ and is a cut point, since it separates $0$ and $c$.
The set of its preimages is dense, which proves point 1 of \cref{def.airplane}.

The point $\alpha_1$ is mapped to the point with external arguments $\pm1/7$ which is mapped to the point with external arguments $\pm2/7$.
This last point is mapped to $\alpha_1$, which is thus 3-periodic.
Let $U_0$ be the component that contains $0$ and let $U_2$ be the image of $U_1$.
These components form a 3-cycle and $U_0$ lies between $U_1$ and $U_2$.
For two components $U,V$ of $\mathring{K_c}$, we denote by $\delta(U,V)$ the number of components that the legal arc from any point in $U$ to any point in $V$ crosses minus 1 (this number can be infinite).
Since $U_0$ is between $U_1$ and $U_2$, we have $\delta(U_1,U_2)=\delta(U_1,U_0)+\delta(U_0,U_2)$. Because of \cite[Exposé IV, Lemme 1]{douady-hubbard-I}, we have $\delta(U_2,U_1)=\delta(f_c(U_1),f_c(U_0))=\delta(U_1,U_0)$ and $\delta(U_0,U_1)=\delta(f_c(U_2),f_c(U_0))=\delta(U_2,U_0)$.
Since all these numbers are at least $1$, they must all be infinite.

Assume that two circles $C_1,C_2$ have a common point.
Let $V_1,V_2$ be the two components of $\mathring{K_c}$ such that $C_i=\partial V_i$.
There is a large enough $k\in\mathbf{N}$ such that $f_c^k$ maps $V_i$ to $U_{j_i}$ (for $i=1,2$ and $j_i \in \{0,1,2\}$) . By \cite[Exposé IV, Lemme 1]{douady-hubbard-I}, $f_c^k$ maps the legal arc between $V_1$ and $V_2$ to the legal arc between $U_{j_1}$ and $U_{j_2}$, so it maps the point shared by $V_1$ and $V_2$ to a point shared by two of the $U_i$'s.
We have a contradiction, since the $U_i$'s have at least one circle that separates them (in fact, infinitely many).
This proves point 3 of \cref{def.airplane}.

Assume that two circles $C_1,C_2$ are not separated by a third one.
Let $V_i$'s be components as above.
The legal arc $A$ joining the centers does not contain $0$ and the same is true for all $f_c^n(A)$.
Let  $n$ be a large enough number so that $f^n_c(V_i)$'s are some $U_j$'s.
The images of $V_1$ and $V_2$ are distinct and thus separated by infinitely many circles.
This implies that $A$ meets infinitely many circles.
A contradiction, which proves point 2 of \cref{def.airplane}.

It remains to prove that every cut point has order $2$.
Since every cut point has order at least $2$, it suffices to prove there is no cut point of order at least $3$.
By contradiction, let $x$ be such a point.
Since $K_c$ is obtained by filling circles from $J_c$, the complement $K_c\setminus\{x\}$ has at least $3$ connected components.
Choose components $V_1,V_2,V_3$ of $\mathring{K_c}$ lying in distinct components of $K_c\setminus\{x\}$.
Let $T$ the legal hull of the centers of the $V_i$'s.
Up to choosing a smaller tree, we may assume that $T$ does not contain $0$ in the interior of any of its arcs.
It is a tripod with $3$ arcs and $x$ as center.
We define a sequence of legal trees in the following way.
Let $T^{(0)}=T$ and assume that $T^{(n)}$ has been constructed and does not contain $0$ in the interior of any of its arcs.
Then, by \cite[Exposé IV, Lemme 1]{douady-hubbard-I}, $f_c$ maps $T^{(n)}$ injectively to $f_c(T^{(n)})$.
If the latter does not contain $0$ in the interior of its arcs, we set $T^{(n+1)} = f_c(T^{(n)})$.
If instead $0$ lies in the interior of one of the arcs of $f_c(T^{(n)})$, we replace this arc by $[f_c^n(x),0]$ and we let $T^{(n+1)}$ be the result of this replacement.
By construction, at any step $n$, $T^{(n)}$ is a tripod with 3 arcs and $f_c^n(x)$ as center.
Since each $V_i$ is eventually mapped to the attractive cycle $U_0,U_1,U_2$, for $n$ large enough the ends of $T^{(n)}$ eventually all fall in some $U_i$.
Since the $U_0$ is between $U_1$ and $U_2$, we have a contradiction with the fact that $T^{(n)}$ is a tripod.
This proves point 4 of \cref{def.airplane}.
\end{proof}

\subsection{Extending homeomorphisms to the complex plane}
\label{sub.extending.homeomorphisms}

Since Julia sets lie by definition in the complex plane, it is natural to determine which homeomorphisms of a Julia extend to a global homeomorphism of the complex plane.
Let us show that such homeomorphisms correspond to the elements of $\Homeo^+(\R_n)$ and $\Homeo^+(\A)$ from \cref{sec.lamination.automorphisms}.

We denote by $\Homeo^+_J(\C)$ the group of orientation preserving homeomorphisms of $\C$ stabilizing $J$ setwise, i.e. $g(J)=J$.

\begin{theorem}
Let $J$ be the airplane Julia set or a rabbit Julia set. The subgroup of $\Homeo(J)$ of homeomorphisms that are the restrictions of elements of $\Homeo^+_J(\C)$ coincide with $\Homeo^+(J)$.
\end{theorem}

\begin{proof} Let $g\in \Homeo(J)$ that is the restriction to $J$ of some $\overline{g}\in\Homeo^+_J(\C)$. Since the Riemann sphere $\hat{\C}$ is the Alexandroff compactification of $\C$, $\overline{g}$ extends uniquely as a homeomorphisme of $\hat{\C}$ fixing $\infty$. Take three points $x,y,z$ on some circle $C$ of $J$ and arcs $I_x,I_y,I_z$ such that $\infty$ is one end and $x$, respectively $y$ or $z$, is the other end, and such that $I_x\cap J=\{x\}$, respectively $\{y\}$ or $\{z\}$.
We may moreover assume that these arcs intersect only at $\infty$. For example, external rays do the job.

Since $\overline{g}$ is a homeomorphism of $\hat{\C}$, the arcs $\overline{g}(I_x),\overline{g}(I_y)$ and $\overline{g}(I_z)$ satisfy the same conditions and thus $O^C_\varphi(x,y,z)=O^{g(C)}_\varphi(g(x),g(y),g(z))$. So $g\in \Homeo^+(J)$.

Conversely, consider $g\in\Homeo^+(J)$.
Then there is a $\gamma\in\Homeo^+(\mathbf{S}^1)$ preserving the lamination $\mathcal{L}J$ such that $g\circ\varphi=\varphi\circ\gamma$. Let $\mathbf{D}$ be the closed unit disk in $\C$ and $\Phi\colon\hat{\C}\setminus\mathbf{D}\to\hat{\C}\setminus FJ$ be the Böttcher coordinates of the complement of the filled in Julia set $FJ$.
In particular, for any $z$ in the unit circle $\mathbf{S}^1$ of $\C$, $\varphi(z)=\lim_{r\to1}\Phi(rz)$ and this defines a continuous extension of $\Phi$ to the complement of the open unit disk.
We define $\overline{g}$ on $\hat{\C}\setminus FJ$ ($r>1$ and $z\in\mathbf{S}^1$) via the formula 
$$\overline{g}(\Phi(rz))=\Phi(r\gamma(z)).$$

This defines a homeomorphism of $(\hat{\C}\setminus FJ)\cup J$. Any component $U$ of the interior of $FJ$ is homeomorphic to an open disk with boundary $\partial U$ which is a circle of $J$.  For such $U$ we choose a homeomorphism between $U$ and the component $V$ of the interior of $FJ$ bounded by $g(\partial U)$ which extends to a homeomorphism between $\overline{U}$ and $\overline{V}$ that coincides with $g$ on $\partial U$. Patching together all these homeomorphisms of interior components we get a global homeomorphism $\overline{g}$ of $\hat{\C}$ whose restriction to $J$ is $g$.
\end{proof}

\begin{remark} The construction of the extension done in the above proof is not unique at all, but one can make it canonical. For each component $U$ of the filled-in Julia center, one can choose a center $c_U$ and this defines internal rays and arguments. Doing a similar construction to the one we do for external rays, one can define uniquely $\overline{g}$ on each component $U$ such that it maps internal rays of $U$ to internal rays of $V$.
This way, one can see that $\Homeo^+_J(\C)$ splits as permutational wreath product
$$\left(\Homeo_0(\C\setminus FJ)\times\prod_{U\in\pi_0(\interior{FJ})}\Homeo_0(U)\right)\rtimes \Homeo^+(J),$$
where $\Homeo_0(U)$ is the subgroup of $\Homeo(U)$ that admits a continuous extension to $\overline{U}$ that is the identity on $\partial U$.
\end{remark}


\section{Rabbits and airplane as limit spaces of replacement systems}
\label{sec.limit.spaces}

In this section we show that the basilica, rabbits and airplane limit spaces built with the replacement systems introduced in \cite{Belk-Forrest} satisfy our \cref{def.rabbit} and \cref{def.airplane}  thus are homeomorphic to the basilica, rabbits and airplane Julia sets.
In order for this section to be as concise as possible, the introduction to replacement systems and limit spaces is reduced to its essential parts and the proofs do not delve into details.
We refer to the paper which introduced this topic \cite{Belk-Forrest} and the second author's PhD thesis \cite{thesis} for a much more thorough introduction to replacement systems, limit spaces and (although we will not use them) rearrangement groups.

\subsection{Limit spaces of replacement systems}

The essential idea is that a \textit{replacement system} defines a rewriting system of graphs whose graphs, under mild conditions, expand and converge towards a \textit{limit space}, which is the quotient space of a Cantor set under an equivalence relation.
This is all entirely codified by graphs.
Let us briefly describe how.

\subsubsection{Replacement systems}

In this section, by \textbf{graph} we mean a quadruple $(V,E,\iota,\tau)$, where $V$ and $E$ are finite sets (of vertices and edges) and $\iota$ and $\tau$ are maps $E \to V$ that associate to each edge its initial and terminal vertices, respectively.

\begin{definition}[Definition 1.4 \cite{Belk-Forrest}]
A \textbf{replacement system} consists of the following collection of graphs, all edge-colored by a finite set $C$ of colors:
\begin{itemize}
    \item a \textbf{base graph};
    \item for each color $c \in C$, a \textbf{replacement graphs} $R_c$ equipped with distinct vertices $\iota_c$ and $\tau_c$, which we call \textbf{initial} and \textbf{terminal vertices}.
\end{itemize}
\end{definition}

Two examples of replacement systems are depicted in \cref{fig.basilica.replacement.system,,fig.airplane.replacement.system}.

\begin{figure}
\centering
\begin{tikzpicture}[font=\small]
    \useasboundingbox (-1.4,-.7) rectangle (6.5,1.35);
    \node at (0,1.1) {Base graph};
    \node[vertex] (s) at (0,0) {};
    \draw[edge] (s) to[loop,out=135,in=225,min distance=2cm,looseness=10] node[above left]{$L$} (s);
    \draw[edge] (s) to[loop,out=-45,in=45,min distance=2cm,looseness=10] node[above right]{$R$} (s);
    \begin{scope}[xshift=5cm,yshift=-.5cm]
    \node at (0,1.6) {Replacement graph};
    \node[vertex] (l) at (-1.3,0) {};
    \draw (l) node[above]{$\iota$};
    \node[vertex] (c) at (0,0) {};
    \node[vertex] (r) at (1.3,0) {};
    \draw (r) node[above]{$\tau$};
    \draw[edge] (l) to node[above]{$0$} (c);
    \draw[edge] (c) to[loop,out=130,in=50,min distance=1.3cm,looseness=10] node[above]{$1$} (c);
    \draw[edge] (c) to node[above]{$2$} (r);
    \end{scope}
\end{tikzpicture}
\caption{The basilica replacement system.}
\label{fig.basilica.replacement.system}
\end{figure}
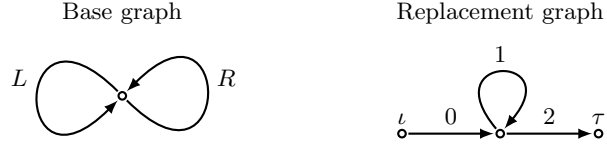

\begin{figure}
\centering
\begin{tikzpicture}[font=\small]
    \node at (0,1.4) {Base graph};
    \node[vertex] (l) at (-.667,0) {};
    \node[vertex] (r) at (.667,0) {};
    \draw[edge,blue] (l) to node[above]{$s$} (r);
    \begin{scope}[xshift=3.75cm]
    \node at (0,1.4) {\textcolor{blue}{Blue} replacement graph};
    \draw[edge,red,domain=5:175] plot ({.5*cos(\x)}, {.5*sin(\x)});
    \draw (90:.5) node[above,red] {$b_2$};
    \draw[edge,red,domain=185:355] plot ({.5*cos(\x)}, {.5*sin(\x)});
    \draw (270:.5) node[below,red] {$b_3$};
    \node[vertex] (l) at (-1.75,0) {}; \draw (-1.75,0) node[above]{$\iota_{\text{\textcolor{blue}{b}}}$};
    \node[vertex] (cl) at (-.5,0) {};
    \node[vertex] (cr) at (.5,0) {};
    \node[vertex] (r) at (1.75,0) {}; \draw (1.75,0) node[above]{$\tau_{\text{\textcolor{blue}{b}}}$};
    \draw[edge,blue] (cl) to node[above]{$b_1$} (l);
    \draw[edge,blue] (cr) to node[above]{$b_4$} (r);
    \end{scope}
    \begin{scope}[xshift=8cm]
    \node at (0,1.4) {\textcolor{red}{Red} replacement graph};
    \node[vertex] (l) at (-1.25,0) {}; \draw (-1.25,0) node[above]{$\iota_{\text{\textcolor{red}{r}}}$};
    \node[vertex] (r) at (1.25,0) {}; \draw (1.25,0) node[above]{$\tau_{\text{\textcolor{red}{r}}}$};
    \node[vertex] (c) at (0,0) {};
    \node[vertex] (ct) at (0,1) {};
    \draw[edge,red] (l) to node[below]{$r_1$} (c);
    \draw[edge,red] (c) to node[below]{$r_2$} (r);
    \draw[edge,blue] (c) to node[left]{$r_3$} (ct);
    \end{scope}
\end{tikzpicture}
\caption{The airplane replacement system.}
\label{fig.airplane.replacement.system}
\end{figure}
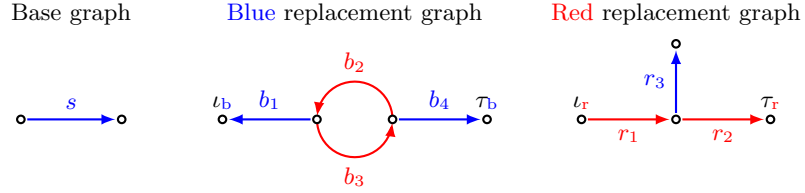

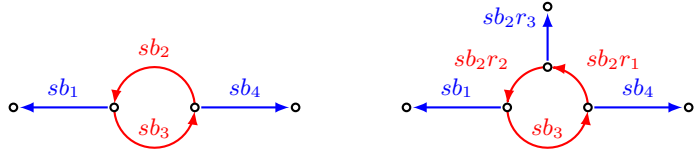
\begin{figure}
\centering
\begin{tikzpicture}[font=\small,scale=.8]
    \draw[edge,red,domain=5:175] plot ({.667*cos(\x)}, {.667*sin(\x)});
    \draw (90:.667) node[above,red] {$sb_2$};
    \draw[edge,red,domain=185:355] plot ({.667*cos(\x)}, {.667*sin(\x)});
    \draw (270:.667) node[above,red] {$sb_3$};
    \node[vertex] (l) at (-2.333,0) {};
    \node[vertex] (cl) at (-.667,0) {};
    \node[vertex] (cr) at (.667,0) {};
    \node[vertex] (r) at (2.333,0) {};
    \draw[edge,blue] (cl) to node[above]{$s b_1$} (l);
    \draw[edge,blue] (cr) to node[above]{$s b_4$} (r);
    \begin{scope}[xshift=6.5cm]
    \draw[edge,red,domain=95:175] plot ({.667*cos(\x)}, {.667*sin(\x)});
    \draw (135:.667) node[above left,red] {$sb_2r_2$};
    \draw[edge,red,domain=185:355] plot ({.667*cos(\x)}, {.667*sin(\x)});
    \draw (270:.667) node[above,red] {$sb_3$};
    \draw[edge,red,domain=5:85] plot ({.667*cos(\x)}, {.667*sin(\x)});
    \draw (45:.667) node[above right,red] {$sb_2r_1$};
    \node[vertex] (l) at (-2.333,0) {};
    \node[vertex] (cl) at (-.667,0) {};
    \node[vertex] (c) at (0,.667) {};
    \node[vertex] (ct) at (0,1.667) {};
    \node[vertex] (cr) at (.667,0) {};
    \node[vertex] (r) at (2.333,0) {};
    \draw[edge,blue] (cl) to node[above]{$s b_1$} (l);
    \draw[edge,blue] (cr) to node[above]{$s b_4$} (r);
    \draw[edge,blue] (c) to node[above left]{$s b_2 r_3$} (ct);
    \end{scope}
\end{tikzpicture}
\caption{Two graph expansions of the airplane replacement system.}
\label{fig.airplane.expansions}
\end{figure}

Let us fix a replacement system.
A $c$-colored edge $e$ of a graph can be \textbf{expanded} by replacing $e$ with the replacement graph $R_c$, identifying $\iota(e)$ and $\tau(e)$ with $\iota_c$ and $\tau_c$, respectively.
A \textbf{graph expansion} is a graph that can be obtained from the base graph by applying a finite sequence of expansions.
For example, \cref{fig.airplane.expansions} depicts two graph expansions of the airplane replacement system.

\subsubsection{The symbol space}

Let $e_1, \dots, e_k$ be all the edges of a replacement graph $R_c$.
When expanding a $c$-colored edge $f$, we denote the newly added edges by $f e_i$ (see for example \cref{fig.airplane.expansions}).
Thus, every edge of a graph expansion corresponds to some finite word in the alphabet of the edges of the base and replacement graphs.

\begin{definition}[Definitions 1.6 and 1.7 \cite{Belk-Forrest}]
Fix a replacement system.
\begin{itemize}
    \item Its \textbf{symbol space} $\mathfrak{C}$ is the set of infinite sequences $\alpha = x_1 x_2 \dots$ such that each finite prefix $x_1 \dots x_k$ is an edge of some graph expansion.
    \item Its \textbf{gluing relation} $\sim$ is the binary relation on $\mathfrak{C}$ defined by setting $x_1 x_2 \cdots \sim y_1 y_2 \cdots$ when, for all $k$ large enough, $x_1 \dots x_k$ and $y_1 \dots y_k$ are incident edges in all graph expansions that feature them both.
\end{itemize}
We equip $\mathfrak{C}$ with the subspace topology induced from the Cantor space $E^\infty$ (which is the set of all infinite sequences in the alphabet $E$), where $E$ is the set of edges of the base and replacement graphs.
\end{definition}

Under mild conditions, a replacement system is said to be \textit{expanding} (see \cite[Definition 1.8]{Belk-Forrest}).
For expanding replacement systems, $\mathfrak{C}$ is a Cantor space and the gluing relation is an equivalence relation by \cite[Proposition 1.9]{Belk-Forrest}.
It is straightforward to see that the basilica and airplane replacement systems are expanding replacement systems, so we can consider the quotient of $\mathfrak{C}/\sim$.

\begin{definition}[Definition 1.7 \cite{Belk-Forrest}]
The \textbf{limit space} of a replacement system is the quotient space $\mathfrak{C}/\sim$.
\end{definition}

Edges of graph expansions correspond to sets that tessellate the limit space:

\begin{definition}[Definition 1.13 \cite{Belk-Forrest}]

To each edge $e$ (which corresponds to some finite word) we associate a \textbf{cell}, which is the set of those points of $\mathfrak{C}/\sim$ some of whose representative begins with $e$.
\end{definition}

Although we do not discuss this here, it is worth recalling that replacement systems allow to easily define countable groups of homeomorphisms of limit spaces, the so-called \textit{rearrangement groups}.
Such groups include Thompson groups $F$, $T$ and $V$ together with new groups acting on fractals such as the basilica \cite{BF15}, the airplane \cite{airplane} and Wa\.zewski dendrites \cite{dendriterearrangement}.

\subsubsection{Topology of the limit space}

Two elements of $\mathfrak{C}$ are equivalent under $\sim$ when their prefixes approach a common vertex.
Essentially, vertices are where ``some topology happens.''

\begin{definition}
A point of a limit space is a \textbf{gluing vertex} if it corresponds to a vertex of some graph expansion and it is a \textbf{regular point} otherwise.
\end{definition}

Note that each regular point always has a unique representative in $\mathfrak{C}$.

\begin{proposition}
\label{prop.limit.spaces.facts}
Each limit space $X$ satisfies the following statements.
\begin{enumerate}
    \item $X$ is compact and metrizable \cite[Theorem 1.25]{Belk-Forrest}.
    \item If the base and replacement graphs are connected, then $X$ is connected \cite[Remark 1.27]{Belk-Forrest} and locally connected \cite[Corollary 2.33]{thesis}.
    \item A subset is dense if and only if it intersects every cell \cite[Remark 3.2]{IG}.
\end{enumerate}
\end{proposition}

\subsection{Arcs and circles in limit spaces}

Here we briefly develop some tools that allow to link the shape of graph expansions to the topology of the basilica and airplane limit spaces.

\begin{definition}
In a directed graph, an \textbf{undirected walk} is a sequence $(e_1, \dots, e_k)$ of edges such that $e_i$ is adjacent to $e_{i+1}$ for all $i\in\{1, \dots, k-1\}$.
Each undirected walk identifies a sequence of vertices $(v_0, v_1, \dots, v_k)$ that the walk travels through.
An undirected walk is an \textbf{undirected path} if there is no $i \neq j$ such that $v_i = v_j$ and it is an \textbf{undirected cycle} if $v_i = v_j$ only for $i=0$ and $j=k$.
\end{definition}

\begin{proposition}
\label{prop.arcs.paths}
Let $p$ and $q$ be gluing vertices of a limit space.
Every arc between $p$ and $q$ (if $p \neq q$, otherwise every circle that contains $p=q$) corresponds to a subset of the limit space obtained as follows.
\begin{enumerate}
    \item Choose a graph expansion where $p$ and $q$ appear as vertices and choose an unoriented path $(e_1, \dots, e_k)$ joining them (an unoriented cycle if $p=q$).
    \item For each edge $e_i$, choose an unoriented path $(e_{i1}, \dots, e_{ik_{i}})$ in its replacement graph joining the initial and terminal vertices.
    \item Iterate the previous point for each edge $e_{ij}$.
\end{enumerate}
\end{proposition}

Note that, if $p$ and $q$ are distinct and not joined by an arc, then they lie in distinct connected components.
In this case, after finitely many (possibly zero) iterations of this procedure the vertices corresponding to $p$ and $q$ will not be able to be joined by an unoriented path.

\begin{proof}
First note that, for each graph expansion featuring $p$ and $q$ as vertices, an arc joining $p$ and $q$ in the limit space is included in the union of cells corresponding to a unique unoriented path between $p$ and $q$ (this is because non-empty intersections of cells consist of one or two gluing vertices, so intervals can only travel between cells corresponding to adjacent edges by going through such gluing vertices).
By Proposition 1.24 of \cite{Belk-Forrest} a cell associated to an edge $e_1$ is included in one associated to an edge $e_2$ if and only if $e_2$ is a prefix of $e_1$ (as finite words).
Thus, given a graph $\Gamma_1$ that features $p$ and $q$ as vertices and a graph $\Gamma_2$ that is an expansion of $\Gamma_1$, if an arc joining $p$ and $q$ is included in the union of cells of unoriented paths $P_1$ in $\Gamma_1$ and $P_2$ in $\Gamma_2$, then $P_2$ must be obtained as an expansion of edges of $P_1$.
Thus, every arc is as in the statement.

For the converse, the idea is that the gluing relation in a sequence of ``nested'' paths as in the statement behaves like a generalization of $n$-ary expansions of numbers in $[0,1]$ with no fixed arity (see \cite{Belk-Forrest}, where Proposition 2.6 explains how to obtain the interval as a limit space with dyadic numbers as gluing vertices and Example 2.14 features examples whose limit spaces are intervals and the gluing vertices have no fixed arity).
Consider a sequence $P_1, P_2, \dots$ of paths as in the statement.
Let $Y$ be the set of classes of infinite words whose every long enough prefix, as an edge, lies in some $P_i$.
By \cite[Proposition 1.22]{Belk-Forrest}, two elements of the symbol space are equivalent under the gluing relation if and only if they represent the same gluing vertex.
Thus, for any $x = x_1 x_2 \dots$ and $y = y_1 y_2 \dots$ in the symbol space whose all long enough prefixes lie in some $P_i$, one has that $x \sim y$ if and only if all the long enough prefixes $x_1 \dots x_k$ and $y_1 \dots y_k$ are adjacent edges.
Since these prefixes all lie in the paths $P_i$'s, the gluing relation on $Y$ is non-trivial precisely on any vertex on the $P_i$'s, where it glues the two sequences of edges approaching from the two sides (compare with the binary gluing relation $w0111\ldots \sim w1000\ldots$ described in \cite[Proposition 2.6]{Belk-Forrest}).
Hence $Y$ is this homeomorphic to $[0,1]$ if $p \neq q$ and to $\mathbf{S}^1$ if $p=q$.
\end{proof}

In many tame yet interesting cases, the previous procedure will stabilize after finitely many steps and will thus identify all of the finitely many arcs joining the two points.
For example, in the basilica replacement system (\cref{fig.basilica.replacement.system}) there is a sole path joining the initial and terminal vertices of each replacement graph, so each unoriented path $(e_1, \dots, e_k)$ between two distinct vertices of a graph expansion identifies a unique arc between the two vertices.
The same holds for any rabbit limit space (see \cref{fig.n.rabbit.replacement.system}).

Instead, in the airplane replacement system (\cref{fig.airplane.replacement.system}) there are infinitely many arcs joining for example the gluing vertices corresponding to the initial and the terminal vertices of a blue edge, since each blue edge requires the choice of one of two red edges and there are infinitely many blue edges at which this choice needs to be made.
However, if two distinct vertices belong to a red cycle of some graph expansion, then by \cref{prop.arcs.paths} there are precisely two arcs joining them and there is a unique circle that corresponds to that cycle.

\begin{remark}
\label{rmk.unique.circles}
A circle of a limit space must contain some gluing vertex (else it would be totally disconnected).
Thus, each circle in the basilica or in the airplane limit space corresponds to a cycle in a graph expansion.
\end{remark}

\begin{figure}
\centering
\begin{tikzpicture}
    \node at (0,1.4) {Base graph};
    \node[vertex] (t) at (0,.8) {};
    \node[vertex] (b) at (0,-.8) {};
    \draw[edge] (b) to[out=180,in=180,looseness=2] node[left]{$l$} (t);
    \draw[edge] (b) to[out=90,in=270] node[left]{$c$} (t);
    \draw[edge] (b) to[out=0,in=0,looseness=2] node[left]{$r$} (t);
    \begin{scope}[xshift=5.5cm]
    \node at (0,1.4) {Replacement graph};
    \node[vertex] (l) at (-1.75,0) {}; \draw (-1.75,0) node[above]{$\iota$};
    \node[vertex] (cl) at (-.5,0) {};
    \node[vertex] (cr) at (.5,0) {};
    \node[vertex] (r) at (1.75,0) {}; \draw (1.75,0) node[above]{$\tau$};
    \draw[edge] (cl) to node[above]{1} (l);
    \draw[edge] (cr) to node[above]{4} (r);
    \draw[edge] (cr) to[out=90,in=90,looseness=1.4] node[above]{2} (cl);
    \draw[edge] (cl) to[out=270,in=270,looseness=1.4] node[above]{3} (cr);
    \end{scope}
\end{tikzpicture}
\caption{A replacement system for the Julia set of the rational map $z \mapsto z^{-2}-1$.}
\label{fig.bubble.bath.replacement.system}
\end{figure}
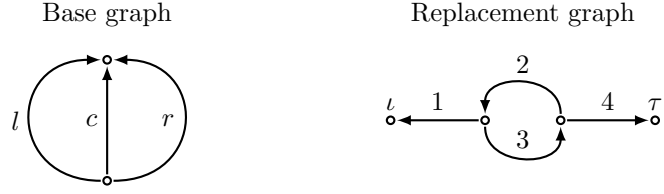

The previous remark does not hold unless the procedure described in \cref{prop.arcs.paths} stabilizes after finitely many steps.
For example, this does not work in the replacement system depicted in \cref{fig.bubble.bath.replacement.system}, whose limit space is homeomorphic to the Julia set for the complex rational map $z \mapsto z^{-2}-1$.

\subsection{The rabbit limit spaces}

\begin{figure}
\centering
\begin{tikzpicture}[font=\small]
    \useasboundingbox (-1.85,-.5) rectangle (6.75,1.85);
    \node at (0,1.65) {Base graph};
    \node[vertex] (s) at (0,0) {};
    \draw[edge] (s) to[loop,out=145,in=215,min distance=2cm,looseness=10] node[above left]{$X_1$} (s);
    \draw[dotted,gray] (140:.8) to[out=40,in=140] node[above]{$X_i$} (40:.8);
    \draw[edge] (s) to[loop,out=-35,in=35,min distance=2cm,looseness=9] node[above right]{$X_n$} (s);
    \begin{scope}[xshift=5cm]
    \node at (0,1.65) {Replacement graph};
    \node[vertex] (l) at (-1.3,0) {};
    \draw (l) node[above]{$\iota$};
    \node[vertex] (c) at (0,0) {};
    \node[vertex] (r) at (1.3,0) {};
    \draw (r) node[above]{$\tau$};
    \draw[edge] (l) to node[below]{$0$} (c);
    \draw[edge] (c) to[loop,out=165,in=125,min distance=1.3cm,looseness=9] node[above left]{$1$} (c);
    \draw[dotted,gray] (120:.8) to[out=30,in=150] node[above]{$i$} (60:.8);
    \draw[edge] (c) to[loop,out=55,in=15,min distance=1.3cm,looseness=9] node[above right]{$n-1$} (c);
    \draw[edge] (c) to node[below]{$n$} (r);
    \end{scope}
\end{tikzpicture}
\caption{The $n$-rabbit replacement system.}
\label{fig.n.rabbit.replacement.system}
\end{figure}
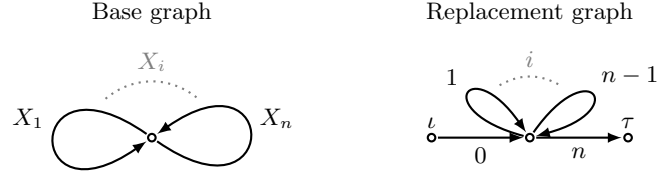

Let us now consider the \textbf{rabbit limit space} $X_n$, which is the limit space of the rabbit replacement system schematically depicted in \cref{fig.n.rabbit.replacement.system}:
its base graph is a bouquet of $n$ loops and its replacement graph is a path of length $2$ joining $\iota$ and $\tau$, with $n-1$ loops attached at the central vertex.

\begin{proposition}
The rabbit limit space $X_n$ is an $n$-regular rabbit in the sense of \cref{def.rabbit}.
\end{proposition}

\begin{proof}
By points (1) and (2) of \cref{prop.limit.spaces.facts}, each rabbit limit space is a Peano continuum.
Each regular point is either an end point or lies on a unique circle (this only depends on whether its unique representative has or does not have infinitely many digits in the set $\{1,2,\dots,n-1\}$), so gluing vertices are the only cut points of the basilica limit space.
Every vertex always separates a graph expansion in $n$ connected components, so the order of each cut point is $2$.
We now need to show that conditions 1, 2 and 3 of \cref{def.rabbit} hold.

1.
In every graph expansion of the rabbit replacement system, each edge belongs to some cycle and its expansion always features a loop, so each cell contains a point that belongs to multiple circles.
By point (3) of \cref{prop.limit.spaces.facts}, the set of points that belong to multiple circles is dense in the rabbit limit space.

2.
Two cycles of any graph expansion can always be separated by removing a vertex, so any two circles can be separated by a cut point.

3.
Consider two points $p$ and $q$ that belong to circles,
Suppose first that both of them are gluing vertices.
Then every unoriented path joining the two vertices travels through edges belonging to finitely many unoriented cycles, so each arc that joins $p$ and $q$ is included in the union of finitely many circles.
If $p$ is not a gluing vertex and belongs to a circle $C$, then for every gluing vertex $x$ on $C$ there are two arcs joining $p$ to $x$ that are entirely included in $C$, so the desired conclusion can also be drawn when $p$ and/or $q$ are not gluing vertices, as long as they belong to circles.
\end{proof}

\subsection{The airplane limit space}

Now consider the \textbf{airplane limit space}, which is the limit space of the replacement system depicted in \cref{fig.airplane.replacement.system}.

\begin{proposition}
The limit space of the airplane replacement system is an airplane in the sense of \cref{def.airplane}.
\end{proposition}

\begin{proof}
By points (1) and (2) of \cref{prop.limit.spaces.facts}, the airplane limit space is a Peano continuum.
We need to show that it satisfies conditions 1, 2, 3 and 4 of \cref{def.airplane}.

1.
In every graph expansion of the airplane replacement system, each red edge belongs to some unoriented cycle and its expansion always features a vertex between that cycle and a blue edge.
Each such vertex corresponds to a cut point of the limit space, so each red cell contains a cut point.
Since the expansion of a blue edge always features two red edges, each blue cell includes a red one and thus contains a cut point.
Thus, the set of all cut points that belong to circles is dense by point (3) of \cref{prop.limit.spaces.facts}.

2.
Every blue cell includes a circle and, in all graph expansions, any two distinct unoriented cycles are separated by some blue edge.
Hence, any two circles are separated by a third circle.

3.
Clearly, in all graph expansions, every two unoriented cycles are disjoint, so every two distinct circles of the limit space are disjoint.

4.
As noted above, those gluing vertices that belong to some circle are cut points.
The removal of the corresponding vertex from a graph expansion breaks the limit space in two connected components, so these cut points have order $2$.
Vertices that do not belong to unoriented cycles instead are end points, so they are not cut points.
Regular cut points of a limit space can only have order $1$ or $2$, so we are done.
\end{proof}

\end{appendices}


\begingroup
\emergencystretch=1em
\printbibliography[heading=bibintoc]
\endgroup

\end{document}